\font\bigbf=cmbx10  scaled \magstep1 

\def\1{\rlap 1\kern.4pt1}
\def\<{\langle}\def\>{\rangle}
\def\\{\qquad}
\def\={\equiv}

\def\and{\,\&\,}	  	\def\cite#1{{\rm[{\bf#1}]}}
	  
  \def\disjointsum{\rlap{$+$}\kern.15em\cdot}
\def\disjointunion{\rlap{$\union$}\kern1.8pt\raise.5em\hbox{$\cdot$}\kern2pt}

 		\def\elementary{\prec}
	  \def\false{\neg}	
	  \def\implies{\Longrightarrow}
\def\includedin{\subseteq}\def\intersect{\cap} 	
\def\iso{\simeq}	  
\def\norm#1{\Vert#1\Vert} \def\normof #1{\Vert#1\Vert}
\def\normal{\triangleleft} 
\def\normalne{\raise1pt\hbox{\rlap{$\normal$}}_\ne}
\def\notiso{{\rlap{$\iso$}\kern 2 pt /}}

\def\onto{\rlap{$\to$}\kern.05pt\to}
\def\orthosum{\rlap{\hbox{$\oplus$}}\kern .3 ex\raise 1.6 ex
   \hbox{$\scriptstyle\perp$}\kern.1ex}
\def\proves{\vdash}	  \def\restrictedto{|\!^{\scriptscriptstyle\backslash}}
\def\qed{\hfill{\vrule height 8pt depth .05pt width 5pt}}
\def\satisfies{\models}	  \def\sdprod{{>}{\kern -2pt \triangleleft}}
\def\splitline#1//#2//#3{\vskip .1 in{\noindent\raise 10 pt\hbox{$#3$}\vbox{\hbox{#1}\hbox{#2}}\par}\vskip 4 pt}
\def\sqr#1#2{{\vbox{\hrule height .#2pt\hbox{\vrule width .#2pt height #1pt
		\kern#1pt\vrule width .#2pt}\hrule height .#2pt}}}

\def\to{\,\longrightarrow\,}		
\def \union{\cup}  	  \def\up#1{\,{}^{#1}\!}

\def\with{\leftrightarrow}\def \union{\cup}

\def\And{\bigwedge}
\def\Disjointunion{\rlap{$\Union$}\kern 3pt\raise.6em\hbox{$\cdot$}\kern3pt}
\def\Intersection{\bigcap}	
\def\Union{\bigcup}

\def\vline{\kern 1 pt{\vrule depth 1pt height .84 em width .4pt}\kern.5pt}


\def\fakebf#1{\rlap{$#1$}\kern.25pt\rlap{$#1$}\kern.25pt#1}

\def\bb#1#2{\rlap{\rm #1}\kern #2pt{\rm #1}}
\def\bbrule{\vrule width .3 pt height .68em depth 0 pt}
\def\scriptbb#1{\rlap{$\rm \scriptstyle #1$}\kern.5pt\hbox{$\rm \scriptstyle#1$}}
\def\midline{\vrule width .3 pt height 5.5 pt depth -.4 pt}

\def\littleslant{\hbox{$\scriptscriptstyle/$}}
\def\slantline{\raise 1.5 pt \hbox{\rlap{\littleslant}\kern.1pt\raise.2pt\littleslant}}
\def\add#1#2#3{\leavevmode\hbox to .8em{\rlap{\kern #3 pt #1}{\rm #2}\hfil}}

\def\Aa{{\raise .39ex \hbox{$\scriptstyle/$}\kern-.9 ex {\rm A}}}
\def\Cc{\add\midline C{2}} 	\def\Ff{\add\bbrule F{2}}
\def\Hh{\add\bbrule H{2}} 	
\def\Nn{\add\bbrule N{0}} 	\def\Pp{\add\bbrule P{2}}
\def\Qq{\add\midline Q{2}}
\def\Rr{\add\bbrule R{3}}
\def\Ss{\leavevmode\rlap{S}\kern .4ex
\rlap{\vrule width .4 pt height .68em depth -.42em}\kern
.35ex\vrule width .4 pt height 3 pt depth -.1 pt\kern 3pt}
\def\Tt{\add\bbrule T{2.4}}

\def\Zz{\bb Z{1}}

\def\UseBbb{\def\Aa{{\Bbb A}}\def\Cc{{\Bbb C}}\def\Ff{{\Bbb F}}
   \def\Nn{{\Bbb N}} \def\Pp{{\Bbb P}}  \def\Qq{{\Bbb Q}} \def\Rr{{\Bbb R}}
   \def\Tt{{\Bbb T}}\def\Zz{{\Bbb Z}}}

\def\AA{{\cal A}}\def\BB{{\cal B}}
\def\FF{{\cal F}}
\def\LL{{\cal L}}
\def\MM{{\cal M}}\def\PP{{\cal P}}
\def\QQ{{\cal Q}}
\def\TT{{\cal T}}\def\VV{{\cal V}}


\def\dom{{\rm  dom\,}}
\def\len{{\rm  len\,}}\def\max{{\rm  max\,}}
\def\mod{{\rm  mod\ }}\def\rg{{\rm  range\,}}
\def\rk{{\rm  rk\,}}\def\sup{{\rm  sup\,}}

\def\Br{{\rm  Br \,}}
\def\min{{\rm min \,}}
\def\Max{{\rm Max\,}}
\def\cf{{\rm  cf\,}}
\def\ess{{\rm  ess\,}}

\def\br{{\rm  Br\,}}

\def\Rang{{\rm Rang}}
\def\sucs{{\rm succ}} 
\def\Dom{{\rm Dom\,}}
\def\Min{{\rm Min\,}}



\def\numbersattop{\footline={\hfil}\headline={\hfil\tenrm\folio}}

\def\doublespace{\baselineskip=24pt}

\def\page{\vfill\break}
\def\singlespace{\baselineskip=\normalbaselineskip}

\def\case#1.{\noindent{\bf Case #1.}} \def\claim#1\par#2\par{\noindent#1{\bf Claim.}\par#2\par}
\def\proof#1:{\noindent{\sl Proof#1\/}:\par}
\def\section#1.#2\par{\bigskip\noindent{\bf \S#1. #2}\par\medskip}
\def\shelahtitle{\vbox {\hbox{\bigbf\title}\vskip.6true in\hbox{Saharon Shelah}
			\medskip\hbox{\sl The Hebrew University}
			\hbox{\sl Rutgers University}
                        \hbox{\sl MSRI}}}

\def\shelahfirst{
                {\topskip= 2.5 true in \centerline\shelahtitle
		\vskip 1 true in \itemitem{{\sl Abstract\/}:}\abstract\par\vfill
		\noindent\thanks\par}\break}

\def\contents{\parindent=0pt\bigskip\noindent{\bf
Contents}\bigskip\leftskip=2em}
\def\endcontents{\parindent=20pt\leftskip=0pt}

\def\part#1.#2\par{\bigskip\noindent\S#1.\\{\sl #2\par}}
\def\proclaim#1{\bigbreak\medskip\noindent{\bf
#1}\par\nobreak\medskip\nobreak\sl} 
\def\endproclaim{\rm\nobreak\smallskip\nobreak}

\def\dom{{\rm  dom\,}}
\def\len{{\rm  len\,}}

\def\epsilon{\varepsilon}
\def\sss{\scriptscriptstyle}
\def\midbar{\raise .4em\hbox{\vrule depth 0pt height .4pt width 4pt}\kern1pt}

\def\before{<\kern-5pt^{\sss|}_{\sss|}}
\def\cof{{\rm cof}}
\def\disjointunion{\rlap{$\union$}\kern1.8pt\raise.5em\hbox{$\cdot$}\kern2pt}
\def\diamond{\diamondsuit}
\def\down#1{\,_#1\!\!}
\def\forces{\kern2pt{\vrule depth 1pt height .8 em width .4pt}\kern.5pt\proves}
\def\name#1{\hbox to 0pt{\parindent=0pt\vbox to 0pt{\vbox to-5pt{}%
\char'176\vss}\hss}#1}
\def\nname#1{\hbox to 0pt{\parindent=0pt\vbox to 0pt{\vbox to-8pt{}%
\char'176\vss}\hss}#1}
\def\supername#1{\hbox to 0pt{\parindent=0pt\vbox to 0pt{\vbox to-5.5pt{}%
\char'176\vss}\hss}#1}
\def\notforces{\rlap/\kern-2pt\forces}
\def\otp{{\rm otp\,}}
\def\reals{\up\omega\omega}
\def\up#1{\,{}^{#1}\kern -1.5pt}


\def\condition#1{{\textstyle \tt #1}}
\def\pp{\condition p}

\def\AA{{\cal A}}
\def\BB{{\cal B}}
\def\DD{{\cal D}}
\def\Hh{{\tt H}}
\def\Pp{{\textstyle \tt P}}
\def\shorterline{\rlap{\vrule width .6em height .4pt depth 0pt}}
\def\shortline{\rlap{\vrule width .75em height .4pt depth 0pt}}
\def\slidedown{\kern.07em\lower.4pt}
\def\slideup{\kern.07em\raise.4pt}
\def\bethtop{\raise.8em\hbox{\shorterline\slidedown\shorterline\slidedown\shorterline}}
\def\bethbottom{\lower.3em\hbox{\shortline\slideup\shortline\slideup\shortline}}

\def\beth{\hbox{\rlap{\bethtop}\rlap{\bethbottom}\kern.63em{]}}}

\doublespace
\tolerance=500
\magnification=\magstep1

\def\1{\rlap 1\kern.4pt1}

\def\upower#1{#1^\omega/\FF}
\def\uprod#1{\prod #1/\FF}

\def\AP{{\cal AP }}
\def\BB{{\cal B}}
\def\DD{{\cal D}}
\def\LT{{\cal LT}}
\def\QQ{{\cal Q}}

\def\spacing{\singlespace}

\spacing
\def\unA{{\AA}}
\def\unF{{\bf F}}

\def\unP{\PP}
\def\Nnorm{N}
\def\rsucc{{\rm succ}}
\font\msym=msym10
\textfont9=\msym
\mathchardef\Nsp="94E
\mathchardef\Q="951
\def\LTfd{\LL\TT^f_d}
\def\LTdf{\LL\TT^f_d}
\def\LTdfat#1{(\LTfd)^{[#1]}}
\def\LTfdat#1{(\LTfd)^{[#1]}}
\def\uprodname{\prod}

\def\PpA#1{$(\Pp\restrictedto \AA_{#1})$}
\def\can{\name {can}}
\def\ref#1.#2\par{\vskip .1 in\leavevmode\noindent \hangindent 2.8 em
\hangafter =1
\llap{\hbox to 2.8em{[#1]\hfil}} #2\par}
\def\same{--------}

\input mssymb

\def\bigdisjointunion{\mathop{\dot{{\bigcup}}}}

\def\restrictedto{{\restriction}}
\def\forces{\Vdash}
\def\uhr{\restrictedto}
\def\lk{\langle}
\def\rk{\rangle}
\def\nl{\par\noindent}
\def\conc{\hat{~~}}
\font\msxm=msxm10

\def\xxop{\mathrel{\hbox{\msxm\char"43}}}
\UseBbb


\def\title{Vive la Diff\'erence I: Nonisomorphism  of
 ultrapowers of countable models}
\def\abstract{We show that it is not provable in ZFC that any two countable
elementarily equivalent structures have isomorphic ultrapowers relative to
some ultrafilter on $\omega$.}

\def\thanks{The author thanks the BSF and the Basic Research Fund, Israeli  
Academy of Sciences and the NSF
for partial support of this
research.\par\medskip \noindent This paper owes its existence to
Annalisa Marcja's hospitality in Trento, July 1987; van den Dries'
curiosity about Kim's conjecture; the willingness of Hrushovski and
Cherlin to look at \S3 through a dark glass; and most of all to
Cherlin's insistence that this is one of the fundamental problems of
model theory.

\noindent Number 326 in the publication list.

}

\shelahfirst

\contents

\pageno=-1
\numbersattop
\doublespace
\part 1. Elementarily equivalent structures do not have isomorphic
ultrapowers.
\numbersattop

If $V$ is a model of CH then in a generic extension we make
$2^{\aleph_0}=\aleph_2$ and we find countable elementarily equivalent graphs
$\Gamma$, $\Delta$ such that for {\sl every} ultrafilter $\FF$ on $\omega$, $\upower
\Gamma\not\iso\upower \Delta$.
In this model there is an ultrafilter $\FF$ such that
any ultraproduct  with respect to
$\FF$ of finite structures is saturated.

\part 2. The case of finite graphs.

By a variant of the construction in \S1 we show that there is a generic extension of $V$ in which
for some explicitly defined sequences of finite graphs $\Gamma_n$, $\Delta_n$,
all nonprincipal ultraproducts $\prod_n {\Gamma_n}/\FF_1$ or
$\prod_n{\Delta_n}/\FF_2$,  are elementarily equivalent,  but no countable 
ultraproduct of the 
$\Gamma_n$ is isomorphic to a countable ultraproduct of the $\Delta_n$.

\part 3. The effect of $\aleph_3$ Cohen reals.

We prove that if we simply add $\aleph_3$ Cohen reals to a model of
GCH, then there is at least one ultrafilter $\FF$ such that for
certain pseudorandom finite graphs $\Gamma_n$, $\Delta_n$, the
ultraproducts $\prod_n {\Gamma_n}/\FF$, $\prod_n {\Delta_n}/\FF$ are
elementarily equivalent but not isomorphic.  This implies that there
are also countably infinite graphs $\Gamma$, $\Delta$ such that for
the same ultrafilter $\FF$, the ultrapowers $\upower \Gamma$, $\upower
\Delta$ are elementarily equivalent and not isomorphic.

\part A.  Appendix.

We discuss proper forcing, iteration theorems, and the use of
$(Dl)_{\aleph_2}$ in \S3.
\endcontents

\page

\spacing
\openup 2 pt

\beginsection 0. Introduction.

\\Any two elementarily equivalent structures of cardinality $\lambda$ have
isomorphic ultrapowers (by [Sh 13], in 1971) with respect to an
ultrafilter on 
$2^\lambda$.  
Earlier, as the culmination of work in the sixties, Keisler
showed, assuming $2^\lambda=\lambda^+$, that the ultrafilter may be
taken to be on $\lambda$ [Keisler].  In particular, assuming the
continuum hypothesis,
for countable structures any nonprincipal ultrafilter on $\omega$ will
do.  As a special case,  the continuum hypothesis implies that an ultraproduct
of power series rings over prime fields $F_p$ is isomorphic to the
ultrapower of the corresponding rings of $p$-adic integers ; this has
number-theoretic consequences [AxKo].  Kim has conjectured that the
isomorphism $\prod_p{F_p[[t]]}/\FF\iso\prod_p{\Zz_p}/\FF$ is valid for any
nonprincipal ultrafilter over $\omega$, regardless of the status of the
continuum hypothesis.  In fact it has not previously been clear what
could be said about isomorphism of nonprincipal ultrapowers or
ultraproducts over $\omega$ {\sl in general}, in the absence of the
continuum hypothesis; it has long been suspected that such questions
do involve set theoretic issues going beyond ZFC, but there have been
no concrete results in this area.   For the case of two different
ultrafilters and on higher cardinals, see [Sh a VI].
In particular, ([Sh a VI, 3.13]) if $M=(\omega, <)^\omega/D$
($D$ an ultrafilter on $\omega$), the cofinality of
$(\{a\in M:a>n$ for every natural number $n\}, >)$ can be any regular
$\kappa\in (\aleph_0, 2^{\aleph_0}]$.

It does follow from the results of [Sh 13] that
there is always an ultrafilter $\FF$ on $\lambda$ such that for any two
elementarily equivalent models $M,N$ of cardinality $\lambda$, $\upower
M$ embeds elementarily into $\upower N$.  
On the other hand,  we show here that it is easy to find some countable elementarily equivalent
structures with nonisomorphic ultrapowers relative to a certain
nonprincipal ultrafilter on $\omega$: given enough Cohen reals, some
ultrafilter will do the trick (\S3), and with more complicated forcing any
ultrafilter will do the trick (\S1, refined in \S2).  The (first order
theories of the) models involved have the
independence property but do not have the strict order property. Every
unstable theory either has the 
independence property or the strict order property (or both) (in nontechnical
terms, in the theory we can interprate in a way the theory 
of the random graph or
the theory of a linear order), and our work here clearly makes use of
the independence property. The rings occurring in the Ax-Kochen
isomorphism are unstable, but do not have the independence property,
so the results given here certainly do not apply directly to Kim's
problem.  However it does appear that the methods used in \S3 can be
modified to refute Kim's conjecture, and we intend to return to this
elsewhere [Sh 405].  

A final technical remark: the forcing notions used here are
${<}\omega_1$-proper, strongly proper, and Borel.  Because of
improvements made in the iteration theorems for proper forcing [Sh
177, Sh f], we just need the properness; in earlier versions
$\omega$-properness was somehow used.

In the appendix we give a full presentation of a less general variant
of the preservation theorem of [Sh f] VI \S1. 

The forcing notions introduced in \S1, \S2 here (see 1.15, 1.16) are
of interest per se.  Subsequently specific cases have found more
applications; see Bartoszynski, Judah and Shelah [BJSh 368], Shelah
and Fremlin [ShFr 406]. 

\page

\numbersattop
\pageno=1

\section 1.  All ultrafilters on $\omega$ can be inadequate

\bigskip

Starting with a model $V$ of CH, in a generic extension we will make
$2^{\aleph_0}=\aleph_2$ and find countable elementarily equivalent
graphs $\Gamma,\Delta$ such that for any pair of ultrafilters $\FF,\FF'$ on
$\omega$, $\upower \Gamma \not\iso \Delta^\omega/\FF'$.
   More precisely:

\proclaim{1.1\\Theorem}

Suppose $V\satisfies$CH.  Then there is a proper forcing notion $\PP$
with the $\aleph_2$-chain condition, of cardinality $\aleph_2$ (and
hence $\PP$ collapses no cardinal and changes no cofinality) which makes
   $2^{\aleph_0}=\aleph_2$ and has the following effects on ultraproducts:

\item{(i)} There are countable elementarily equivalent graphs
   $\Gamma, \Delta$
   such that no ultrapowers $\upower \Gamma_1$,
   $\upower \Delta_2$ are isomorphic.

   \sl\item{(ii)} There is a nonprincipal ultrafilter $\FF$ on $\omega$ such
   that for any two sequences $\Gamma_n,\Delta_n$ of finite models for a countable
   language, if their ultrapowers with respect to $\FF$ are elementarily
   equivalent, then these ultrapowers are isomorphic, and in fact saturated.

\endproclaim

\proclaim{1.2\\Remark}
   \endproclaim

  The two properties {\sl (i,ii)} are handled quite independently by
  the forcing, and in particular {\sl (ii)} can be obtained just by
  adding random reals.

\proclaim{1.3 \\Notation}

\endproclaim

We work with the language of bipartite graphs (with a specified
bipartition $P,Q$).  $\Gamma_{k,l}$ is a bipartite graph with
bipartition
$U=U_{k,\ell},V=V_{k,\ell}$, $|U|=k$ and $V=\bigdisjointunion_{m<l}{U \choose m}$, where
${U\choose m}$ denotes the set of all subsets of $U$ of cardinality
$m$.  The edge relation is membership.  We also let $\Gamma_\infty$ be
the bipartite graph with $|U|=\aleph_0$ specifically $U=\omega$
and $V$ the set of all finite
subsets of $U$.  The theory of the $\Gamma_{k,l}$ converges to that of
$\Gamma_\infty$ as $l,k/l\to\infty$.

\proclaim{1.4\\Remark}
\endproclaim

Our construction will ensure that for any sequence $(k_n,l_n)$ with
$l_n,k_n/l_n\to\infty$ and any ultrafilters $\FF_1,\FF_2$ the
ultraproducts $\prod_i{\Gamma_{k_{n_i},l_{n_i}}}/\FF_1$ and $\upower
{\Gamma_\infty}_2$ are nonisomorphic.  In particular, if $\Gamma_{\rm
fin}$ is the disjoint union of the graphs $\Gamma_{2^n,n}$, and
$\Gamma$ is the disjoint union of $\Gamma_{\rm fin}$ and
$\Gamma_\infty$, then $\Gamma_{\rm fin}$ and $\Gamma$ are elementarily
equivalent, but any isomorphism of $\upower
\Gamma$ and $\upower {\Gamma_{\rm fin}}$ would induce an isomorphism of
an ultrapower of $\Gamma_\infty$ with some ultraproduct
$\prod_i{\Gamma_{2^{n_i},n_i}}/\FF$.  (Note that the graphs under
consideration have connected components of diameter at most 4.)

\proclaim{1.5\\The model}
\endproclaim

We will build a model $N$ of ZFC by iterating suitable proper forcing
notions with countable support [Sh b], see also [Jech].  The model $N$
will have the following combinatorial properties:

\item{P1.}  If 
$(A_n)_{n<\omega}$ is a collection of finite sets with $|A_n|\to\infty$, and
$g:\omega\to\omega$ with $g(n)\to\infty$,
and $f_i\ (i<\omega_1)$ are functions from $\omega$ to
$\omega$ with $f_i\in \prod_n A_n$ for all $i<\omega_1$, then there is a function $H$ from
$\omega$ to finite subsets of $\omega$ such that:  $H(n)$ has
size at most $g(n)$;
$H(n)\includedin A_n$; and for each $i$, $H(n)$ contains
$f_i(n)$ if $n$ is sufficiently large (depending on $i$).

\item{P2.}  $^\omega\omega$ has true cofinality $\omega_1$, that is: there is a sequence
$(f_i)_{i<\omega_1}$ which is cofinal in $^\omega\omega$ with respect to the partial
ordering of eventual domination (given by ``$f(n)<g(n)$ for sufficiently large $n$'').

\item{P3.} For every sequence $(A_k:k<\omega)$ of finite sets, for  any
collection $B_i (i<\omega_1)$ of infinite subsets of $\omega$, and 
for any collection $(g_i)_{i<\omega_1}$ of functions in $\prod_k A_k$, there is a function 
$f\in \prod_k A_k$ such that for all $i,j<\omega_1$, the set $\{n\in
B_i:f(n)=g_j(n)\}$ is infinite.

\item{P4.} $2^{\aleph_1}=\aleph_2$.

Note that (P3,P4) imply $2^{\aleph_0}=\aleph_2$.

\proclaim{1.6 Proposition}

Any model $N$ of ZFC with properties (P1-P2) satisfies part (i) of Theorem 1.1.
More precisely,
 the following weak saturation property holds for any ultraproduct
$\Gamma^*=\prod_n{\Gamma_{k_n,l_n}}/\FF$ for which $l_n\to\infty$, $(\ell_n<k_n)$
and fails in any countably
indexed ultrapower of $\Gamma_\infty${\/\rm:}

$(\dagger)$\\ Given $\omega_1$ elements of $U^{\Gamma^*}$,
 some element of $V^{\Gamma^*}$ is
linked to each of them.

\endproclaim

\proof:

Our discussion in Remark 1.4 shows that it suffices to check the claim
regarding $(\dagger)$. First consider an ultraproduct
$\Gamma^*=\prod_n{\Gamma_{k_n,l_n}}/\FF$ for which $l_n\to\infty$,
$l_n<k_n$.

Given $\aleph_1$
elements $a_i=f_i/\FF\in\Gamma^*$ we apply (P1) with $g(n)=l_n-1,
A_n=U_{k_n,l_n}$.  $H$ picks out a sequence of small subsets of
$U_{k_n,l_n}$, and if $b\in V^{\Gamma^*}$ is chosen so that its $n$-th
coordinate is linked to all the elements of $H(n)$, then this does the
trick.

Now let $\Gamma^*$ be of the form $\upower {\Gamma_\infty}$.  We will
show that $(\dagger)$ fails in this model.  Let $(f_i:i<\omega_1)$ be
a cofinal increasing sequence in $^\omega\omega$, under the partial
ordering given by eventual domination.
Remember $U^{\Gamma_\infty}=\omega$. Let $a_i=f_i/\FF$ for
$i<\omega_1$.  Let $b\in V^{\Gamma^*}$ be represented by the sequence
$b_n$ of elements of $V$ in $\Gamma_\infty$.  Let $B_n$ be the subset
of $U^{\Gamma_\infty}$ coded by $b_n$; we may suppose it is never empty.  Define
$g(n)=\sup\,B_n$ and let $i$ be chosen so that $f_i$ dominates $g$
eventually.  Then off a finite set we have $f_i(n)\not\in B_n$, and
hence in $\Gamma^*$, $a_i$ and $b$ are unlinked. \qed

\proclaim{1.7\\Proposition}
Any model $N$ of ZFC with the properties (P3,P4) satisfies part (ii)
of Theorem 1.1.
\endproclaim

\proof: 
We must construct an ultrafilter $\FF$ on $\omega$ such that any
ultraproduct of finite structures with respect to $\FF$ is saturated.
The construction takes place in $\aleph_2$ steps; at stage
$\alpha<\aleph_2$ we have a filter $\FF_\alpha$ generated by a
subfilter of at most $\aleph_1$ sets $(B_i)_{i<\omega_1}$ containing
the cobounded subsets of $\omega$, and we have
a type $p=(\varphi_i)_{i<\omega_1}$ over some ultraproduct $\prod_k{A_k}/\FF$
of finite structures to realize.  (More precisely, since the filter
$\FF$ has not yet been constructed, the ``type'' $p$ is given as a set
of pairs $(\varphi_i(x;\bar y),\bar f^{(i)})$ where $\bar f^{(i)}=
\lk f^{(i)}_1,\ldots\rk$ with $f^{(i)}_j\in \prod_k{A_k}$, $\,p$ is closed
under conjunction, and $p$ is consistent in a strong sense: for each
$\varphi_i$ there is a function $g_i$ such that
$\varphi_i(g_i(n);f^{(i)}_1(n),\ldots)$ holds for all $n$ in some set
which has already been put into $\FF$.)  By (P4) we can arrange the
construction so that at a given stage $\alpha$ we only have to deal
with one such type.

By (P3) there is a function $f\in \prod_k{A_k}/\FF$ such that for all
$i,j<\omega_1$, the set $\{n\in B_j:f(n)=g_i(n)\}$ is infinite, where
$g_i$ witnesses the consistency of $\varphi_i$.  We adjoin to $\FF$ all
of the sets $X_i=\{n\in \omega:\varphi_i(f(n);f^{(i)}_1(n),\ldots)\}$.
The resulting filter is nontrivial, and is again generated by at most
$\aleph_1$ sets.  Furthermore our construction ensures that $f/\FF$
will realize the type $p=\{\varphi_i(x;f^{(i)}_1/\FF,\ldots)\}$ in the
ultraproduct.

One may also take care as one proceeds to ensure that the filter
which is being constructed will be an ultrafilter.\qed

\proclaim{1.8\\Outline of the Construction}
\endproclaim

In the remainder of this section we will manufacture a model $N$ of
ZFC with the properties P1-P4 specified in 1.5.  We will use a
countable support iteration of length $\omega_2$ of
$^\omega\omega$-bounding proper forcing notions of
cardinality at most $\aleph_1$, starting from a model $M$ of GCH. (See
the Appendix for definitions and an outline of relevant results.)  By
[Sh 177] or [Sh f] VI\S2 or A2.3 here, improving the iteration theorem of [Sh b, Theorem
V.4.3], countable support iteration preserves the property:
$$\hbox{``$\up\omega\omega$-bounding and proper''.}$$ Thus every
function $f:\omega\to\omega$ in $N$ is eventually dominated by one in
$M$, and property P2 follows: $^\omega\omega$ has true cofinality
$\omega_1$ in $N$.  Our construction also yields P4:
$2^{\aleph_1}=\aleph_2$.  The other two properties are more
specifically combinatorial, and will be ensured by the particular
choice of forcing notions in the iteration.  The next two propositions
state explicitly that suitable forcing notions exist to ensure each of
these two properties; it will then remain only to prove these two
propositions.

\proclaim{1.9\\Proposition}

Suppose that $(A_n)_{n<\omega}$ is a collection of finite sets with $|A_n|\to\infty$, and
$g:\omega\to\omega$ with $g(n)\to\infty$.  Then there is a proper
$^\omega\omega$-bounding forcing notion $\PP$ such that for some $\PP$-name $\name H$ the
following holds in the corresponding 
generic extension:

\item{}$\name H$ is a function with domain $\omega$ with $\name H(n)\includedin
A_n$ and $|\name H(n)|\le g(n)$ for all relevant $n$, and for every
$f\in\prod_nA_n$ in the ground model, we have $f(n)\in \name H(n)$ if $n$ is sufficiently large
(depending on $f$). 
\endproclaim

\page
\proclaim{1.10\\Proposition}
Suppose $M$ is a model of ZFC, and $(A_k:k<\omega)$ is a sequence of
finite sets in $M$.  Then there is an $^\omega\omega$-bounding proper
forcing notion such that in the corresponding generic extension we
have a function $\name \eta\in \prod_k A_k$ satisfying: for all $f\in
\prod_k A_k$ and infinite $B\includedin \omega$, both in $M$, $\name
\eta$ agrees with $f$ on an infinite subset of $B$.

\endproclaim

We give the proof of Proposition 1.10 first.


\proclaim{1.11\\Definition}
\endproclaim

For $\AA=(A_k:k<\omega)$ a sequence of finite sets of natural numbers,
for simplicity $|A_k|\ge 2$ for every $k$,
let $\QQ(\AA)$ be the set of pairs $(T,K)$ where
$T\includedin \up\omega\omega$
is a tree and $K:T\to\omega$, such that for all $\eta$ in $T$ we have:

\item{1.} $\eta(l)\in A_l$ for $l<\len(\eta)$.
\item{2.} For any $k\ge K(\eta)$ and $x\in A_k$
there is $\rho$ in $T$ extending
$\eta$ with $\rho(k)=x$.

We take $(T',K')\ge (T,K)$ iff $T'$ is a subtree of $T$. By abuse of notation, we may
write ``$T$'' for ``$(T,K)$'' with $K(\eta)$ the minimal possible
value, and we may ignore the presence of $K$ in other ways. 

We use $\QQ(\AA)$ as a forcing notion: the intersection of a generic set of conditions
defines a function $\name \eta\in \prod_k A_k$, called the {\sl generic branch}.

We also define partial orders $\le_m$ on $\QQ(\AA)$ as follows.  $T\le_m T'$ iff $T\le T'$ 
and:
\item{1.}  $T\intersect {}^{ m\ge}\omega= T'\intersect {}^{m\ge}\omega$;
\item{2.}  $K(\eta)=K'(\eta)$ for $\eta\in T\intersect {}^{m\ge}\omega$.

Note the {\sl fusion} property: if $(T_n)$ is a sequence of conditions with $T_n\le_n
T_{n+1}$ for all $n$, then $\sup\ T_n$ exists (and is a condition).  We pay attention to
$K$ in this context.

\proclaim{1.12\\Remark}
\endproclaim

With the notation of 1.11, $\QQ(\AA)$ forces:

\item{} For any $f\in \prod_k A_k$ and infinite $B\includedin\omega$, both in the ground
model, the generic branch $\name \eta$ agrees with $f$ on an infinite subset of $B$.


\proclaim{1.13\\Proof of Proposition 1.10}
\endproclaim

It suffices to check that $\QQ(\AA)$ is an $^\omega\omega$-bounding proper
forcing notion.  We claim in fact:

\splitline 
\qquad Let $(T,K)\in \QQ(\AA)$, $m<\omega$, and let
$\name\alpha$ be a $\QQ(\AA)$-name for an ordinal. Then
//\qquad
there is $T', T\le_m T'$ such that for some
finite set $w$ of ordinals, $T'\forces\,``\name\alpha\in
w$''.//{$(*)$}

This condition implies that $\QQ(\AA)$ is $^\omega\omega$-bounding,
since given a name $\name f$ of a function in $^\omega\omega$, we can
find a sequence of conditions $T_n$ and finite sets $w_n$ of integers
such that $(T_n)$ is a fusion sequence (i.e. $T_n\le_n T_{n+1}$ for
all $n$) and $T_n\forces\,\hbox{``$\name f(n)\in w_n$''}$; then
$T=\sup\ T_n$ forces ``$\name f(n)\le \max\ w_n$ for all $n$''.

At the same time, the condition $(*)$ is stronger than Baumgartner's
Axiom A, which implies $\alpha$-properness for all countable $\alpha$.

It remains to check $(*)$. We fix $T$ (and the corresponding function
$K:T\to\omega$), $\name\alpha$, $m$ as in $(*)$.  For $\nu\in T$ let
$T^\nu$ be the restriction of $T$ to the set of nodes comparable with
$\nu$.  For $\nu$ in $T$, pick a condition $(T_\nu,K_\nu)$ by
induction on $\len(\nu)$
such that
$T_\nu\ge T^\nu$ and
$\eta\xxop\nu \& \nu \in T_\eta \implies T_\nu \ge T_\eta$ and
$T_\nu\forces\,\hbox{``$\name\alpha=\alpha_\nu$''}$ for some
$\alpha_\nu$.  We may suppose $K_\nu\ge K$ on $T_\nu$.  Set $k_0=m$,
and define $k_l$ inductively by
$$k_{l+1}=:
\max(k_l+1,\max\{K_\eta(\eta)+1:\eta\in T\intersect
{}^{k_l}\omega\}).$$ Let $(\eta_j)_{j=2,\dots,N}$ be an enumeration of
$T\intersect {}^{\le k_1}\omega$.  (It is convenient to begin counting
with 2 here.)  For $\nu\in T$ with $\nu\restrictedto k_1=\eta_j$, we
will write $j=j(\nu)$.

Let $T'$ be: $$\{\eta:\exists \nu\in T\hbox{ extending }\eta,\ \len(\nu)\ge k_N, \hbox{ and }
\nu\in T_{\nu\restrictedto k_{j(\nu)}}\}$$

Observe that for $\eta$ of length at least $k_N$, the only relevant
$\nu$ in the definition of $T'$ is $\eta$ itself.  That is, $\eta\in
T'$ if and only if $\eta\in T_{\eta\restrictedto k_{j(\eta)}}$.  In particular
$T'$ is a condition (with $K'(\eta)\le K_{\eta\restrictedto
k_{j(\eta)}}(\eta)$ for $\len(\eta)\ge k_N$).  Also, since $T'\intersect
{}^{k_N\ge}\!\omega \includedin \Union \{
T_\nu:\nu\in T\intersect{}^{ k_N\ge}\omega\}
$,
we find
$T'\forces\,\hbox{``}\name\alpha\in\{\alpha_\nu:\nu\in T\intersect
{}^{k_N\ge}\omega\}\hbox{''}$.
Notice also that $T'\restrictedto k_1=T\restrictedto k_1$.

The
main point, finally, is to check that we can take $K'=K$ on
$T'\intersect {}^{ m\ge}\omega$. Fix $\eta_j\in T'\intersect{}
^{m\ge}\omega$, $k\ge K(\eta_j)$, and $x\in A_k$;
we have to
produce an extension $\nu$ of $\eta_j$ in $T'$, with $\nu(k)=x$.
Let $\eta_h$ be an extension of $\eta_j$ of length $k_1$, such that
$\eta_h$ has an extension $\nu\in T$ with $\nu(k)=x$.
If $k< k_h$, then $\nu\restrictedto (k+1)\in T'$, as required.

Now suppose $k\ge k_{h+1}$, and let $\eta$ be an extension of $\eta_h$
of length $k_h$.  Then $T_\eta\includedin T'$, and $k\ge
K_\eta(\eta)$.  Thus a suitable $\nu$ extending $\eta$ exists.

We are left only with the case: $k\in[k_h,k_{h+1})$.  In particular
$k\ge k_2$, so $k>K(\eta_h)$ for all $\eta_h$ in
 $T\intersect \up{ k_1\ge}\omega$.
This means that any extension of $\eta_{h'}$  of $\eta$ of length $k_1$ could be
used in place of our original choice of $\eta_h$.
Easily there is such $h'\neq h$ (remember $|A_k|\ge 2$ and demand on $K$).
But $k$ cannot lie in two intervals of the form $[k_h,k_{h+1})$, so
we must succeed on the second try.\qed

\proclaim{1.14\\Logarithmic measures}
\endproclaim

We will define the forcing used to prove Proposition 1.9 in 1.16
  below.  Conditions will 
be perfect trees carrying extra information in the form of a (very
  weak) ``measure'' 
associated with each node.  These measures may be defined as follows. 

For $a$ a set, we write $P^+(a)$ for $P(a)\setminus\{\emptyset\}$.  A {\sl logarithmic measure} on
$a$ is a function $\norm\ :P^+(a)\to\Nn$ such that:
\item{1.} $x\includedin y\implies \normof x\le\normof y$;
\item{2.} If $x=x_1\union x_2$ then for some $i=1$ or $2$, $\normof
{x_i} \ge \normof x -1$.

By (1), $\norm\ $ has finite range.  If $a$ is finite (as will generally be the case in
the present context), one such logarithmic measure is $\normof
  x=\lfloor \ln_2|x|\rfloor$.

\proclaim{1.15\\The forcing notion $\LT$}
\endproclaim

We will force with trees such that the set of successors of any node carries a specified
logarithmic measure; the measures  will be used to
prevent the tree from being pruned too rapidly.
The formal definition is as follows.

\noindent 1. $\LT$ is the set of pairs $(T,t)$ where:

\item{1.1.} $T$ is a subtree of $^{\omega>}\omega$ with finite stem;
this is the longest branch in $T$ before ramification occurs.  We call
the set of nodes of $T$ which contain the whole stem the {\sl
essential part} of $T$; so $T$ will consist of its essential part
together with the proper initial segments of its stem. We denote the
essential part of $T$ by $\ess(T)$.

\item{1.2.} $t$ is a function defined on the essential part of $T$, with
$t(\eta)$ a logarithmic measure on the set $\sucs_T(\eta)$ of all
successors of $\eta$ in $T$; we often write $\norm\ _\eta$ (or
possibly $\norm\ ^T_\eta$) for $t(\eta)$. For $\eta $ a proper initial
segment of the stem of $T$, we stipulate $t(\eta)[\sucs(\eta)]=0$.

{\parindent=11pt \item{2.} The partial order on $\LT$ is defined by:
$(T_2,t_2)\ge(T_1,t_1)$ iff
$T_2\includedin T_1$, and for $\eta\in T_2$
$t_2(\eta)$ is the restriction of $t_1(\eta)$ to $P^+(\sucs_{T_2}(\eta))$.

\item{3.} We define $\LT^{[(T,t)]}$ to be $\{(T', t')\in\LT: (T',
t')\ge (T,t)\}$ with the induced order.  Similarly for $\LT^f$,
$\LT_d$, and $\LTfd$ (see below).

}

\proclaim{1.16\\The forcing notion $\LTfd$}
\endproclaim

$\LT^f$ is the set of pairs $(T,t)\in\LT$ in which $T$ has only {\sl finite} ramification
at each node.

$\LT_d$ is the set of pairs $(T,t)\in\LT$ such that for any $m$, every
branch of $T$ is almost contained in the set
$\{\eta\in T:\forall \nu \ge \eta\
\normof{\sucs_T(\nu)}_\nu\ge m\}$
(i.e. the set difference is finite).

$\LTfd$ is $\LT^f\intersect\LT_d$.  For $T\in \LT^f$, an equivalent
condition for being in $\LT^f_d$ is: $\lim_k
\inf\{\normof{\sucs_T(\eta)}_\eta: \len(\eta)=k\}=\infty$. Note:
 $\LT^f_d$ is an upward closed subset of $\LT_d$.


We make an observation concerning fusion in this connection.
Define:

\item{1.} $(T_1,t_1)\le^* (T_2,t_2)$
if $(T_1,t_1)\le (T_2,t_2)$ and in addition for all  $\eta\in
\ess T_2$, $\normof{\sucs_{T_2}(\eta)}^{T_2}_\eta\ge\normof
{\sucs_{T_1}(\eta)}^{T_1}_\eta-1$.
\item{2.} $(T_1,t_1)\le_m (T_2,t_2)$
if $(T_1,t_1)\le (T_2,t_2)$ and for all $\eta\in T_2$ with
$\normof{\sucs_{T_1}(\eta)}_\eta\ge m$, (so $\eta \in \ess (T_1)$) we have
$\normof{\sucs_{T_2}(\eta)}_\eta\ge m$ (hence $\eta\in\ess (T_2)$ when $m>0$).
\item{3.} $(T_1,t_1)\le^*_m (T_2,t_2)$
if $(T_1,t_1)\le_m (T_2,t_2)$ and for all $\eta\in T_2$ with
$\normof{\sucs(\eta)}^{T_1}_\eta\le m$, we have
$\sucs_{T_1}(\eta)\includedin T_2$.

If  $(T_n,t_n)$ is
a sequence of conditions in $\LTfd$ with
$(T_n,t_n)\le^*_{n}(T_{n+1},t_{n+1})$ for all~$n$, then $\sup
(T_n,t_n)$ exists in $\LTfd$.

We also mention in passing that a similar statement holds for $\LT_d$, with a more
complicated notation. Using arguments like those given here  one can
show that $\LT_d$ is also proper. This will not be done here.

For $\eta\in T, (T,t)\in \LT$ we let $T^\eta$ be the set of $\nu\in T$
comparable with $\eta$,
$t^\eta=t\uhr\ess (T^\eta$): so $(T, t)\le (T^\eta, t^\eta)$;
we may write $(T, t)^\eta$ or $(T^\eta, t)$ instead of $(T^\eta, t^\eta)$.


We will now restate Proposition 1.9 more explicitly, in two parts.

\page
\proclaim{1.17\\Proposition}

Suppose that $(A_n)_{n<\omega}$ 
is a collection of finite sets with $|A_n|\to\infty$, and
that $g:\omega\to\omega$ with $g\to\infty$.
  Then there is a condition $(T_0,t_0)$ in $\LTfd$
such that $(T_0,t_0)$ forces:
\splitline \qquad There is a function $\name H$
 such that $|\name H(n)|<g(n)$ for all $n$
 [more exactly,
//\qquad $|\name H(n)|<\max\{g(n), 1\}]$ , and for
 every $f$ in the ground model,
//\qquad\qquad $f(n)\in \name H(n)$ for $n$
sufficiently large.
{}
\endproclaim

\proof:

Without loss of generality $g(n)>1$ and $A_n$ is nonempty for every $n$.
Let $a_n=\{A\includedin
A_n:|A|=g(n)-1\}$,
$T_0=\Union_N\prod_{n<N}a_n$, and define
a logarithmic measure $\norm\ _n$ on
$a_n$ by $\norm x_n=\max\{l:$ if
$A'\includedin A_n$ has cardinality
$2^l$, then there is $A\in x$ containing
$A'\}$.  Set $t_0(\eta)=\norm\
_{\len\eta}$.

Obviously $(T_0,t_0)\in \LTfd$, 
(a pedantic reader will note $T_0\not\subseteq {}^{\omega>}\omega$ and rename)
 For a generic branch $\name\eta$ of $T_0$:
\splitline \qquad$(T_0,t_0)\forces_{\LTfd}$``$|\name\eta(n)|<g(n)$ for
all $n$;''\hfill 
//\qquad
$(T_0,t_0)\forces_{\LTfd}$``For $f$ in the ground model, $f(n)\in\name\eta(n)$ 
for all large $n$.''//{}
\qed

\proclaim{1.18\\Proposition}
The forcing notion $\LTfd$ is $^\omega\omega$-bounding and proper.
\endproclaim

It remains only to prove this proposition.

\proclaim{1.19\\Lemma}

If $(T,t)\in \LT_d$ and $W$ is a subset of $T$, then
there is some $(T',t')\in \LT_d$ with
$(T,t)\le^* (T',t')$ such that either:
\splitline $(+)$\qquad every branch of $T'$  meets $W$;
or else
//$(-)$\qquad $T'$ is disjoint from $W$.
\hfill//{}

\endproclaim

\proof: Let $T^W$ be the set of all $\eta\in T$ for which there is a condition
$(T',t')$ such that $T'$ has stem
$\eta$, $(T^{\eta
},t)\le^* (T',t')$, and every infinite branch
of $T'$ meets $W$.
($T^{\eta}$ is the set of $\nu\in T$ comparable to $\eta$; so it is a tree whose stem contains $\eta$.)

If the stem of $T$ is in $T^W$ we get (+).  Otherwise we will
construct $(T',t')\in \LT_d$ such that $(-)$ holds, $(T,t)\le^*
(T',t')$, and $T'\intersect T^W=\emptyset$.  For this we define
$T'\intersect {}^{n}\omega$ (and $t'=t\restrictedto \ess(T')$)
inductively.

If $n\le\len({\rm stem}(T))$ then we let $T'\intersect {}^{n}\omega$
be $\{{\rm stem}(T)\restrictedto n\}$.

So suppose that $n\ge \len({\rm stem\ } T)$ and that we have defined
everything for $n'\le n$.  Let $\nu\in T'\intersect {}^n\omega$, and in
particular, $\nu\not\in T^W$.  Let $a=\sucs_T(\nu)$, $a_1=a\intersect
T^W$, $a_2=a{\setminus}a_1$. Then for some $i=1$ or $2$,
$\norm{a_i}_\nu\ge\norm 
a_\nu-1$.

Since $\nu\notin T^W$, it follows easily that $\norm{a_1}_\nu<\norm
a_\nu-1$; otherwise one pastes together the conditions
$(T_{\nu'},t_{\nu'})$ associated with $\nu'\in a_1$ to show $\nu\in
T^W$.  Thus $\norm{a_2}_\nu\ge\norm{a}_\nu-1$.  Let
$T'\intersect(\sucs_T(\nu))$ be $a_2$.  As we can do this for all
$\nu\in T'\intersect {}^n\omega$, this completes the induction
step.\qed


\proclaim{1.20\\Lemma}

If $\name \alpha$ is an $\LTfd$-name of an ordinal, $(T,t)\in\LTfd$,
$m<\omega$, and $\normof{\sucs_T \eta}_\eta>m$ for $\eta\in \ess(T)$,
then there is $(T',t')\in\LTfd$ with $(T,t)\le_m(T',t')$,
and a finite set $w$ of ordinals, such that
$(T',t')\forces_{\LTfd}$``$\name\alpha\in w$''.
\endproclaim

\proof:

Let $W$ be the set of nodes $\nu$ of $T$ for which 
there is a condition $(T_\nu,t_\nu)$ with $(T_\nu,t_\nu) \down m\ge
(T^{\nu},t^{\nu})$ such that
$(T_\nu,t_\nu)$ forces a value on $\name \alpha$.
We claim that for any
$(T_1,t_1)\,{}\up*\mkern-5mu\ge (T,t)$,
 $T_1$ must meet $W$. Indeed, fix $(T_2,t_2) \ge
(T_1,t_1)$ forcing ``$\name \alpha=\beta$'' for some $\beta$. Then for some
$\nu \in T_2$, all extensions $\eta$ of $\nu$ in $T_2$ will satisfy
$\norm{\sucs_{T_2}(\eta)}_\eta \ge m$, and $(T_2,t_2)^{\nu}$
witnesses the fact that $\nu\in W$.  Thus if we 
apply Lemma 1.19, the alternative $(-)$ is not possible. 

Accordingly we have some $(T_1,t_1) \,{}\up*\mkern-6mu\ge (T,t)$ such that every
branch of $(T_1,t_1)$  meets $W$.
Let $W_0$ be the set of
minimal elements of $W$ in $T_1$.  Then $W_0$ is
finite.  For $\nu\in W_0$ select $(T_{\nu},t_{\nu})$ with
$(T_{\nu},t_{\nu}) \down m\!\ge(T,t)^\nu$ and
$(T_{\nu},t_{\nu})\forces$ ``$\name\alpha=\alpha_\nu$''
for some $\alpha_\nu$.  Form $T'=\Union\{T^{\nu}:\nu\in W_0\}$.\qed

\proclaim{1.21\\Lemma}

If $(T,t)\in \LTfd$, $\name \alpha$ is an $\LTfd$-name of an ordinal,
$m<\omega$,  then there is $(T',t')\in \LTfd$ with
$(T,t)\le_m^* (T',t')$, and a finite set of ordinals
$w$, such that
$(T',t')\forces\,\hbox{``$\name\alpha\in w$''}$.
\endproclaim

\proof:

Fix $k$ so that $\normof{\sucs(\eta)}_\eta>m$ for $\len(\eta)\ge k$.
Apply 1.20 to each $T^\nu$ for $\nu\in T$ of length $k+1$.\qed

\proclaim{1.22\\Proof of 1.18}
\endproclaim

As in 1.13, using 1.21.\qed
\vskip 5 pt

This completes the verification that the desired model $N$ can be constructed by iterating forcing.
\par

\page
\section 2.  Nonisomorphic ultraproducts of finite
models.

We continue to use the bipartite graphs $\Gamma_{k,l}$ introduced in
1.3. 
Varying the forcing used in \S1, we will get:

\proclaim{2.1\\Theorem}

Suppose that $V$ satisfies CH, and that $(k_n,l_n),(k_n',l_n')$ are
monotonically increasing sequences of pairs (and
$2<l_n'<k'_n<l_n<k_n<l'_{n+1}$)  such that:
\medskip
\settabs\+xxxxxxxxxx&\cr
\+$(1)$&\hfill$k'_n/l'_n\to\infty$;\hfill\cr
\medskip
\+$(2)$&\hfill$(k_n/l_n) > (k'_n)^{ndl'_n}$, for each $d>0$, for $n$
large enough;\hfill\cr
\medskip
\+$(3)$&\hfill$\ln l'_n>k_{n-1}^n$.\hfill\cr
\medskip

Then there is a proper forcing $\PP$ satisfying the $\aleph_2$-cc, of
size $\aleph_2$, such that in $V^\PP$ no two ultraproducts $\uprod{\Gamma_{k_i,l_i}}_1$,
$\uprod{\Gamma_{k'_i,l'_i}}_2$ are isomorphic.
\endproclaim

More precisely, we will call a bipartite graph with bipartition $(U,V)$
$\aleph_1$-{\sl complete} if every set of $\omega_1$ elements of $U$
is linked to a single common element of $V$ (property $(\dagger)$ of
Proposition 1.6), and then our claim is that in $V^\PP$, no nonprincipal
ultraproduct of the first sequence $\Gamma_{k_n,l_n}$ is
$\aleph_1$-complete, and every nonprincipal ultraproduct of the second sequence
$\Gamma_{k_n',l_n'}$ is; furthermore, as indicated, this phenomenon can be
controlled by the rates of growth of $k$ and of $l/k$.

\proclaim{2.2\\Definition}
\endproclaim

Let $f,g$ be functions in $^\omega\omega$.  A model $N$ of ZFC is
$(f,g)$-{\sl bounded} if for any sequence $(A_n)_{n<\omega}$ of finite
sets with $|A_n|=f(n)$, there are $\aleph_1$ sequences
$\BB_i=(B_{i,n}:n<\omega)$, indexed by $i<\omega_1$,  with:
\medskip
\+$(1)$&$B_{i,n}\includedin A_n$ for all $n$;\cr
\medskip
\+$(2)$&For all $i<\omega_1$, $|B_{i,n}|<g(n)$ eventually;\cr
\medskip
\+$(3)$&$\Union_i\prod_n B_{i,n}=\prod_n A_n$ in $N$.\cr
\medskip

\proclaim{2.3\\Lemma}
Let $(k_n),(l_n)$ be sequences with $l_n,k_n/l_n\to\infty$, and
let $f(n)={k_n\choose l_n}$, $g(n)=k_n/l_n$.
Suppose that $N$ is a model of ZFC which is $(f,g)$-bounded.
Then no ultraproduct
 $\prod_n{\Gamma_{k_n,l_n}}/\FF$
can be $\aleph_1$-complete.
\endproclaim

\proof:
Let $\BB_i$ have properties (1-3) of 2.2 with respect to
$A_n=V_{k_n,l_n}$.  For each $i$, choose $a_i\in \prod_n U_{k_n,l_n}$ so that
$a_i(n)$ is not linked to any $b\in B_{i,n}$, as long as $|B_{i,n}|< g(n)$
(so $l_n|B_{i,n}|<k_n$). Then $a_i/\FF (i<\omega_1)$ cannot all be linked to
any single $b$ in $\prod_n{\Gamma_{k_n,l_n}}/\FF$, for any
ultrafilter $\FF$. \qed

\page
\proclaim{2.4\\Definition}
\endproclaim

For functions $f,g\in{} ^\omega\!\omega$ we say that a forcing notion
$\PP$ has the $(f,g)$-{\sl bounding property} provided that:

$$\vcenter{\hsize=0.7\hsize 
For any sequence $(A_k:k<\omega)$ in the ground model,
with $|A_k|=f(k)$, and any $\name\eta\in\prod_kA_k$ in the generic
extension, there is a ``cover'' $\BB=(B_k:k<\omega)$ in the ground
model with $B_k\includedin A_k$, 
$|B_k|<g(k)$ (more exactly, $<\Max\{g(k),2\})$,
 and $\name\eta(k)\in B_k$ for each $k$.}\leqno(*)$$

Similarly a forcing notion has the  $(\unF,g)$-bounding property, for
$\unF$ a collection of functions, if it has the $(f,g^\epsilon)$-bounding
property for each $f\in \unF$ and each $\epsilon>0$.
In this terminology, notice  that $(\{f\},g)$-bounding is a stronger
condition than $(f,g)$-bounding.

\proclaim{2.4A\\Definition}
\endproclaim
Call a family $\unF$ {\sl $g$-closed} if it satisfies the following
two closure conditions:

\item{1.} For $f\in \unF$, the function $F(n)=\prod_{m<n}(f(m)+1)$
 lies in
$\unF$;
\item{2.} For $f\in \unF$, $f^g$ is in $\unF$.

\proclaim{\\Proof of 2.1}
\endproclaim
We build a model $N$ of ZFC by an iteration of length $\omega_2$ with
countable support of 
proper forcing notions with the $(\unF,g)$
bounding property for a suitable family $\unF$, all of which are of the
form $\LTfdat {(T,t)}$; and
we arrange that all of the forcing notions of this form which are
actually $(\unF,g)$-bounding will occur cofinally often.  (In order to
carry this out one actually makes use of auxiliary functions
$(f_1,g_1)$ with $f_1$ eventually dominating $\unF$ and $g_1$
eventually dominated by any positive power of $g$, but these details
are best left to the discussion after 2.5.)

One can show that a countable support iteration of proper
$(\unF,g)$-bounding
forcing notions is again $(\unF,g)$-bounding.  This is an instance of a
general iteration theorem of [Sh f, VI] but we make our presentation 
self-contained by giving a proof in the appendix-A2.5.
If we force over a ground model with CH (so that CH holds at
intermediate points in the iteration) then our final model is
$(\unF,g)$-bounded,
and by 2.3  no ultraproduct of the
$\Gamma_{k_n,l_n}$ can be $\aleph_1$-complete.

One very important point still remains to be checked. It may be
formulated as follows.

\proclaim{2.5\\Proposition}
Let $f_0,g_0,h:\omega\to\omega\setminus\{0,1\}$ and suppose that
$(A_n)_{n<\omega}$ is a
sequence of finite nonempty sets with $|A_n|\to\infty$. Assume:

$$\prod_{m\le n}|A_m|^{h(m)}< g_0(n){{\rm\;for\; every}}\; n
{{\rm \; large \; enough}}
;\leqno(1)$$
$${\ln h(n)\over\ln\prod_{i<n} f_0(i)}\to \infty.\leqno(2)$$

Then there is a condition $(T,t)\in \LTdf$ such that $\LTfdat {(T,t)}$
is $(f_0,g_0)$-bounding and $(T,t)$ forces:
\item{}There is a function
$\name H$ such that $\name H(n)\includedin A_n$, $|\name H(n)|<h(n)$
for all $n$, and for every $f$ in the ground model, $f(n)\in \name
H(n)$ for $n$ sufficiently large.  \endproclaim

\proclaim{Continuation of the proof of 2.1:}\endproclaim
We will now check that the proof of theorem 2.1 can be completed using
this proposition.

We set
$f^*(n)={k_n\choose l_n}$, $g(n)=k_n/l_n$, $h(n)=l_n'$, and
$A_n=U_{k'_n,l'_n}$.  (So $|A_n| = k_n'$.)
Let $\unF_0$ be the set of increasing functions $f$ satisfying
$$\lim_{n\to\infty} \ln h(n)/(g^{d}(n-1)\ln f(n-1))\to \infty \qquad
\hbox{for all $d>0$}.$$

If $f_0\in \unF_0$ and $g_0$ is a positive power of $g$, then conditions
(1,2) of 2.5 hold by condition (2) of 2.1
(for (2) of 2.5 note for $d=2$ that $g^d(n-1)>n$).
Furthermore $\unF_0$ is $g$-closed (this uses the fact that
$g(n)\ge n$ eventually by (2) of 2.1), and $f^*\in\unF_0$.
 By diagonalization
find $f_1,g_1$ satisfying $(1,2)$ of 2.5
so that $f_1$ eventually dominates
any function in the $g$-closure of $f^*$, and $g_1$ is eventually
dominated by any positive power of $g$.  Apply the proposition to
$(f_1,g_1,h)$ and observe that an $(f_1,g_1)$-bounding forcing notion
is  $(g\hbox{-closure of $f^*$},g)$-bounding.
We let $\unF=g$-closure of $\{f^*\}$.

Forcing with the corresponding
$\LTdfat{(T,t)}$ produces a branch  $\name H$ so that if $\name H(n)$ is
thought of as an element $b_n\in V_{k'_n,l'_n}$, then for all $f\in
\prod_n A_n$ in the ground model, and any ultrafilter $\FF$ on
$\omega$, $f/\FF$ is linked to $\name H(n)/\FF$ in
$\prod_n{\Gamma_{k'_n,l'_n}}/\FF$.

\proclaim{2.6\\Terminology}
\endproclaim

A logarithmic measure $\norm\ $ on $a$ is called $m$-{\sl additive} if
for every choice of $(a_i)_{i<m}$ with
$\Union_i a_i=a$, there is $i<m$ with $\norm{a_i}\ge\norm{a}-1$.

\proclaim{2.7\\Lemma}
Suppose $f,g:\omega\to\omega\setminus\{0,1\}$, $(T,t)\in \LTfd$, and:
\item{i.} for every
$\eta\in \ess (T)$, $t(\eta)$ is $\prod_{i<\len\eta}f(i)$-additive;
\item{ii.} for every $n$ we have $|T\intersect\,{}^{(n+1)}\omega|<g(n)$.

\noindent Then $\LTfdat{(T,t)}$ is $(f,g)$-bounding.
\endproclaim

\proof:
Let $F(n)=\prod_{i<n}f(i)$.  Suppose that $(A_n)_{n<\omega}\in V$,
$|A_n|\le f(n)$, and
$(T,t)\forces\,$``$\name\eta\in\prod_n A_n$''.  By fusion as in 1.19--1.22 there is $(T',t')\in
\LTfd$ with $(T',t')\ge(T,t)$ such that for every $n$ the set
$$W=:
\{\nu\in T':(T'^{\nu},t') \hbox{ forces a value on } \name \eta(n)\}$$
meets every branch of $(T',t')$.

For each $n$, choose $N(n)$ large enough that $(T'^{\nu},t')$ forces a value
$\eta^n_\nu$ on $\name\eta\restrictedto n$ for each $\nu\in T'\intersect
{}^{N(n)}\omega$.  Thus $\eta^n_\nu\in \prod_{i<n} A_i$.  By downward
induction on $k<N(n)$,
for $\nu\in T'\intersect {}^k\omega$ choose
$\eta^n_\nu\in {}^n\omega$ and $s(\nu,n)\includedin \sucs_{T'}(\nu)$
so that:
$$\normof{s(\nu,n)}_\nu \ge \normof{\rsucc_{T'}(\nu)}_\nu-1;
\quad \eta^n_\nu\restrictedto\min\{k,n\}=\eta^n_{\nu'}\restrictedto
\min\{k,n\} \hbox{ for } \nu'\in s(\nu,n).$$
Since $|\{\eta_{\nu'}\restrictedto \min(k,n):\nu'\in\sucs_{T'}(\nu)\}|\le F(k)$
and $\norm\ _\nu$ is $F(k)$-additive, this is easily done.
Let $T'_n=\{\nu\in T':(\forall l<\len(\nu)\cap N(n))\ \nu\restrictedto (l+1)\in
s(\nu\restrictedto l,n)\}$.

We now define $T''\includedin T'$ so that for all $k$ the set $X_k$ of
$n$ for which
$T''\intersect \up{k\ge}\omega = T'_n\intersect\up{k\ge}\omega$
is infinite.
For this we proceed by induction on $k$.
If $T''\intersect \up{ k\ge}\omega$ has been defined, then we can select
$X\includedin X_k$ infinite such
that for $n\in X$ and $\nu\in T''\intersect  \up{k}\omega$, $s(\nu,n)=s(\nu)$ is
independent of $n$.  We then define
$$T''\intersect \up{(k+1)}\omega=\{\nu\in T'\intersect\up{k+1}\omega:
\nu\restrictedto k\in T''\intersect \up {k}\omega\hbox{ and }\nu\in s(\nu\restrictedto k)\}$$

Observe that $(T'',t\restrictedto T'')\;{}^*\mkern-6mu\ge (T',t')$, and $(T'',t\restrictedto
T'')$ forces:
$$\hbox{``For any $k$, if $n\in X_{N(k)}$ and $n\ge k$,
then $\name\eta\restrictedto k=\eta^n_\nu\restrictedto k$ for some
$\nu\in T''\intersect \up k\omega$''.}$$
Indeed, for any $\nu'$ of length $N(k)$ in $T''$,
if $\nu'\in T'_n$ then $\eta^k_{\nu'}=\eta^n_{\nu'}\restrictedto
k=\eta^n_{\nu'\restrictedto k}\restrictedto k$.
Since $|T''\intersect {}^{k+1}\omega|\le |T\intersect {}^{k+1}\omega|<g(k)$,
this yields the stated bounding principle.
\qed

\page
\proclaim{2.8\\Proof of 2.5}
\endproclaim

Let $F_0(n)=\prod_{i<n} f_0(i)$.
Let $a_n=\{A\includedin
A_n:|A|=h(n)-1\}$,
$T_0=\Union_N\prod_{n<N}a_n$, and define
a logarithmic measure $\norm\ _n$ on
$a_n$ by: for $a\subseteq a_n$
$$\norm a_n=\max\{l: \hbox{
for all $A'\includedin A_n$ of cardinality $\le F_0(n)^l$, there is
$A\in a$ containing $A'\}$.}$$
Set $t_0(\eta)=\norm\
_{\len\eta}$.

Obviously $\norm\ _n$ is $F_0(n)$-additive and $|T\intersect
{}^{(n+1)}\omega|=\prod_{m\le n}(|A_m|)^{(h(m)-1)}$ which
is (by condition (1) of 2.5)
$<g_0(n)$, so $\LTdfat{(T_0,t_0)}$ is $(f_0,g_0)$-bounding
by lemma 2.7.

 We need to check that
$\norm{a_n}_n\to\infty$:
$$\norm{a_n}_n=\max\{l:F_0(n)^l<h(n)\}\sim {\ln h(n)\over\ln F_0(n)}.\eqno$$
So (2) from 2.5 guarantees it.
\qed

\page
\section 3.  Adding Cohen reals creates a bad ultrafilter.

In this section we show that a weaker form of the results in
$\S\S1,2$ is obtained just by adding $\aleph_3$
Cohen reals to a suitable ground model.
This result was actually the first one obtained in this direction.
This construction is also used in [Sh 345] and again in [Sh 405].

\proclaim{3.1\\Theorem}
If we add $\aleph_3$ Cohen reals to a model of
 $[2^{\aleph_i}=\aleph_{i+1} \ (i=1,2)\ \&\
\diamond_{\{\delta<\aleph_3:\cof \delta=\aleph_2\}}]$, then there
will be a nonprincipal ultrafilter $\FF$ on $\omega$ and two
sequences of pseudorandom finite graphs $(\Gamma^1_n)$, $(\Gamma^2_n)$ such that
$\prod_n{\Gamma^1_n}/\FF\not\iso\prod_n {{\Gamma^2_n}}/\FF$. In fact the
same result will apply if the sequences $\Gamma^1_n$, $\Gamma^2_n$ are
replaced by any subsequences.
\endproclaim

Here we call a sequence $(\Gamma^1_n)$ of finite graphs {\sl
pseudorandom} if the theory of $\Gamma^1_n$ converges fairly rapidly
to the theory of
the random infinite graph; cf. 3.4 below.  The only condition needed
on the two sequences in Theorem 3.1 is that the $\Gamma^1_m$ and
$\Gamma^2_n$ are of radically different sizes (3.5 below). As a variant
(with very much the same proof) we can take all $\Gamma^2_n$ equal to the
random infinite graph, keeping $(\Gamma^1_n)$ a sequence of pseudorandom
finite graphs, and obtain the same result for a suitable ultrafilter.

\proclaim{3.2\\ Corollary}

Under the hypotheses of Theorem 3.1 there
are elementarily equivalent countable graphs $\Gamma^1_\omega$,
$\Gamma^2_\omega$ and a nonprincipal ultrafilter $\FF$ on $\omega$
with $\upower {(\Gamma^1_\omega)} \not\iso \upower {(\Gamma^2_\omega)}$.
\endproclaim

This is proved much as in Remark 1.4, noting that large pseudorandom graphs
are connected of diameter 2.

\proclaim {3.3\\Remark}
\endproclaim

With more effort we can replace the hypotheses on the ground model in
Theorem 3.1 by:
$$2^{\aleph_i}=\aleph_{i+1}\ (i=0,1)\,\&\,
\diamond_{\{\delta<\aleph_2:\cof\delta=\aleph_1\}},$$ adding only
$\aleph_2$ Cohen reals.  In the definition of $\AP$ below, \name$\FF$
would then not be an arbitrary name of an ultrafilter; instead $\AP$
would be replaced by a family of $\aleph_1$ isomorphism types of
members of $\AP$, (using $\aleph_0$ in place of $\aleph_1$ in clause
3.8 (i) below) which is closed under the operations used in the proof.

The same approach allows us to eliminate the $\diamond$ from Theorem
3.1.  With the modified version of $\AP$ and $\aleph_3$ Cohen reals,
we can replace $\diamond_{\{\delta<\aleph_3:\cof \delta=\aleph_2\}}$
by $\diamond_{\{\delta<\aleph_3:\cof \delta=\aleph_1\}}$, which in
fact follows from the other hypotheses [Gregory, Sh 82].

We will not enlarge on these remarks any further here.

\proclaim{3.4 \\ Definition}
\endproclaim

A finite graph $\Gamma$ on $n$ vertices is {\sl sufficiently random} if:

\item{i.} For any two disjoint sets of vertices $V_1,V_2$ with
$|V_1\union V_2|\le  (\log n)/3$, there is a vertex $v$ linked to all
vertices of $V_1$, and none in $V_2$;

\item{ii.} For any sets of vertices $V_1,V_2$ with $|V_i|>3\log n$
there are adjacent and nonadjacent pairs of vertices in $V_1\times V_2$.
\endproclaim

\item{iii.} If $V_1,V_2,V$ are three disjoint sets of vertices and
$P\includedin V_1\times V_2$, with $|P|,|V|> 5 \log n$, and if all pairs
in $P$ have distinct first entries, 
then some $v\in V$ {\sl separates} some pair $(v^1,v^2)\in P$ in the
sense that:\ [$R(v^1,v)\iff\neg R(v^2,v)$].  Here $R$ is the edge
relation (in the appropriate graph).

For sufficiently large  $n$ most graphs of size $n$ are  sufficiently
random. We call any 
sequence of sufficiently random graphs of size tending to infinity a
sequence of pseudorandom graphs.  

(See [Bollobas] for background on random graphs.) 

\proclaim{3.5\\ Notation}
\endproclaim

\item {i.} $(\Gamma^1_n)$, $(\Gamma^2_n)$ are two sequences of
sufficiently random graphs such that for any $m,n$ we have
$\norm{\Gamma^1_m}>{\norm{\Gamma^2_n}}^5$ or
$\norm{\Gamma^2_n}>{\norm{\Gamma^1_m}}^5$.  ($\norm \Gamma$ is
the number of vertices of $\Gamma$.) 
These sequences are
kept fixed.  $\Gamma$ is the infinite random (homogeneous) graph.  If
we replace $\uprod{\Gamma^2_n}$ by $\upower{\Gamma}$ throughout, the
argument is much the same, with slight simplifications.

\medskip
\item {ii.} $\Pp$ is the forcing notion that adds $\aleph_3$ Cohen
reals to $V$.  $\name x_\alpha$ is the name of the $\alpha$-th Cohen
real as an element of $\up\omega\omega$.  For $\AA\includedin \aleph_3$,
$\Pp\restrictedto \AA$ denotes $\{\pp\in \Pp:\dom
\pp\includedin 
\AA\}$. 

\proclaim{3.6\\Discussion}
\endproclaim

Working in the ground model we will build a $\Pp$-name for a suitable
nonprincipal ultrafilter $\name \FF$.  We will view the reals $\name
x_\alpha$ as (for example) potential members
of the ultraproduct $\uprodname {\Gamma^1_n}$.  We will 
consider candidates $\name {y_\alpha}$ for (representatives of) their images under a putative
isomorphism, and defeat them by arranging (for example) that the set of $n$ for which
$$R(\name{x_\alpha}(n),\name{x_\beta}(n)) {\rm \ iff\ }
\false R(\name{y_\alpha}(n),\name{y_\beta}(n))$$
gets into $\name \FF$.  

Note however that this must be done for every two potential sequences
$(\name k^1(n))$ and $(\name k^2(n))$ indexing the ultraproducts
$\prod_n{\Gamma^1_{\name k^1(n)}}/\name\FF$,
$\prod_n{\Gamma^2_{\name k^2(n)}}/\name\FF$ to be formed.  At stage $\alpha$ we
deal with sequences $\name k^1_\alpha(n),\name k^2_\alpha(n)\in
V^{\Pp\restrictedto\alpha}$ (which are guessed by the diamond).  We require
$\{n:\name{x_\alpha(n)}\in \Gamma^{\epsilon_\alpha}_{\name
k^{\epsilon_\alpha}_\alpha(n)}\}\in \name\FF$ where
$\epsilon_\alpha\in\{1,2\}$ is a label, and another very important
requirement is that for any sequence $(\name{A_n}:n<\omega)\in
V^{\Pp\restrictedto\alpha}$ with $\name A_n\includedin
\Gamma^{\epsilon_\alpha}_{\name k^{\epsilon_\alpha}_\alpha}(n)$ and $|\name
A_n|/\norm{\Gamma^{\epsilon_\alpha}_{\name k^{\epsilon_\alpha}_\alpha(n)}}$ small enough, the set $\{n:\name
x_\alpha(n)\not\in \name A_n\} \in\name\FF$.  (This sort of condition
is an analog of the notion of a $\Gamma$-{\sl big} type in [Sh 107].)
It will be used in combination with clause (ii) in the definition
of sufficient randomness.  

The name $\name\FF$ is built by carefully amalgamating a large set of
approximations to the final object, using the combinatorial principle
$\diamond_{\aleph_2}$, which follows from the cardinal arithmetic [Gregory];
this method, which  was
illustrated in [Sh 107], is based on the theorem from  [ShHL 162].  (The
comparatively elaborate tree construction of [ShHL 162] can be
simplified in the presence of $\diamond$; it is designed to work when
$\aleph_2$ is replaced by a limit cardinal and $\diamond$ is weakened
to the principle $\rm Dl_\lambda$.)  In what follows, the connection
with [ShHL 162] is left somewhat vague; the details will be found in \S A3
of the Appendix. In particular, in \S A3.5 we show how the present
 $\AA\PP$ fits the framework of \S A3.1-3.

\proclaim{3.7\\A notion of smallness}
\endproclaim

If $\FF$ is a filter on $\omega$, $k\in\up\omega\omega$,
$\epsilon\in\{1,2\}$, then a sequence $(A_n:n<\omega)$ of subsets of
the $\Gamma^\epsilon_{k(n)}$
(i.e. $A_n\subseteq \Gamma^\epsilon_{k(n)})$
is $(\FF,k,\epsilon)$-{\sl slow} if there
is some $d$ such that $\FF$-lim
$\left[|A_n|/\left(\sqrt{\norm{\Gamma^\epsilon_{k(n)}}}\cdot
(\log\norm{\Gamma^\epsilon_{k(n)}})^d\right)\right] = 0$.  Later on we will deal
primarily with the case $\epsilon=1$, to lighten the notation, and we
will then write ``$(\FF,k)$-slow'' in place of ``$(\FF,k,1)$-slow''.

It should perhaps be emphasized that here (as opposed to \S2)
$\epsilon$ is merely a label.

\proclaim{3.8\\ Definition}
\endproclaim

We define the partially ordered set $\AP$ of approximations as
follows.  The intent is that the approximations should build the name of a
suitable ultrafilter $\name \FF$.
Recall that the sequences $(\Gamma^\epsilon_n)$ (with
$\epsilon\in\{1,2\}$) are fixed (3.5(i)).  Also bear in mind that the
ultrafilter must eventually ``defeat'' a potential isomorphism between
two ultraproducts
$\prod_n{\Gamma^\epsilon_{\name k^\epsilon(n)}}/\name \FF $.
\bigskip

1. An element $q \in \AP$ is a quadruple
$(\AA,
\name\FF,
{\fakebf\epsilon},
\name
{\hbox{\bf k}})
=(\AA^q,\name\FF^q,{\fakebf \epsilon}^q,\name
{\hbox{\bf k}}^q)$ where

\item{i.} $\AA\includedin
\aleph_3$ has cardinality $\aleph_1$; ${\fakebf
\epsilon}=(\epsilon_\alpha:\alpha\in \AA)$ with each
$\epsilon_\alpha$ an
element of $\{1,2\}$;
\item{ii.} \name $\FF$ is a
$\Pp\restrictedto \AA$-name of a nonprincipal ultrafilter on
$\omega$, and if we set
$\name\FF\restrictedto(\AA\intersect\alpha)=:\\\\\name\FF\restrictedto\{\name
X: \name X$ is a $\Pp\restrictedto (\AA\intersect \alpha)$-name for a  subset of $\omega\}$,
then $\name\FF\restrictedto (\AA\intersect\alpha)$
is a $\Pp\restrictedto (\AA\intersect\alpha)$-name for all $\alpha$;
\item{iii.} $\name {\bf k}${}$= (\name k_\alpha:\alpha \in
\AA)$ with $\name k_\alpha$ a $\Pp\restrictedto (\AA\intersect\alpha)$-name of a
function from $\omega$ to $\omega$;
\item{iv.} For each $\alpha\in \AA$,
and each $\Pp\restrictedto(\AA\intersect\alpha)$-name $(\name
A_n:n<\omega)$;

\itemitem{}if $\forces_{\Pp\restrictedto (\AA\intersect\alpha)}$
``$(\name A_n)_{n<\omega}$ is
$(\name\FF\restrictedto\alpha$,$\name k_\alpha$,$\epsilon_\alpha)$-slow''
then
\itemitem{}\ \ \ $\forces_{\Pp}\hbox{``$\{n:\name x_\alpha(n)\in
\Gamma^{\epsilon_\alpha}_{\name k_\alpha(n)}\setminus A_n\}\in \name\FF$''}$.

We write $\AA=\AA^q$, $\FF=\FF^q$, and so on, when necessary.
\medskip
2. We take $q\le q'$ if $\AA^q\includedin \AA^{q'}$ and
$q'\restrictedto \AA^q=q$.
\bigskip

Some further comment is in order here.  When we begin to check
that $\name\FF$ is indeed  the name of an ultrafilter such that for
any pair of sequences $\name k^1(n),\name k^2(n)$,
the ultraproducts $\uprodname
{\Gamma^\epsilon_{\name k^\epsilon(n)}}/\name\FF$
 are nonisomorphic, we will notice that
there is
an automatic asymmetry because the
sequences $(\Gamma^1_n)$ and $(\Gamma^2_n)$ are so different: 
on some set in $\name\FF$ we will have
$|\Gamma^\epsilon_{\name k^\epsilon(n)}| >
|\Gamma^{\epsilon^*}_{\name k^{\epsilon^*}(n)}|^5$ holding with
$\{\epsilon, \epsilon^*\}=\{1,2\}$ in some order.  The
parameter $\epsilon_\alpha$ in an approximation can be viewed as a guess as
to the direction in which this asymmetry goes (after adding Cohen
reals); the 
notion of an approximation includes a clause (iv) designed to be
useful when $\name k_\alpha$ coincides with a particular $\name k^\epsilon$
in the context just described.  

On the other hand, we could first use $\diamond$ to guess
$\epsilon_\alpha$, $\name k^{\epsilon_\alpha}_\alpha$, and many other
things; in this case we do not actually need to include these kinds of
data in the approximations themselves, though it would still be
necessary to mention them in clause (iv).  Alternatively, the set
$\AP$ could also be used as a forcing notion, without $\diamond$, and
in this case the $\epsilon$ and $\name k$ would have to be included.
So the version given here is the most flexible one.

\proclaim{3.9\\Claim (Amalgamation)}

1. Suppose that $q_0,q_1,q_2\in \AP$, 
$\AA^{q_1}\includedin \delta$,
$\AA^{q_2}=\AA^{q_0}\union\{\delta\}$, and $q_0\le q_1,q_2$.
Then we can find
$r\ge q_1,q_2$ in $\AP$.

2. If $q_1,q_2\in\AP$, 
$\alpha<\aleph_3$,
$\dom q_1\includedin \alpha$, and $q_2\restrictedto \alpha \le q_1$,
then there is
$r\ge q_1,q_2$ in $\AP$.

\endproclaim

\proof:

1:\\Let $\AA_i=\AA^{q_i}$, $\name\FF^i=\name\FF^{q_i}$,
$\AA=\AA_1\union\{\delta\}, \varepsilon=\varepsilon^{q_2}_\delta$ and
$\name k=\name k^{q_2}_\delta$.
 In particular $\name\FF^0\includedin
\name\FF^1,\name\FF^2$, and we have to combine them into one
ultrafilter $\name\FF$ in $V^{\Pp\restrictedto \AA}$.  The point
is to preserve 3.8(iv), that is to ensure that $\Pp\restrictedto \AA$ forces the
relevant family of
sets (namely, $\name\FF^1,\name\FF^2$, and sets imposed on us by 3.8(iv))
to have  the finite intersection property.

If $\pp\in
\Pp\restrictedto\AA$ forces the contrary, then after extending
$\pp$ suitably we may suppose that there is
a $(\Pp\restrictedto \AA_1)$-name  $\name a$ of a
member of \name $\FF^1$, a $(\Pp\restrictedto
\AA_2)$-name $\name b$ of a member of $\name
\FF^2$,
and
-- since $\AA_1=\AA\intersect
\delta$ --
a
$(\Pp\restrictedto \AA_1)$-name
$(\name A_n:n<\omega)$ forced by $\pp$ to be
$(\name\FF^1,\name k,\epsilon)$-slow  (as
in (iv) of 3.8) so that letting $\name c =
\{n<\omega:\name x_\delta(n)\in
\name \Gamma^{\epsilon_\delta}_{\name k_\delta(n)}
\setminus\name A_n\}$ we have:
$$\pp\forces_{\Pp\restrictedto \AA} \hbox{``$\name a
\intersect \name b \intersect \name c = \emptyset$''}.$$

\noindent (i.e. we used the fact that there are three kinds of
requirements of the form ``a set belongs to $\FF$ '',
each kind is closed under finite intersections).

Let $\pp_i=\pp\restrictedto \AA_i$ for $i=0,1,2$. To clarify the
matter choose $\Hh^0\includedin \Pp\restrictedto \AA_0$
generic over $V$ so that $\pp_0\in \Hh^0$.
Note that $\name k$ is a $(\Pp\restrictedto \AA_0)$-name (3.8(iii)).

In $V[\Hh^0]$, for each $n<\omega$ let
$$\name B_n[\Hh^0] = \{v\in
\Gamma^\epsilon_{\name k(n)}[\Hh^0]:\hbox{ For some
$\pp'_2\in\Pp\restrictedto \AA_2$ with $\pp'_2\ge\pp_2$ and
$\pp'_2\restrictedto \AA_0\in\Hh^0$,}$$
$$\pp'_2
\forces_{\Pp\restrictedto
\AA_2}\hbox{``$\name x_{\delta}(n)= v$ and $n\in\name b$''}\}.$$
Then $(\name B_n:n<\omega)$ is not $(\name
\FF^0\restrictedto\delta,\name k,\epsilon)$-slow, since $(\name
B_n:n<\omega)$ is a $\Pp\restrictedto \AA_0$-name,
 $q_2\in\AP$, and
 $\pp_2\forces\hbox{``For $n\in\name b$,
 $\name x_\delta(n)\in\name B_n$''}$ (and (iv) of 3.8(1)).

Also in $V[\Hh^0]$, let
$\name b^+[\Hh^0]=\{n:\hbox{ for every $\pp'_0\in \Hh^0$,
$\pp'_0\union\pp_2\notforces$}
\hbox{``$n\not\in\name b$''}\}$.  As $q_2\in
\AP$, we have $\name b^+\in\name \FF^0[\Hh^0]$.
For each $n\in \name b^+ [H^\delta]$ let
$$\name A^1_n[\Hh^0]=: \{v\in \Gamma^{\epsilon}_{\name
k(n)}[\Hh^0]:
\hbox{ For no $\pp'_1\ge\pp_1$ in $\Pp\restrictedto \AA_1$ with
$\pp'_1\restrictedto \AA_0\in \Hh^0$,}$$
$$\pp'_1\forces_{\Pp\restrictedto \AA_1}
\hbox{``$n\in\name a$ and $v\notin \name A_n$.''}\}$$
Let $\name A^1_n[\Hh^0]=\emptyset$ if $n\notin \name b^+$.

Easily $(\name A^1_n:n<\omega)$ is $\name\FF^0$-slow.  Hence
in $V[\Hh^0]$ the sequence
$(\name B_n \setminus \name A^1_n:n<\omega)$ is not $(\name
\FF^0[\Hh^0])$-slow.  
We can compute the values of $\name B_n$ and $\name A^1_n$ in
$V[\Hh^0]$. 
So we can find $n\in \name b^+ [H^0]$ with
$\name B_n \setminus \name A^1_n\ne \emptyset$, and choose $v\in \name
B_n \setminus \name A^1_n$. 
Then there are $\pp_1'\in \Pp\restrictedto \AA_1/\Hh^0$,
$\pp_1'\ge \pp_1$, and $\pp_2'\in \Pp\restrictedto \AA_2/\Hh^0$, with
$\pp_2'\ge \pp_2$, so that: 
$$\pp'_1\forces\hbox{``$n\in \name a$ and
$v\not \in\name A_n$''.}$$ $$\pp_2'\forces
\hbox{``$n\in \name b$ and $\name x_\delta(n)=v$''}$$

Now $\pp\le\pp_1'\union \pp_2'\in \Pp\restrictedto \AA$ and
$\pp'_1\union\pp'_2$ forces ``$n\in\name a\intersect\name b\intersect\name c$''
(over $\Hh^0$), contradicting the choice of $\pp$.
This completes the proof of 3.9 (1).
\medskip
2: Let $[(\AA^{q_2} \setminus \alpha)\Union\{\sup\,\AA^{q_2}\}] =
\{\delta_i:i\le 
\gamma\}$ in increasing order.  Define inductively $r_i\in \AP$,
increasing in $i$, with $q_2\restrictedto (\AA\intersect\delta_i)\le
r_i$, $\dom r_i\includedin \delta_i$, $r_0=q_1$; then let $r=r_\gamma$.

At successor stages $i=j+1$ we
apply  3.9 (1) to $q_2\restrictedto (\AA^{q_2}\intersect \delta_j)$, $r_j$, $q_2\restrictedto [\AA^{q_2}\intersect(\delta_j+1)]$.

If  $i$ is a limit of uncountable cofinality,
we just take unions:
$$\AA^{r_i}=\Union_{\zeta<i} \AA^{r_\zeta}; \name
\FF^{r_i}=\Union_{\zeta<i}\name\FF^{r_\zeta};
{\fakebf\epsilon}^{r_i}=\Union_{\zeta<i}{\bf\epsilon}^{r_\zeta}; {\bf\name k}^{r_i}=\Union_{\zeta<i}{\bf\name k}^{r_\zeta};$$
while if  $i$ is a limit of cofinality $\aleph_0$,
we have actually to extend
$\Union_{\zeta<i}\name\FF^{r_\zeta
}$ to a $\Pp\restrictedto \AA^{r_i}$-name of an
ultrafilter in $V^{\Pp\restrictedto \AA^{r_i}}$.
However, in $V^{\Pp\restrictedto \AA^{r_i}}$, $\Union_{\zeta<i}\name\FF^{r_\zeta}$
is interpreted as a filter including all cofinite subsets of $\omega$, hence can be completed
to an ultrafilter.   \qed

\proclaim{3.10\\ Claim}

\item{1.}If $q_i$ $(i<\delta)$ is an increasing sequence of members of $\AP$, with
$\delta<\aleph_2$, then for some $q\in \AP$, $q\ge q_i$ for all
$i<\delta$.

\item{2.} If $q_1,q_2\in\AP$,
$\alpha<\aleph_3$,
$q_2\restrictedto \alpha \le q_1$, and $\dom q_1\intersect \dom
q_2=\dom q_1\intersect \alpha$,
then there is
$r\ge q_1,q_2$ in $\AP$.

\endproclaim

\proof:

1: We may suppose $\delta=\aleph_0$ or $\aleph_1$.
Let $\AA=:\Union_i \AA^{q_i}$ be enumerated in increasing order as
$\{\alpha_j:j<\gamma\}$ for the appropriate $\gamma$, and set
$\alpha_\gamma=\sup\,\AA$.  We define an increasing sequence of members
$r_j$ of $\AP$  for $j\le \gamma$ by induction on $j$ so that:$$\AA^{r_j}=\{\alpha_\zeta:\zeta<j\};$$
$$q_i\restrictedto\alpha_j\le r_j {\rm\ for\ all\ } i<\delta.$$

In all cases we proceed as in the proof of Claim 3.9.  The only difference is
that we deal with several $q_i$, but as they are
linearly ordered there is no difficulty.

2: This is proved similarly to part (1): 
let $\gamma=\sup(\dom q_1\cup \dom q_2$). 
 Choose by induction on $\beta\in(\dom q_1\cup \dom q_2\cup\{\gamma\}
)\setminus\alpha$ an upper bound $r_\beta$ of $q_1\uhr\beta$ and $q_2\uhr \beta$, increasing with $\beta$, with $\dom r_\beta=\beta\cap (\dom q_1\cup
\dom q_2)$. The successor step is by 3.9(i). 
The limit is easy too. 
Note: if $\dom q_1/E$ has only 
finitely many classes, when $\beta_1\,E\, \beta_2$ iff 
$\bigwedge_{\gamma\in \dom q_2}[\gamma<\beta_1\Leftrightarrow \gamma <\beta_2]$, then  
3.9(ii) suffices. 
 \qed

\proclaim{3.11\\Proof of Theorem 3.1: The construction}
\endproclaim

We define an increasing sequence $G^\alpha\includedin\{q \in\AP:\AA^q
\includedin \alpha \}$ of $\aleph_2$-directed sets increasing in $\alpha$,
and a set of at most $\aleph_2$ ``commitments'' which $G^\alpha$ will
meet.  In particular we require
that $\forall\beta<\alpha\,\exists q \in G^\alpha\ (\beta\in
\AA^q)$, and at each stage $\alpha$ we may make new commitments to
``enter some collection of dense sets'' -- in set theoretic
terminology -- or equivalently, to ``omit some type'' -- in model
theoretic terms. We make use of
$\diamond_{\{\delta<\aleph_3:\cof\delta=\aleph_2}\}$ to choose the
commitments.  The combinatorics involved in meeting the commitments are
treated in [ShLH 162], and are reviewed in \S A3 of the Appendix.  Our
summary of the construction in the present section will be less
formal.

At a stage $\delta<\aleph_3$ with $\cof\,\delta=\aleph_2$, we will
``guess'' $\Pp$-names $\name k^1_\delta, \name k^2_\delta, \name F_\delta$,
a condition $\pp^\delta\in\Pp\restrictedto\delta$ and a parameter
$\epsilon_\delta\in\{1,2\}$, explained in connection with (4) below,
and attempt to ``kill'' the possibility that  $\pp^\delta$ forces:
\medskip

\item{}``$\name
F_\delta:\prod_n \Gamma^1_{\name k^1_\delta(n)} \to \prod_n
\Gamma^2_{\name k^2_\delta(n)}$ induces a map which can be extended to
an  isomorphism:
$$\uprodname
{\Gamma^1_{\name k^1_\delta(n)}}/\name\FF\iso \uprodname {\Gamma^2_{\name
k^2_\delta(n)}}/\name\FF\hbox{''.}$$

\noindent
(Here we have taken $\epsilon_\delta=1$; otherwise the roles of 1 and
2 in this -- and in all that follows -- must be reversed.)

We will refer to the genericity game of [ShHL 162], as described in \S A3
of the Appendix.  In that game the Ghibellines can accomplish the
following.  For $\delta<\aleph_3$, they determine a set of compatible
approximations $G^\delta$ which together will determine an ultrafilter
$\FF\restrictedto \delta$ on $\omega$ in $V^{\Pp\restrictedto
\delta}$ (specifically, $G^\alpha$ is a subset of $\{r \in \AA\PP:
\Dom r \subseteq \alpha\}$ which is directed, increasing in $\alpha$).
 The Guelfs set them tasks which ensure that the ultrafilter
$\FF$ which is gradually built up by the Ghibellines has all the
desired properties.

Let $\FF_0$ be a fixed nonprincipal ultrafilter on $\omega$, in the
ground model and without loss of generality there is $q\in G\up 0$
with $\FF\up q=\FF_0$.
For $\delta <\aleph_3$ of cofinality $\aleph_2$,
let $q^*_\delta$ be an  approximation
$(\{\delta\},\name\FF^\delta,(\epsilon_\delta), (\name
k^{\epsilon_\delta}_\delta))$, where $\name\FF^\delta$ is the
$\Pp\restrictedto \{\delta\}$-name of some ultrafilter on $\omega$
extending $\FF_0$ such that
\medskip
\+$(1)$&$\{n:\name x_\delta(n)\in
\Gamma^{\epsilon_\delta}_{\name k^{\epsilon_\delta}_\delta(n)}\}\in
\name\FF^\delta$;\cr
\medskip
\+$(2)$&$\{n:\name x_\delta(n)\notin A_n\}\in\name\FF^\delta$ for any
$(\FF_0,\name k^{\epsilon_\delta}_\delta,\epsilon_\delta)$-slow
sequence $(A_n)$\cr
\+&in the universe $V$\cr
\medskip

The Ghibellines will be required (by the Guelfs) to put $q^*_\delta$
in $G^{\delta+1}$.  The Ghibellines
are  also obliged to make commitments of the following form, which
should then be
respected throughout the rest of the construction. (These commitments
involve a parameter $\alpha>\delta$ to be controlled by the
Ghibellines as the play progresses: of course these commitments have
to satisfy density requirements.)
$$\eqalign{&\hbox{For every $\alpha>\delta$, every $q\in G^\alpha$ with
$\delta\in \dom q$,}\cr
&\hbox{every $k_\delta^{1-\varepsilon_\delta}(n)$ 
(really a ($\Pp\restrictedto \delta)$-name)}\cr
&\hbox{and every $(\Pp\restrictedto \AA^\delta)$-name $\name z$}\cr
 & \hbox{of a member of
 $\prod_n \Gamma_{k_\delta(n)}^{1-\epsilon_\delta}:$}\cr}
 \leqno(*)^\delta_{q^*,\name z^*}:$$
\itemitem{}\\{\it if}
$(q,\name z)\iso (q^*,\name z^*)$ over $\delta+1$, {\it then} there
will be some $r$ in $G^\alpha$, some
\itemitem{}\\$\pp'\in \Pp\restrictedto \AA^r$, and some
$\Pp\restrictedto(\AA^r\intersect \delta)$-name
$\name x$ of a member of
$\prod_n\Gamma^{\epsilon_\delta}_{\name k^{\epsilon_\delta}_\delta(n)}$,
\itemitem{}\\with $r\ge q$, $\pp'\ge \pp^\delta$,
$\name
F_\delta(\name x)$ is a $\Pp\restrictedto
(\AA^r\intersect\delta)$-name, and:
$$\eqalign{\pp'\forces_{\Pp\restrictedto \AA^r}\,\hbox{``}
\{n:\Gamma ^{\epsilon_\delta}_{k^{\epsilon_\delta}_\delta(n)}&\models
R(\name x(n),\name x_\delta(n))\iff\cr
&\Gamma^{1-\epsilon_\delta}_{k^{1-\epsilon_\delta}_\delta(n)}\models
\neg R(\name
F_\delta(\name x)(n),\name z(n))\}\in\name\FF^r\hbox{''}\cr}\leqno(\dagger)$$

\noindent There is such a commitment for each $q^*,\name z^*$ with
$q^*_\delta\le q^*\in \AP$, $q^*\restrictedto \delta\in G^\delta$, and
$\name z^* $ a $(\Pp\restrictedto \AA^{q^*})$-name of a member of
$\prod_n \Gamma^2_{\name k^2_\delta(n)}$.  So apparently we are making
$\aleph_3$ commitments, which is not feasible, but as we are using
isomorphism types this amounts to only $2^{\aleph_1}=\aleph_2$
commitments, and this {\sl is} feasible.  This is formalized in  \S A3.6 in the 
Appendix.

These commitments can only be met when the corresponding set of
approximations is dense, but on the other hand we
have a stationary set $\delta$ of opportunities to meet such a
commitment,
and we will show that for any candidate $\name F$ for an isomorphism,
either we kill it off as outlined above (by making it obvious that
$\name F(\name x_\delta)$ cannot be defined), or else
-- after
failing to do this on a stationary set -- that $\name F$ must be quite
special (somewhat definable) and hence even more easily dealt with, as
will be seen  in detail in the next few sections.

After we have obtained $G^\alpha$ for all $\alpha$, we will let $\name
\FF^\alpha$ be  $\Union\{\name \FF^q :q \in G^\alpha\}$  (that is, the
appropriate  $(\Pp\restrictedto\alpha)$-name of a uniform ultrafilter
on $\omega$). Letting $G=:G^{\aleph_3}=:\Union_\alpha G^\alpha$, also
$\name \FF=\name \FF^{\aleph_3}$ is defined.
\medskip

\proclaim{3.12\\Proof of Theorem 3.1: The heart of the matter}
 \endproclaim

Now suppose toward a contradiction that after $\name \FF$ has been
constructed in this way, there are 
$\Pp$-names $\name F$, $\name k^1$, $\name k^2$, and a
condition $\pp \in \Pp$ such that: 
$$\pp\forces_{\Pp}``\name F \hbox{ is a function from } 
\textstyle{\prod_n} \Gamma^1_{\name k^1(n)} 
\hbox{ onto }\textstyle{\prod_n}
\Gamma^2_{\name k^2(n)}\hbox{ which induces an isomorphism }\leqno(3)$$
\item{}\\of the corresponding ultraproducts with respect to 
$ \name \FF$''.  
\medskip

Actually, we will want to assume in addition that
$\pp$ forces:
$$\hbox{``$\{n:\norm{\Gamma^1_{\name k^1(n)}}>\norm{\Gamma^2_{\name k^2(n)}}\}\in
\name\FF$,''}\leqno(4)$$

\noindent which could force us to increase $p$ 
and to switch the roles of 1 and 2 in all
that follows; this is why we have carried along a parameter $\epsilon$
in our definition of $\AP$.  

We will say that a set
$\AA\includedin\aleph_3$ is  $(\name F,\name k^1,\name
k^2,\pp)$-{\sl closed} if: 
\item{i.}$\name k^1,\name k^2$
are $(\Pp\restrictedto \AA)$-names; $\name F\restrictedto \AA$ is a
$(\Pp\restrictedto \AA)$-name;
\item{ii.}$\pp\forces_{\Pp\restrictedto \AA}$:
\itemitem{``}$\name F\restrictedto \AA$
is a function from $\prod_n \Gamma^1_{\name k^1(n)}$
onto $\prod_n \Gamma^2_{\name k^2(n)}$ which
(interpreted in $\Pp\restrictedto \AA$) induces an isomorphism from
$\prod_n \Gamma^1_{\name k^1(n)}/(\name\FF\restrictedto \AA)$ onto
$\prod_n \Gamma^2_{\name k^2(n)}/(\name\FF\restrictedto \AA)$''.
\item{iii.}$\pp\forces_{\Pp\restrictedto \AA}$:
``$\{n:\norm{\Gamma^1_{\name k^1(n)}}>\norm{\Gamma^2_{\name k^2(n)}}\}\in
\name\FF\restrictedto \AA$.''

\noindent Properly speaking, the only actual closure condition here is
clause (ii).  Note that the condition in (iii) can be strengthened to:
$$\hbox{``$\left\{n:\norm{\Gamma^1_{\name k^1(n)}}>
{\norm{\Gamma^2_{\name k^2(n)}}}^5 
\right\}\in \name\FF\restrictedto \AA$,''}$$ 
by the choice of
the sequences $(\Gamma^i_n)$ $(i=1,2)$.

 Let $C$ be 
$\{\delta<\aleph_3:\cof (\delta) = \aleph_2,\, \delta $ is
$(\name F,\name k^1,\name k^2,\pp)$-closed$\}$.
Clearly the set $C$
is unbounded and is closed under $\aleph_2$-limits.  By our
construction, for a stationary subset $S_C$ of $C$ we may suppose
that for $\delta\in S_C$:
$\name F_\delta=\name F\restrictedto \delta$, $\pp^\delta=\pp$,
$\epsilon_\delta=1$, $\name k_\delta=\name k^1$, and that $\delta$ was $(\name F,\name
k^1,\name k^2,\pp)$-closed.
So $q^*_\delta\in G^{\delta+1}$, and we
can find $q\in G$ such that $\name z=:\name F(\name x_\delta)$ is
a $(\Pp\restrictedto \AA^{q})$-name,
$\delta\in \AA^{q}$.

At stage $\delta$ in the
construction, the Ghibellines had tried to make the commitment
$(*)^\delta_{q^*,\name z^*}$, with $(q^*,\name z^*)=(q,\name z)$. They
later failed to meet this commitment, since otherwise there would be
some $r\ge q$ in $G$,
some $\pp'\ge \pp$ in $\Pp\restrictedto \AA^r$, and some
$[\Pp\restrictedto (\AA^r\intersect\delta)]$-name of a member $\name
x$ of $\Gamma^1_{\name k^1(n)}$, for which $(\dagger)$ holds:
$$\pp'\forces_{\Pp\restrictedto \AA^r}
\hbox{``$\{n:[\Gamma^1_{k^1_\delta}\models R(\name x(n),\name
x_\delta(n))\iff \Gamma^2_{k^2_\delta}\models\neg R(\name F_\delta(\name x)(n),\name z(n))]\}\in
\FF^r$''.}$$
 and  $\name z$ is $\name F(\name x_\delta)$. 
But $\pp$ forced $\name F$ to induce an
isomorphism,
 so we have a contradiction.

The failure to make the commitment
$(*)^\delta_{q,\name z}$, implies a failure of density,
which means that for
some $(q',\name z')\iso (q,\name z)$ over $\delta+1$ -- and hence also
for $(q,\name z)$ -- taking
$q_0=q\restrictedto \delta$, we will have:

\item{(i)} $\delta$ is $(\name F,\name
k^1,\name k^2,\pp)$-closed.
\item{(ii)} $\pp\in\Pp\restrictedto \AA^{q_0}$, $\delta\in\AA^q$,
$\epsilon^q=1$, $\name k^q_\delta=\name k^1$, $\name
F_\delta=\name F\restrictedto \delta$;
\item{(iii)} $\name z$ is a $(\Pp\restrictedto \AA^q)$-name for a
 member of $\prod_n \Gamma^2_{k^2_\delta(n)}$;

\item{(iv)} For all $r\ge q$ in $\AP$ such that
$r\restrictedto \delta\in G^\delta$, and
$\name x$ a $(\Pp\restrictedto \AA^{r\restrictedto \delta})$-name,
with $\name y=: \name F(\name x)$  a
$(\Pp\restrictedto \AA^{r\restrictedto \delta})$-name, we have:
\itemitem{$(*)_{\name x,\name y}$\ }
$\pp\forces$
	 ``The set $\{n:\Gamma^1_{\name k^1(n)}\satisfies
	 R(\name x(n),\name x_\delta(n))\hbox{ iff }\Gamma^2_{\name
	 k^2(n)}\satisfies
	 R(\name y(n),\name z(n))\}$
\itemitem{}\\\kern 5pt is in $\name \FF^r$''.
\nl
(Note: another possibility of failure, $q\notin G^\alpha$, is ruled out 
by the choice of $q$).

Now we analyze the meaning of $(*)_{\name x,\name y}$.  Consider the
following property of $(\Pp\restrictedto \delta)$-names $\name x,\name y$
for a fixed choice of $\delta\in C$, $q\in \AP$ with $\delta\in
\AA^{q}$, and $\name z$ a $(\Pp\restrictedto \AA^q)$-name.

\medskip
\+$(**)_{\name x,\name y}$&For all $r\ge q$ in $\AP$ such that
$r\restrictedto \delta\in  G^\delta$ and $\name x$, $\name y$ are
$(\Pp\restrictedto \AA^{r\restrictedto \delta})$-names,\cr
\+&$(*)_{\name x,\name y}$ holds.\cr
\medskip
 We explore the meaning of this property when $\name y$ is not
 necessarily $\name F(\name x)$.

Clearly, 
\item{$(\otimes_1)$} If $\name x$ is a $(\Pp\restrictedto \delta)$-name, 
$\name y=\name F(\name x)$, then 
$(**)_{\name x, \name y}$.

To simplify the analysis, let $\Hh$ be generic for $\Pp\restrictedto \delta$. 
 Let $\name x$ be a $\Pp\restrictedto \delta$-name of a real,
 $\AA\includedin \delta$. 
 We say $\name x$ is
 {\sl unrestricted} for $(\Hh,\AA,\name k^1)$ if:
 \itemitem {}There is no $(\name \FF\restrictedto \AA,\name k^1)$-slow 
 sequence $(\name B_n)_{n<\omega}$ in $V[\Hh\restrictedto \AA]$ such
 that:

 \itemitem{} $\{n:\name x(n)\in \Gamma^1_{\name k^1(n)}\backslash
 \name B_n\}\equiv \emptyset 
 \\\mod \name\FF^\delta[\Hh]$.

 Observe that if $\sup\, \AA<\gamma<\delta$ and $\name
 k_\gamma=\name k^1$, then the Cohen real $\name x_\gamma$ is forced (in 
 $\Pp\restrictedto \delta)$ to be unrestricted for
 $(\Hh,\AA,\name k^1)$.

 \bigskip

 \proclaim{3.12A Claim}
\endproclaim
If $\name x^1,\name x^2$ are $(\Pp\restrictedto
 \delta)$-names of functions in $\prod_n\Gamma^1_{\name k^1(n)}$, 
 \ $\name y$ is a
 $(\Pp\restrictedto \delta)$-name of a member of
 $\prod_n\Gamma^2_{\name k^2(n)}$, 
 and both pairs $(\name x^1,\name y)$ and $(\name x^2,\name y)$ satisfy the
 condition (**)  above, then:

 $$\pp\forces_{\Pp\restrictedto \delta}\hbox{``$\name x^1=\name x^2\mod
 \name \FF\restrictedto \delta[\Hh]$ or
 both are restricted for $(\Hh,\AA^{q_0},\name k^1)$.''}\leqno({\rm Clm})$$

 We will give the proof of this, which contains one of the main
 combinatorial points,  in paragraph 3.13.  For the present we
 continue with the proof of the theorem.
 We first record a consequence of the claim.

$$\vcenter{\hbox{If $\name x,\name y$ are $(\Pp\restrictedto
\delta)$-names with 
 $\name x$ forced by $\Pp\restrictedto \delta$ to be unrestricted for 
$(\Hh,\AA^{q_0},\name k^1)$,\strut}
\hbox{\strut and the pair $(\name x,\name y)$
satisfies $(**)_{\name x,\name y}$, then
$p\forces_{\Pp\restrictedto \delta}$ ``$\name F(\name x)=\name y$
 mod $\name \FF^\delta$''}}\leqno(\otimes_2)$$

 Indeed, if
 $\Hh\includedin \Pp\restrictedto
 \delta$ is generic over $V$,
 and $\name F(\name x)[\Hh]=\name y_1[\Hh]\ne \name y[\Hh]$ mod $\name
 \FF^\delta$, then since $\name F$ is {\sl onto} (in $V[\Hh]$, as $\delta$ is
 $(\name F,\name k^1,\name k^2,\pp)$-closed),
 there is a $(\Pp\restrictedto \delta)$-name $\name x'$ with
 $\name F(\name x')[\Hh]=\name y[H]$, so $\name x'[\Hh]\ne\name x[\Hh]$
 mod $\name \FF^\delta$.
 Now $\name x, \name x', \name y$ contradict (Clm).  Thus $(\name F)$ holds.
 As ($\otimes_1)+(\otimes_2)$
 holds for stationarily many $\delta$'s, it holds for $\delta=\aleph_3$
 (in the natural interpretation).

In what follows, we use the statements
 $(\otimes_1)+(\otimes_2)$ as a kind of
``definability'' condition on $\name F$; but we deal with the current
concrete case, rather than seeking an abstract formulation of the situation.

Let 
$S =\{\gamma\in S_C:\name F(\name x_\gamma)$ is (forced by $p$ to be equal to) 
a
$[\Pp\restrictedto (\gamma+1)]$-name $\}$.  
We claim that $S$ is
stationary.  Let $C'\includedin \aleph_3$ be closed unbounded, and let
$\delta\in S_C$ be taken with $C'\intersect S_C$ unbounded below
$\delta$. Let $q\in G$ be chosen so that $\name F(\name x_\delta)$ is
a $(\Pp\restrictedto \AA^q)$-name, let $q_0=q\restrictedto \delta$, and $\gamma_0=\sup\,
\AA^{q_0}$. It suffices to check that for $\gamma_0<\gamma<\delta$
with $\gamma\in S_C$, we have $\gamma\in S$.
So let $r_1\in G^\delta$ be chosen so that $\name y_1=:\name F(\name
x_\gamma)$ is a $(\Pp\restrictedto \AA^{r_1})$-name.  It suffices to
show that $\name y_1$ is (forced by $p$ to be equal to) a $(\Pp\restrictedto
[\AA^{r_1}\intersect(\gamma+1)])$-name.  Otherwise, 
by a density requirement (Appendix, \S A3) we can find a 1-1 order preserving
function $h$ with domain $\AA^{r_1}$, $h$ is the identity on
$\AA^{r_1}\intersect(\gamma+1)$,
$h(\min(\AA^{r_1}\backslash(\gamma+1)))>\sup\,\AA^{r_1}$, with $r_2=:
h(r_1)$ in $G^\delta$.  Let $\name y_2=h(\name y_1)$.  Then
$(**)_{\name x_{\gamma},\name y_i}$ holds for $i=1,2$, so
$p\forces_{\Pp\restrictedto \delta}$ ``$\name y_1=\name y_2
\mod\name\FF^\delta$'',   
but by 3.14 below we can ensure that this is not the case (by making
additional commitments, cf. \S A3).

Now for $\gamma\in S$ let $q_\gamma\in G^{\gamma+1}$ be chosen so that
$\name z_\gamma=\name F(\name x_\gamma)$ is a $(\Pp\restrictedto
\AA^{q_\gamma})$-name, and let $\hat
\gamma=\sup(\AA^{q_\gamma}\intersect\gamma)$. By Fodor's lemma we can
shrink $S$ so that $\hat \gamma$ and $\AA_0=\AA^{q_\gamma}\intersect
 \hat\gamma$ and $q_\gamma\restrictedto \hat \gamma$ 
are constant for $\gamma\in S$.  Now choose $\delta_1<\delta_2$ in
$S$, and let $q_i=q_{\delta_i}$, 
$\AA_i=\AA^{q_i}$ for $i=1,2$, so $\AA_1=\AA^{q_1}=\AA_0\union
\{\delta_1\}$, $\AA_2=\AA^{q_2}=\AA_0\union\{\delta_2\}$; also let 
$\AA=:\AA_1\cup \AA_2$; we now let $q_i\restrictedto \hat \gamma$ be called $q_0$.  
Let $\name \FF^i=\name \FF^{q_i}$, and set
$$\name d =:\{n:\Gamma^1_{\name k^1(n)}\satisfies
R(\name x_{\delta_1}(n), \name x_{\delta_2}(n))\iff
\Gamma^2_{\name k^2(n)}\satisfies
\neg R(\name z_{\delta_1}(n), \name z_{\delta_2}(n))\}.\leqno(\name d)$$

We want to find $r\in \AP$ with $\AA^r=\AA$ so that
$r\ge q_1,q_2$, and 
$\pp\forces \hbox{``$\name d\in \name \FF^r$''}$. 
This will then mean that $\name F$ could have been ``killed'', after all,
and will complete the argument.

Suppose this is not possible, and thus as in 3.9 (1) for some $\pp'\ge\pp$ in
$\Pp\restrictedto \AA$, if $\pp'_i=\pp'\restrictedto \AA_i$ for
$i=0,1,2$, we have:

\item{} a \PpA1-name $\name a$ of a member of $\name \FF^1$;
\item{} a \PpA2-name $\name b$ of a member of $\name \FF^2$; and
\item{} a $\Pp$-name $\name c
=:\{n: \name x_{\delta_2}(n)\in\Gamma^1_{\name k^1(n)}\backslash\name A_n\}$
associated with a \PpA1-name $(\name A_n)_{n<\omega}$ of an $(\name
\FF^1,\name k^1)$-slow sequence; with
\itemitem{}$\pp'\forces_{\Pp\restrictedto \AA}$
``$\name a \intersect\name b\intersect\name c\intersect \name d
= \emptyset$''

\noindent We shall get a contradiction.
 Let $H^0\subseteq\Pp\restrictedto\unA^0$ be generic over $V$.

We define for every $n$ the following ($\Pp\restrictedto\unA^0$)-names:
$$\eqalign{\can_n^1[H^0]{}=\{(u,v)\in & \Gamma^1_{\nname k^1(n)}\times
\Gamma^2_{\nname k^2(n)}:
 \hbox{For some $\pp''_1\in
\Pp\restrictedto \AA_1$ with $\pp''_1\ge \pp'_1$ and
$\pp''_1\restrictedto \AA_0\in \Hh^0$,}\cr
 &\pp''_1\forces_{\Pp\restrictedto \AA_1/\hbox{$\Hh^0$}}
``\hbox{$[\nname
z_{\delta_1}(n)=u$, $u\notin \nname A_n$, $n\in \nname a$
 and $\nname z_{\delta_1}(n)=v$]''}\}\cr}$$
$$\eqalign{
\can_n^2[H^0]{}=\{(u,v)\in\Gamma^1_{\nname k^1(n)}\times\Gamma^2_{k^2(n)}:
& \hbox{For some $\pp''_2\in
\Pp\restrictedto \AA_2$ with $\pp''_2\ge \pp'_2$ and
$\pp''_2\restrictedto \AA_0\in \Hh^0$,}\cr
&\pp''_2\forces_{\Pp\restrictedto \AA_2/\hbox{$\Hh^0$}}
\hbox{``$[\nname
z_{\delta_2}(n)=u,\;n\in \nname b$ and
$\nname z_{\delta_2}(n)=v]$''}\}\cr}$$
and for
$i=1,2$ and $u\in \Gamma^1_{\name k^1(n)}$ we let
$$\can_n^i{}(u)=:\{v\in\Gamma^2_{\name k^2(n)}:
\;(u,v)\in\can_n^i{}\}$$
$$\name A^i_n=:\{u:(\exists v) (u,v)\in \can_n^i{}\}$$

Now in $V[\Hh^0]$, $(\name A^i_n:n<\omega)$ is not $(\name \FF^i,\name k^1)$-slow, and thus
the set: $$\{n:|\name A^1_n|/\norm{\Gamma^1_{\name k^1(n)}},
|\name A^2_n|/\norm{\Gamma^1_{\name k^1(n)}} \hbox{ are greater than}\;
\norm{\Gamma^1_{\name k^1(n)}}^{-1/2}\}$$ belongs to $\name \FF^0[\Hh]$.
Choose any such $n$, and by finite combinatorics we shall derive a
contradiction.  Remember that we have assumed without loss of
generality that $\norm{\Gamma^1_{\name k^1(n)}}>{ \norm {
\Gamma^2_{\name k^2(n)}
} }^5$
for a large set of $n$ modulo $\name \FF\restrictedto
\AA_0$, so wlog our $n$ satisfies this, too.
 Let $g_i:\name A^i_n\to \Gamma^2_{\name k^2(n)}$ be such that
$g_i(v)\in
\can_n^i{}(v)$.  Now
$|\rg(g_i)|\le \norm{\Gamma^2_{\name k^2(n)}}$, so there are
$b_1,b_2\in \Gamma^2_{\name k^2(n)}$ such that for $i=1,2$:
$$|g^{-1}_i(b_i)|\ge
|\name A^i_n|/\norm{\Gamma^2_{\name k^2(n)}}> \norm{\Gamma^1_{\name
k^1(n)}}^{1/5}.$$
Now
by 3.4(ii)
we find $a_i,a'_i\in g^{-1}(b_i)$ for $i=1,2$
with $\Gamma^1_{\name k^1(n)} \satisfies
R(a_1,a_2)\&\neg R(a'_1,a'_2)$.  As either $\Gamma^1_{\name k^1(n)}\satisfies
R(b_1,b_2)$ or $\Gamma^1_{\name k^1(n)}\satisfies \neg R(b_1,b_2)$, we can show
that it is not forced by $\pp'$  that $n\not\in\name a\intersect\name
b\intersect\name c\intersect\name d$, a contradiction.  \qed

\proclaim{3.13\\Proof of the Claim 3.12A from 3.12}
\endproclaim

We first recall the situation.  We had:

\item{(i)} $\delta$ is $(\name F,\name
k^1,\name k^2,\pp)$-closed; $q_0=q\restrictedto \delta$;
\item{(ii)} $\pp\in\Pp\restrictedto \AA^{q_0}$, $\delta\in\AA^q$,
$\epsilon^q=1$, $\name k^q_\delta=\name k^1$, $\name
F_\delta=\name F\restrictedto \delta$;
\item{(iii)} $\name z$ is a $(\Pp\restrictedto \AA^q)$-name for a real;

\item{(iv)} For all $r\ge q$ in $\AP$ such that
$r\restrictedto \delta\in G^\delta$, and
$\name x$ a $(\Pp\restrictedto \AA^{r\restrictedto \delta})$-name,
with $\name y=: \name F(\name x)$  a
$(\Pp\restrictedto \AA^{r\restrictedto \delta})$-name, we have:
\itemitem{$(*)_{\name x,\name y}$\ }
$\pp\forces$
	``The set $\{n:\Gamma^1_{\name k^1(n)}\satisfies
	R(\name x(n),\name x_\delta(n))\hbox{ iff }\Gamma^2_{\name
	k^2(n)}\satisfies
	R(\name y(n),\name z(n))\}$
\itemitem{}\\\\is in $\name \FF^r$''.

We defined the property $(**)_{\name x, \name y}$ as follows:

\medskip
\+$(**)_{\name x,\name y}$&For all $r\ge q$ in $\AP$ such that
$r\restrictedto \delta\in  G^\delta$ and $\name x$, $\name y$ are\cr
\+&$(\Pp\restrictedto \AA^{r\restrictedto \delta})$-names, $(*)_{\name
x,\name y}$ holds.\cr
\medskip

\proclaim{Claim}

If $\name x^1,\name x^2$ are $(\Pp\restrictedto
\delta)$-names of functions in $\prod_n\Gamma^1_{\name k^1(n)}$,
$\name y$ is a
$(\Pp\restrictedto \delta)$-name of a member of
$\prod_n\Gamma^2_{\name k^2(n)}$,
 and both pairs $(\name x^1,\name y)$ and $(\name x^2,\name y)$ satisfy the
condition $(**)_{\name x,\name y}$  above, then:

\centerline{$\pp\forces_{\Pp\restrictedto \delta}$``$\name x^1=\name x^2\mod
\name \FF\restrictedto \delta[\Hh]$ or
both are restricted for $(\Hh,\AA^{q_0},\name k^1)$.''}
\endproclaim

\proof:
Suppose that
$\pp\le \hat\pp\in \Pp\restrictedto \delta$
and $\hat\pp$ forces
the contrary; so without loss of generality
$$ \hat\pp\forces 
\hbox{``$\name x^1\ne \name x^2 \mod \name \FF\restrictedto \delta [\name
\Hh]$''};\leqno(5)$$
$$ \hat\pp\forces 
\hbox{``$\name x^1\hbox{ is unrestricted for } (\Hh,\AA^{q_0},\name k^1)$.''}\leqno(6)$$
Choose any $q_1\ge q_0$ with $q_1\in G^\delta$ so that
$\name x^1,\name x^2, \name y$ are $\Pp\restrictedto \AA^{q_1}$-names.
Now we will construct $r \ge q_1, q\restrictedto (\delta+1)$, with $r
$ in $\AP$ and
$\AA^r =\AA^{q_1}\union\{\delta\}$, so that:
$$\hat\pp\forces\hbox{``$\{n:\Gamma^1_{\name
k^1(n)}\satisfies$`$R(\name x^1(n),\name x_\delta(n))\iff \neg R(\name
x^2(n),\name x_\delta(n))]\}\in \name \FF^r $.''}\leqno(7)$$
By 3.9(2) we can also find $r'\ge r,q$, and then (7)
contradicts $(**)_{\name x^1,\name y}\and(**)_{\name x^2,\name y}$.
Thus to complete the proof of our claim, it suffices to find $r $.

This is the sort of problem considered in 3.9(1), with an additional
set required to be in $\name \FF\restrictedto (\AA^{q_1}\union \{\delta\})$.
The $q_0,q_1$ under consideration here correspond to the $q_0,q_1$ of
3.9(1), and we
let $q_2$ be $q\restrictedto (\delta+1)$. Following the notation of
3.9(1),  set
$\name \FF^i=\name \FF^{q_i}$, $\AA_i=\AA^{q_i}$ for $i=0,1,2$, and
$\AA=\AA_1\union \{\delta\}=\AA_1\union\AA_2$.  We need to find $r \ge
q_1,q_2$ as in 3.9(1), with (7) holding.

Suppose on the contrary that $\hat \pp\le \pp'\in \Pp\restrictedto
\AA$ and $\pp'$ forces
``There is no $\name \FF$ as required''.  Then extending $\pp'$, we
may suppose that we have a $\Pp\restrictedto\AA_1$-name $\name a$ for
a member of $\name \FF^1$, a $\Pp\restrictedto \AA_2$-name $\name b$
for a member of $\name \FF^2$, a $\Pp\restrictedto \AA_1$-name for an
$(\name \FF^1,\name k^1)$-slow sequence $(\name A_n)$ (associated with a
power $d<\omega$ -- cf. 3.7), such that setting:
$$\name c=\{n:\name x_\delta(n)\in\Gamma^1_{\name
k^1(n)}\backslash\name A_n\}$$
$$\name d=\{n:\Gamma^1_{\name k^1(n)}\satisfies\hbox{``$R(\name
x^1(n),\name x_\delta(n))\iff\neg R(\name x^2(n),\name x_\delta(n))$''}\}$$
we have:
$$\pp'\forces_{\Pp\restrictedto \AA}
\hbox{``$\name a\intersect\name b\intersect\name
c\intersect\name d=\emptyset$''}$$

Let $\pp'_i=\pp'\restrictedto \AA_i$ for $i=0,1,2$, and take
$\Hh^0\includedin \Pp\restrictedto \AA_0$ generic over $V$.  Without
loss of generality,  for some natural number $d$: 
$$\pp'_1\forces\hbox{``$n\in \name a \implies \name
x^1(n)\ne\name x^2(n)$
and $|\name A_n|\le
\sqrt{\norm{\Gamma^1_{\name
k^1(n)}}}\cdot(\log \norm{\Gamma^1_{\name k^1(n)}})^d$ 
 (and $\name A_n\subseteq \Gamma^1_{\name k^1(n)})$.''}$$
We are
interested in $\name B_n[\Hh^0]$=:
$$\{v\in \Gamma^1_{\name k^1(n)}: \hbox{ for some $\pp''_2\ge \pp'_2
$ with $\pp''_2\restrictedto \AA_0\in \Hh^0$,
$\pp''_2\forces_{\Pp\restrictedto \AA_2}$:``$n\in \name b$ and $\name x_\delta(n)=v$''}\}$$
(which is a $(\Pp\restrictedto \AA_0)$-name).  Clearly the sequence $(\name
B_n)$ is not $(\name \FF^1,\name k^1)$-slow in $V[\Hh^0]$.

For each $n$ let us also consider the set
$\name Y_n[\Hh^0]=$:

$$ \{(A,v_1,v_2):
A \cup \{v_1, v_2\}\subseteq \Gamma^1_{\name k^1(n)},
 v_1\neq v_2,
\hbox{ and for some $\pp''_1$ with $\pp''_1\ge \pp'_1, \pp''_1\restrictedto \AA_0\in
\Hh^0,$}$$
$$ \hbox{$\pp''_1\forces\,$``$n\in \name a,\name A_n=A,\name x^1(n)=v_1,\name x^2(n)=v_2$.''}\}$$

For every $(A,v_1,v_2)\in \name Y_n$, we have:
$$|A|\le \sqrt{\norm{\Gamma^1_{\name
k^1(n)}}}\cdot(\log \norm{\Gamma^1_{\name k^1(n)}})^d\hbox{, and }
v_1\ne v_2\leqno(8)$$

As $\name x_\delta$ is unrestricted over $\AA_0$ in $V[\Hh^0]$,
for the $\name\FF^0$-majority of $n$ we have:
$$|\name B_n|\ge \sqrt{\norm{\Gamma^1_{\name k^1(n)}}}\cdot(\log
\norm{\Gamma^1_{\name k^1(n)}})^{d+2}
\leqno(9)$$
Now (by (6)), also for the $\name \FF^0$ majority of $n$ we have:
$$\name C_n=:\{v_1\in \Gamma^1_{\name k^1(n)}:
\hbox
{There are $A,v_2$ so that
$(A,v_1,v_2)\in \name Y_n$}\}
\leqno(10)$$
$$\hbox {has at least}\;
\sqrt{\norm{\Gamma^1_{\name k^1(n)}}} \hbox{ members}$$

Now it will suffice to find $n$, $v\in \name B_n$ and
$(A,v_1,v_2)\in \name Y_n$ so
that
$$\Gamma^1_{\name k^1(n)}\satisfies [R(v_1,v)\iff\neg R(v_2,v)] \and
v\notin A,\leqno(11)$$
as we can then choose
$\pp''_1\in \Pp\restrictedto \AA_1$,
$\pp''_2\in \Pp\restrictedto \AA_2$ with $\pp''_i\ge\pp'_i$,
$\pp''_i\restrictedto \AA_0\in \Hh^0$ for $i=1,2$, so that:
$$\pp''_1\forces\hbox{``$n\in \name a, \name A_n=A$,
 $\name x^1(n)=v_1$, $\name x^2(n)=v_2$'';\\ }
\pp''_2\forces\hbox{``$n\in \name b$ and $\name x_\delta(n)=v$''}$$
and hence
$\pp''_1\union \pp''_2\forces$ ``$n\in\name a\intersect\name b\intersect\name
c\intersect\name d$'', a contradiction.

So it remains to find $n,\;v$ and $(A,v_1,v_2)$.  For $n$ sufficiently
large satisfying (8-10), we can choose triples
$t_i=(A^i,v^i_1,v^i_2)\in \name Y_n$ for $i<5\log
\norm{\Gamma^1_{\name k^1(n)}}$ with all vertices $v^i_1$ distinct from
each other and from all $v^i_2$. By the
pseudorandomness of $\Gamma^1_{\name k^1(n)}$
(more specifically 3.4(iii)), the set
\item{}$\name S=\{v\in\Gamma^1_{\name k^1(n)}:\hbox{ For no } i<5\log
\norm{\Gamma^1_{\name k^1(n)}}$
do we have $R(v^i_1,v)\iff \neg R(v^i_2,v)\}$

\noindent has size at most $5 \log \norm{\Gamma^1_{\name k^1(n)}}$.
So if $\name S'=: \name S \union
\Union\{A^i:i<5\log\norm{\Gamma^1_{\name k^1(n)}}\}$, then we will have:
$|\name S'|\ll \sqrt{ {\norm {\Gamma^1_{\name k^1(n)}} } }
(\log \norm{\Gamma^1_{\name k^1(n)}})^{d+2}$, so there is $v\in \name
B_n\setminus \name S'$.  Since $v\notin \name S'$, for some $i$ (11) will hold
with $(A,v_1,v_2)=(A^i,v^i_1,v^i_2)$.\qed

{\bigskip\noindent \bf3.14\\The last detail}

\noindent The following was used in the proof of 3.12
(after 3.12A slightly before ($\name d$)).
\proclaim{Claim.}
Assume $q_2\restrictedto \beta\le q_1$, $\AA^{q_1}\includedin
\beta$. Let $q_0=q_2\restrictedto \beta$, and write $\AA_i$ for $\AA^{q_i}$, $\AA=\AA_1\union\AA_2$,
and $\FF^i$ for $\FF^{q_i}$.  Let $\pp\in\Pp\restrictedto \AA$ and
$\pp_i=\pp\restrictedto \AA_i$.
Then we can find $r$ with $\AA^r=\AA$ and
$r\ge q_1,q_2$,  so that for any $(\Pp\restrictedto \AA_i)$-names
$\name y_i$ ($i=1,2$) of members of $\prod_n\Gamma^2_{\name k^2(n)}$
if:
$$\pp_i\forces_{\Pp\restrictedto \AA_i}\hbox{``$\name y_i\ne\name y'
\mod \name \FF^i$''}$$
for $(i=1,2)$ and for all $(\Pp\restrictedto \AA_0)$-names $\name y'$,
then we have:
$$\pp\forces_{\Pp\restrictedto \AA}\hbox{``$\name y_1\ne\name y_2\\\mod\name\FF^r$''}
$$
\endproclaim

Hence $p\forces_{\Pp\restrictedto \AA} $ 
 ``if $\name y_i\neq \name y'\mod \name \FF^i$ 
for $i=1,2$ and $\name y'$ a ($\Pp\restrictedto \AA_0)$-name 
 then $\name y_1\neq \name y_2\,\mod \name \FF^r$''.

{\it Proof}: We use induction construction. Much
as
in the proof of 3.9, we must deal primarily with the case in
which $\AA^{q_2}=\AA^{q_0}\union\{\beta\}$.
Suppose toward a contradiction
that $\pp\le\pp'\in\Pp\restrictedto \AA$, and with
$\pp'_i=\pp'\restrictedto \AA_i$ for $i=0,1,2$ we have:
\item{i.}a \PpA1-name $\name a$ of a member of $\name \FF^1$;
\item{ii.}a \PpA2-name $\name b$ of a member of $\name \FF^2$;
\item{iii.}a $(\Pp\restrictedto \AA)$-name 
$\name c=\{n:\name x_\beta(n)\in \Gamma^2_{\name
k^2(n)}\backslash\name A_n\}$ associated with a
\PpA1-name $(\name A_n)_{n<\omega}$ of a $(\name \FF^1,\name
k^1)$-slow sequence; and

\item{iv.}a $(\Pp\restrictedto \AA)$-name $\name
d=\Intersection_{j=1}^N\name d_j$, for a finite intersection of sets of
the form $\name d_j=:\{n:\name
y^1_j(n)\ne\name y^2_j(n)\}$, with 
each $\name y^i_j$ a $\Pp\restrictedto \AA_i$-name of a member of
$\prod_n\Gamma^2_{\name k^2(n)}$, such that for each $i=1,2$ and
$j=1,\ldots,N$:
$$\pp_i\forces\hbox{``$\name y^i_j\ne\name y' \\\mod \name \FF^i$
for any
$(\Pp\restrictedto \AA_0)$-name $\name y'$ of a member of
$\textstyle{\prod_n}\Gamma^2_{\name k^2(n)}$.''}$$

and that $\pp'\forces\hbox{``$\name a\intersect\name b\intersect\name
c\intersect\name d=\emptyset$''.}$
Let $\Hh^0$ be generic over $V$,
$\pp_0\in \Hh^0$, and 
let us define in $V[\Hh^0]$:
$$\name \Delta^1_n[\Hh^0]=: \hbox{$\{(A,u_1,\ldots,u_N): A\includedin
\Gamma^2_{\name k^2(n)}$, $u_1,\ldots,u_N\in \Gamma^2_{\name
k^2(n)}$,}$$
$$\hbox{ and there is $\pp'_1\in \Pp\restrictedto \AA_1$,
$\pp'_1\ge\pp_1$, $\pp'_1\restrictedto \AA_0\in\Hh^0$, and }$$
$$\hbox{$\pp'_1\forces_{\Pp\restrictedto \AA_1}$ ``$\name A_n=A$,
$\name y^1_1(n)=u_1,\ldots,\name y^1_N(n)=u_N$, and $n\in \name a$''.}\}$$

$$\name \Delta^2_n[\Hh^0]=:\{(v_0,v_1,\ldots, v_N):
  \hbox{ all $v_j\in \Gamma^2_{\name
  k^2(n)}$ and there is }$$
$$\pp'_2\in \Pp\restrictedto \AA_2,\, \pp'\ge\pp_2,\,
\pp'_2\restrictedto \AA_0\in \Hh^0\hbox{ and}$$
$$\hbox{$\pp'_2\forces_{\Pp\restrictedto \AA_2}$ 
``$\name
x_\beta(n)=v_0,\name y^2_1(n)=v_1\ldots,\name y^2_N(n)=v_N$ and $n\in
\name b$''}\}$$

Without loss of generality, for some $d$,

\noindent $$\pp_1\forces:\;\hbox{``For $n\in
\name a$, $|\name A_n|\le \sqrt{\norm{\Gamma^2_{\name k^2(n)}}}\cdot (\log
\norm{\Gamma^2_{\name k^2(n)}})^d$.''}$$
Thus:
$$(A,u_1,\ldots,u_N)\in \name \Delta^1_n
\implies
|A|\le
\sqrt{\norm{\Gamma^2_{\name k^2(n)}}} \cdot(\log
\norm{\Gamma^2_{\name k^2(n)}})^d.\leqno(1)$$

By the assumption on
$\name y^1_1,\ldots,\name y^1_N$,
$$
\hbox
  {If
     $e<\omega$,
     $C_n\includedin \Gamma^2_{\name k^2(n)}$,
     $|\name C_n|\le e$,
     and
     $(\name C_n:n<\omega)\in V[\Hh^0]$,
     then }
  \leqno(2)
$$

$$\{n:
  \hbox{
    there is
       $(A,u_1,\ldots,u_N)\in \name\Delta^1_n$,
       $u_1,\ldots,u_N\notin \name C_n$
    }
  \} \in \name \FF^0.
$$

Hence without loss of generality:
$$\hbox{For $n\in\name a$, there are $(A,u^j_1,\ldots,u^j_N)\in
\name \Delta^1_n$, for $j\le N+1$, with}$$
$$\hbox{ The sets $\{u^j_1,\ldots,u^j_N\}$ (for $j\le N+1)$
pairwise disjoint.}\leqno(3)$$

As $q_2\in \AP$,
$$\hbox{If $(\name C_n:n<\omega)\in V[\Hh^0]$ is $(\name \FF^0,\name
k)$-slow then} \leqno(4)$$
$$\{n:
\hbox{There is $(v_0,v_1,\ldots,v_N)\in\name\Delta^2_n$
with $v_0\notin \name C_n$}\}\in\name \FF^0$$

Let $\name a^+=:\{n:\name \Delta^1_n\ne\emptyset,
\name \Delta^2_n\ne\emptyset$, moreover, $\name\Delta^1_n$ satisfies (3)$\}$
(a $\Pp\restrictedto\unA^0$-name of a member of $\FF^0)$.
So for $n\in \name a^+$, there are $(N+1)$-tuples
$(A^{n,j},u^{n,j}_1,\ldots,u^{n,j}_N)$ for $j\le N+1$ with the sets
$\{u^{n,j}_1,\ldots,u^{n,j}_N\}$ pairwise disjoint.  Let $\name C_n=
\Union_{j\le N}A^{n,j}$ for $n\in\name a^+$, $\name C_n=\emptyset$ for
$n\notin \name a^+$.  So $(\name C_n)_{n<\omega}\in V[\Hh^0]$ is
$(\name \FF^0,\name k^1)$-slow, hence for some $n\in \name a^+$, there is
$(v_0,v_1,\ldots,v_N)\in \name \Delta^2_n$, with $v_0
\notin\name C_n$.  Now for some $j\le N+1$ we
have $\And_{i=1}^N v_i\ne u^{n,j}_i$.
Choose
$\pp'_2\in\Pp\restrictedto \AA_2$, $\pp'_2\ge\pp_2$, with
$\pp'_2\restrictedto \AA_0\in\Hh^0$ and $\pp'_2\forces$ ``$n\in \name
b$,
$\name
x_\beta(n)=v_0,\And_{i=1}^N \name y^2_i(n)=v_i$''.  
Choose $\pp'_1\in
\Pp\restrictedto \AA_1$, $\pp'_1\ge\pp_1$, with $\pp'_1\restrictedto
\AA_0\in\Hh^0$ and $\pp'_1\forces$ ``$n\in \name a$, $\name
A_n=A^{n,j}$, and for all
$i=1,\ldots,N$ $\name y^1_i(n)=u^{n,j}_i$.''  
Now
$\pp'_1\union\pp'_2\forces$ ``$n\in\name a\intersect\name b\intersect
\name c\intersect\name d$'', a contradiction.

This finishes the case $\AA_2=\AA_1\union\{\beta\}$.  The general case
follows as in 3.9(2). At successors we apply the case just treated.
Limits of uncountable cofinality are handled by taking unions.  At limits of
cofinality $\omega$ we have to repeat the first argument with some variations;
we do not have to worry about $\name c$, so the fact that there are several
$\name x_\beta$ involved is not  a problem.  The problem in this case is of
course to extend the union of the ultrafilters constructed so far to an
ultrafilter in a slightly larger model of set theory, while retaining the main
property for new names $\name y^2_i$.\qed

\page

\proclaim{Appendix.\\Background material.}
\endproclaim

\proclaim{A1. Proper and $\alpha$-proper forcing}
\endproclaim

\proclaim{A1.1.\\Proper Forcing}
\endproclaim
Let $\PP=(P,\le)$ be a partially ordered set.  A cardinal $\lambda$ is
{\sl $\PP$-large} if the power set of $P$ is in $V_\lambda$ (the
universe of all sets of rank less than $\lambda$). With
$\PP$ fixed and $\lambda$ $\PP$-large, let $\VV_\lambda$ be the
structure $(V_\lambda;\in,P,\le)$.

\item{1.} For $\MM\elementary \VV_\lambda$ and $p\in P$, $p$ is $\MM$-{\sl
generic} iff for each name of an ordinal $\name \alpha$ with $\name
\alpha\in M$, $p\forces$ ``$\name \alpha\in \MM$''.

\item{2.} $\PP$ is {\sl proper} iff for all $\PP$-large $\lambda$ and all
countable elementary substructures $\MM$ of
$\VV_\lambda$ with $\PP\in \MM$, each $p\in \MM$ has an $\MM$-generic extension
in $P$.

\proclaim{A1.2.\\Axiom A}
\endproclaim
$\PP$ satisfies {\sl Axiom} A if there is a collection $\le_n$
($n=1,2,\ldots$) of
partial orderings on the set $P$
with $\le_1$ coinciding with the given ordering $\le$, and $\le_{n+1}$
finer than $\le_n$ for each $n$, satisfying the following two
conditions:

\item{1.} If $p_1\le_1 p_2\le_2\le p_3\le_3\ldots$ then there is some
$p\in P$ with $p_n\le_n p$ for all $n$;
\item{2.} For all $p\in P$, any name $\name \alpha$ of an ordinal, and
any $n$, there is a condition $q\in P$ with $p\le_n q$, and a
countable set $B$ of ordinals, such that $q\forces \name \alpha\in B$.

The forcings used in \S\S1,2 were seen to satisfy Axiom A, and the
following known result was then applied.

\proclaim{A1.3.\\Proposition}
If $\PP$ satisfies Axiom A then $\PP$ is
proper.
\endproclaim
\proof:
 Given a countable $\MM\elementary \VV_\lambda$ and $p\in P\intersect
M$, let $\name \alpha_n$ be a list of all ordinal names in $\MM$, and use
clause (2) of Axiom A to find $q_n,B_n\in \MM$ with $q_n\in P$, $B_n$ countable.
$p\le_1q_1\le_2q_2\le\ldots$ and $q_n\forces$ ``$\name \alpha_n\in
B_n$. Then use clause (1) to find $q\ge$ all $q_n$; this $q$ will be
$\MM$-generic.\qed

\proclaim{A1.4\\Countable Support Iteration}
\endproclaim
Our notation for iterated forcing is as follows.  $\name \QQ_\alpha$
is the name of the $\alpha$-th forcing in the iteration, and
$\PP_\alpha$ is the iteration up to stage $\alpha$. The sequence
$\PP_\alpha$ is called the iteration, and the $\name \QQ_\alpha$ are
called the factors. It is assumed that
$\name \QQ_\alpha$ is a $\PP_\alpha$-name for a partially ordered set
with minimum element $0$,
and that $\PP_{\alpha+1}$ is $\PP_\alpha*\name\QQ_\alpha$.  

In general
it is necessary to impose some further conditions at limit ordinals.
We will be concerned exclusively with countable support iteration: at
a limit ordinal $\delta$, $\PP_\delta$ consists of $\delta$-sequences
$p$ such that $p\restrictedto \alpha\in \PP_\alpha$ for
$\alpha<\delta$, and $\forces_{\PP_\alpha}p(\alpha)=0$ for all but
countably many $\alpha <\delta$. 

\proclaim{A1.5\\Proposition}
Let $\PP_\alpha$ be a countable support
iteration of length $\lambda$ with factors $\name\QQ_\alpha$ such that for
all $\alpha<\lambda$, $\forces_{\PP_\alpha}$ ``$\name\QQ_\alpha$ is
proper.'' Then $\PP_\lambda$ is proper.
\endproclaim

See [Sh b, Sh f, or Jech] for the proof.

In \S\S1,2 we need additional iteration theorems discussed in [Sh b]
in the context of $\omega$-proper forcing. Improvements in [Sh 177] or
[Sh f] make this unnecessary, but we include a discussion of the
relevant terminology here.  This makes our discussion compatible
with the contents of [Sh b].

\proclaim{A1.6\\$\alpha$-Proper Forcing}
\endproclaim
Let $\alpha$ be a countable ordinal.  Then $\PP$ is $\alpha$-{\sl proper}
iff for every $\PP$-large $\lambda$, every continuous increasing $\alpha+1$-sequence
$(\MM_i)_{i\le\alpha}$ of countable elementary substructures of
$\VV_\lambda$ with $\PP\in \MM_0$, every $p\in P\intersect M_0$ has an
extension $q\in P$ which is $\MM_i$-generic for all $i\le \alpha$.

Axiom A implies $\alpha$-properness for $\alpha$ countable. For
example we check $\omega$-properness. So we consider a condition $p$
in $M_0$, where $(\MM_i)_{i<\omega}$ is a sequence of suitable
countable models satisfying, among other things, $\MM_i\in \MM_{i+1}$.
There is an $\MM_0$-generic condition $p_1$ above $p$, and we can take
$p_1\in \MM_1$, since $\MM_1\elementary \VV_\lambda$. Similarly we can
successively find $p_{n+1}\in P\intersect M_{n+1}$ with $p_{n+1}$
$\MM_n$-generic, and $p_n\le_n p_{n+1}$. A final application of Axim
A yields $q$ above all the $p_n$.

Countable support iteration also preserves $\alpha$-properness for
each $\alpha$ [Sh b].  Furthermore it is proved in [Sh b, V4.3] that
countable support iteration preserves the following conjunction of
two properties: $\omega$-properness and $\up\omega\omega$-bounding.
So [Sh b] contains most of the information needed in \S\S1,2, though
we will need to add more concerning the iteration theorems below.
\par


\proclaim{A2. Iteration theorems}
\endproclaim

\proclaim{A2.1\\Fine* Covering Models}
\endproclaim

We recall the formalism introduced in [Sh b, Chap. VI] for proving
iteration theorems.  We consider collections of subtrees of
$\up{\omega>}\omega$ that cover $\reals$ in the sense that
every function in $\up\omega\omega$ represents a branch of one of the
specified trees, and iterate forcings that do not destroy this
property.  Of course the precise formulation is considerably more restrictive.
 See discussion A2.6.
\medskip

\noindent {\bf Weak covering models.}

A structure $(D;R)$ consisting of a set $D$ and a binary relation $R$ on
$D$ is called a {\sl weak covering model} if:
\item{1.} For $x,t\in D$, $R(x,t)$ implies that $t$ is a (nonempty)
subtree of ${}^{\omega>}\omega$, with no terminal nodes (leaves); we
denote the set of branches of $t$ by $\br(t)$.
\item{2.} For every $\eta\in \up\omega\omega$, and every $x\in \dom
R$, there is some $t\in D$ with $R(x,t)$ and $\eta\in \br(T)$. In this
case, we say: $(D,R)$ covers $\up\omega\omega$.

$(D;R)$ should be thought of as a suitable small fragment of a universe
of sets, and $R(x,t)$ is to be thought of intuitively as saying, in some manner,
that the tree $t$ has ``size'' at most $x$.  In the next definition we
introduce an ordering on the ``sizes'' and exploit more of our
intutition, though certain intuitively natural axioms are omitted, as
they are never needed in proofs.
\medskip

\noindent {\bf Fine* Covering Models.}

A structure $\DD=(D;R,<)$ is called a {\sl fine* covering model} if $(D;R)$
is a weak covering model, $<$ is a partial order on $\dom R$ with no
minimal element,  and:

\noindent\item{(1)} If $x,y\in \dom R$ with $x<y$, then there is $z\in \dom R$
with $x<z<y$ (and $D\neq \emptyset$ and for every $y\in D$
there is $x<y$ in $D$).
\noindent\item{(2)} $x<y\and R(x,t)$-implies $R(y,t)$.
\item{(3)} In any generic extension $V^*$ in which $(D;R)$ is a weak
covering model we have:

\itemitem{$(*)$} for $x<y$ (from $\dom R$) and $t_n\in D$ with $R(x,t_n)$
for all $n$ there is $t\in D$ with $R(y,t)$ holding and there are indices
$n_0<n_1<\ldots$ such that: for all $\eta\in \reals$: if
$\eta\restrictedto n_i \in\bigcup_{j\le i} t_j$ for all $i$ then
$\eta\in\Br(t)$.

\itemitem{$\otimes$} if $\eta\in \reals$, $\eta_n\in \reals$,
$\eta_n\restrictedto n=\eta\restrictedto n $ for $n<\omega$ and $x\in\dom R$
then for some $t$, $R(x,t)$, $\eta\in\Br(t)$ and for infinitely many
$n$ we have $\eta_n\in \Br(t)$.

\noindent In particular we require $(*)$ and $\otimes$ to hold
in the original universe $V$.
Observe also that in (3$*$) we have in particular $t_0\subseteq t$.

\medbreak
\noindent Note that (3)$^+$ below implies (3).

\noindent\item{(3)$^+$} In any generic extension $V^*$ (of $V$) in which $(D,R)$
is a weak covering model we have:

\itemitem {$(*)^+$} For $x<y$ and $t_n\in D$ with $R(x,t_n)$ for all $n$,
there is $t\in D$ with $R(y,t)$ holding and there are indices  
$0=n_0<n_1<\ldots$ such that for all $\eta\in \reals$ if
$\eta\restrictedto n_i \in \bigcup_{j\le i} t_{n_j}$ for all $i$,
then $\eta\in Br(t)$; we let $w=\{n_0, n_1,\ldots\}$.

[Why $(3)^+\Rightarrow (3)$? assume $(3)^+$,
so let a generic extension $V^*$ of $V$ in which $(D, R)$ is a weak covering
model be given,
so in $V^*$, $(*)^+$ holds.
First, for $\otimes$ of $(3)$ let $\eta, \eta_n, y$ be given, let $x<y$;
as $``(D, R)$ is a weak covering model in $V^*\,$'' for each $n<\omega$
there is $t_n\in D$ such that
$R(x, t_n)\&\eta_n\in \Br(t_n)$.
Apply $(*)^+$ to $x,y, t_n$ and get $t$ which is as required there.
Second, for $(*)$ of (3), let $x<y, t_n(n<\omega)$ be given.
Choose inductively $y', x_n, x<x_n<y'<y$, $x_n<x_{n+1}$
(possible by condition (1)).
Choose by induction on $n, k_n, t^*_n$ such that:
$t^*_0=t^*, R(x_n, t^*_n),
t^*_n\subseteq t^*_{n+1}$ and
$[\nu \in t_{n+1}\& \nu\restrictedto k_n\in t^*_n
 \Rightarrow \nu\in t^*_{n+1}]$.
 For $n=0$-trivial, for $n+1$ use $(*)^+$ with
 $\lk x_n, x_{n+1}, t^*_n, t_{n+1}, t_{n+1}, \ldots\rk$ here standing for
 $\lk x,y,t_0, t_1, t_2, \ldots\rk$ there, and we get $t^*_{n+1}$, $w_n$ 
 (for $t, w$ there),
 let $k_n=\Min(w_n\setminus\{0\})$,
 easily $t^*_n$ as required. Now apply $(*)^+$ to
 $\lk y', y, t^*_0, t^*_1, \dots \rk$ and get $t,
 \lk n_i:i<\omega\rk$;
 thining the $n_i$'s we finish].

A forcing notion $\PP$ is said to be $\DD$-{\sl preserving} if $\PP$
forces: ``$\DD$ is a fine* covering model''; equivalently, $\PP$
forces: ``$(D;R)$ covers $\up\omega\omega$.''  So this means that
$\PP$ does not add certain kinds of reals.

In this terminology, we can state the following general iteration
theorem ([Sh 177],[Sh-f]VI\S1, \S2):

\proclaim{A2.2\\Iteration Theorem}
Let $\DD$ be a fine* covering model.
Let $\<\PP_\alpha,\name \QQ_\beta:\alpha\le\delta,\beta<\delta\>$
 be a countable
support iteration of proper forcing notions with each factor
$\DD$-preserving. Then $\PP_\delta$ is $\DD$-preserving.
\endproclaim

\proof:
We reproduce the proof given in [Sh b, pp. 199-202], with the modifications
suggested in [Sh 177]. We note that in the present exposition we have
suppressed some of the terminology in [Sh b] and made other
minor alterations. In particular our statement of the main theorem is
slightly weaker than the one given in [Sh f].  We have also suppressed
the discussion of variants of condition $(3*)$ in the definition of
fine* covering model, which occurs on  pages 197-198 of [Sh b]; as a result
we leave a little more to the reader.

By [Sh b, V4.4], if $\delta$ is of uncountable cofinality then
there is no problem, as all new reals are added at some earlier point.
So we may suppose that $\cf\delta=\aleph_0$ hence by associativity of CS
iterations of proper forcing ([Sh-b], III) without loss of generality
$\delta =\omega$.

We claim that $\forces_{\PP_\omega}$ ``$(D;R)$ covers $\reals$.''
(Note that this suffices for the proof of the iteration theorem.) 

Fix $x\in \dom R$, $p\in \PP_\omega$, $\name f$ a $\PP_\omega$-name
with $p\forces$ ``$\name f\in \reals$.'' We need to find an extension
$p'$ of $p$ and a tree $t\in D$ with $R(x,t)$ such that $p'\forces$
``$\name f\in\br(t)$.'' As in the proof that countable support
iteration preserves properness, we may assume without loss of
generality (after increasing $p$) that $\name f(n)$ is a $\PP_n$-name
for all $n$.

By induction on $n$ we define conditions $p^n\in \PP_n$ and
$\PP_m$-names $\name t_{m,n}$ for $m\le n$ with the following
properties:
\medskip
\+$(1)$	&$\forces_{\PP_i}$ ``$p(i)\le p^n(i)\le p^{n+1}(i)$'' for
$i<n$;\cr

\+$(2)$	&If $G_m\includedin \PP_m$ is generic with $m\le n$, then in
$V[G_m]$ we have \cr
\+&$(p^n(m),\ldots,p^n(n-1))\forces_{\PP_n/\PP_m}$
``$\name f(n)=\name t_{m,n}$.''\cr
This is easily done; for each $n$, we increase $p^n$ $n$ times, once
for each possible $m$. By $(1)$ we have $p\restrictedto n\le p^n\le
p^{n+1}$.

We let $\name f_m$ be the $\PP_m$-name for an element of $\reals$ satisfying:
$\name f_m(n)=\name t_{m,n}$ for $n\ge m$, $\name f_m(n)=\name f(n)$ for
$n<m$. Then we have:

\+$(3)$	&$(0,\ldots,0,p^n(m))\forces_{\PP_{m+1}}$ ``$\name
f_m\restrictedto n=\name f_{m+1}\restrictedto n$''\cr
\+$(4)$&$\forces_{\PP_n}$ ``$\name f\restrictedto n=\name
f_n\restrictedto n$.''\cr

Choose $x_1<x'<x$ and then inductively $x_1<x_2<\ldots$ with all
$x_n<x'$, and choose a countable $\Nnorm\elementary V_\lambda$
(with $\lambda$
$\PP$-large) such that  all the data $(x_n)_{n<\omega}$,
$(\PP_n,\QQ_n)_{n<\omega}$, $\name f$, $(p^n)_{n<\omega}$,
$(\name t_{m,n})_{m\le n<\omega}$  lie in $\Nnorm$. We will define conditions
$q^n\in \PP_n$ and trees $t_n\in D$ (not names!) by induction on $n$ with
$q^{n+1}\restrictedto n=q^n$ (hence we may write:
$q^n=(q_0,q_1,\ldots,q_{n-1})$) and $t_n\includedin t_{n+1}$,
satisfying the following
conditions:

\+$(A)$&$p\restrictedto n\le q^n$;\cr
\+$(B)$&$q^n$ is $(\Nnorm,\PP_n)$-generic;\cr
\+$(C)$	&$q^n\forces$ ``$\name f_n\in \br(t_n)$'';\cr
\+$(D)$	&$R(x_{3n},t_n)$;\cr
\+$(E)$	&For $m<n<\omega$ we have
$q^m\forces_{\PP_m}$ ``$q_m$ and $p^n(m)$ are
compatible in $\name \QQ_m$''.\cr

Suppose we succeed in this endeavour.  Then we can let $q = \Union_n
q^n$.
 By condition (2) in A2.1 for every $n<\omega$ $R(x\up\prime,t_n)$
(as $x_{3n}<x\up\prime)$.
Let
$(n_i:\,i<\omega)$
be a strictingly increasing sequence of natural numbers
and $t$ be as guaranteed by ($*$) of condition (3) of A2.1
(for $\<t_n:n<\omega\>, x\up\prime, x$)
 so $R(x,t)$ and: if $\eta\restrictedto n_i\subseteq \bigcup_{j\le i}t_j$
 for each $i<\omega$ then $\eta\in t$.
 Let $g(i)=:n_i$.

By (E) above there are conditions $q_m'$
 with $q^m \forces_{\PP_m}$ ``$q_m'\in \name \QQ_m$,
$q_m'\ge q_m,p^{g(m)}(m)$.'' Let $q'=(q'_0,q'_1,\ldots)$. Then $q'\ge
q\ge p$
 and for $m\le n\le g(m)$ we will have (if we succeed in defining $q_n, t_n$)
 $q'\restrictedto n\forces_{\PP_n}$
 ``$\name f\restrictedto n=\name f_m\restrictedto n$'', hence:
$$q'\restrictedto n\forces_{\PP_n}\,\hbox{``$\name f\restrictedto n\in
\br\,t_m$''.}$$
Now we have finished proving the existence of 
$p',t$ (see before (1)) as required:   
$q'\forces$ ``$\name f \in \br\, (t)$'', as  $t$ includes  the tree:
$\{\eta\in \up{\omega>}\omega$:  For all $i$, $\eta\restrictedto n_i\in
\Union_{j\le i}t_j\}$; and $R(x,t)$ holds. 
Hence we have finished proving 
$ \forces_{\PP_\omega} ``(D; R)$ covers ${}^\omega\omega$''. 
 So it suffices to carry
out the induction.
\def\fhat{\name{\hat f}_{n+1}}

There is no problem for $n=0$ or $1$. Assume that $q^n$ and $t_n$ are
defined. Let $G_n\includedin P_n$ be generic with $q_n\in G_n$.  Then
$\name f_{n+1}$ becomes a $\name \QQ_n[G_n]$-name $\fhat=\name
f_{n+1}/G_n$ for a member of $\reals$. As $\PP_{n+1}$ preserves
$(D,R)$, for every $r\in \name \QQ_n[G_n]$ and every $y\in \dom R$
there is a condition $r'\ge r$ in $\name \QQ_n[G_n]$ such that
$$r'\forces \hbox{``$\fhat \in \br(t')$'' for some $t'\in D$ with $R(y,t')$.}
\leqno(*)$$
For each $m<\omega$, applying this to $r=: p^m(n)$, $y=x_{3n}$ we get
$r' = r^n_m$, $t'=t^n_{m+1}$;
we could have guaranteed $t^n_{m+1}\subseteq t^n_{m+2}$.
Now choose by induction on $l<\omega$, $r^n_{m,l} \in \name \QQ_n[G_n]$
such that: $r^n_{m,0} = r^n_m$, $r^n_{m,l}\le r^n_{m,l+1}$,
$r^n_{m,l+1}$ forces a value to $\name{\hat f}_{n+1} \restrictedto l$.
So for some $\eta^n_m\in \up{\omega}{\omega}[G_n]$, $r^n_{m,l} \forces
$``$\name{\hat f}_{n+1} \restrictedto l = \eta^n_m\restrictedto l$''.
Note $\eta^n_m\restrictedto m = f_n \restrictedto m$.
%
%
%
 Without loss of generality, $\langle r^n_m, t^n_m, r^n_{m,\ell}, \eta^n_m: 
n,m,\ell<\omega\rangle$ belongs to $N$.  
Applying (3$\otimes$) from A2.1
(to $\eta=\name f_n[G_n]$,
$\eta_m=\eta^n_m$) we can find
$T^{I}_n\in D\cap N[G_n]=D\cap N $
such that $R(x_{3n},T^{I}_n)$,
$\name f_n\in \Br(T^{I}_n)$ and
 $\eta^n_m\in \Br (T^{I}_n)$
 for infinitely many $m<\omega$.
 Applying (3$*$) from A2.1
 (to $T^{I}_n, t^n_1, t^n_2, \ldots$
  and $x_{3n},x_{3n+1}$) we obtain a tree $T^{II}_n$.
 Returning to $V$, we have a $\PP_n$-name $\name T$ for such a tree. For
$s\in \PP_n$, if $s\forces$ ``$\name T=T$'' for some tree $T$ in $V$,
let $T(s)$ be this tree. Let $U$ be the open dense subset of $s\in\PP_n$
for which $T(s)$ is defined.  Some such function $T(\cdot)$ belongs to
$\Nnorm$, and $U\in \Nnorm$. If $q^n$ is in the generic set $G_n$, then some
$s\in U\intersect \Nnorm$ is in $G_n$, by condition $(2)$.  Let
$U\intersect \Nnorm=\{s_i:i<\omega\}$.

\noindent Applying (3$*$) there is a tree $t_{n+1}$ satisfying:

\item {(a)} $R(x_{3n+3}, t_{n+1})$.

\item {(b)} $t_n\subseteq t_{n+1}$.

\item{(c)} for every $T\in (\Rang R)\cap N$ such that $R(x_{3n+2}, T)$
 for some $k_T<\omega$ we have:
 $${\nu}\in T\and {\nu}\restrictedto k_T\in t_n \Rightarrow {\nu}\in t_{n+1}$$

 We shall prove now

\item{(d)} suppose $G_n\subseteq \PP_n$ is generic over $V$
with $q^n\in G_n$, and $k\up *<\omega$.
\item{} {\it Then} there is $q'$,
$p^{k\up *}(n)\le q\up \prime\in \name \QQ_n[G_n]\cap N[G_n]$, such that 
\item{} $q'\forces$``$\name f_{n+1}\in \Br (t_{n+1})$''
\item{} (though $t_{n+1}$ is generally not in $N$).

\proof\ of (d):
As $q^n \in G_n$ necessarily for some $s\in P_n\cap N$ we have
$s\in G_n$ so (c) applies to $T_s$ and $T_s=\name T^{II}_n[G_n]$
(as $T^{II}_n=\name T^{II}_n[G_n]$ is well defined and also $T^{I}_n$
is well-defined and belongs to
$ N\cap D$
 not only $N[G_n]\cap D$, as $D \subseteq V$).
 By the choice of $T^{I}_n$ the following set is infinite
 $$ w=\{ i< \omega: \eta^n_i\in \Br (T^{I}_n)\}$$
 By the choice of $t^n_{i+1}$,  for every $i\in w$ there exists $
k_i< \omega $ such that  $\eta\in t^n_{i+1} \and \eta\restrictedto k_i =
\eta^n_i\restrictedto k_i \implies \eta \in T^{II}_n$.   To show
(d), choose $i\in w \setminus  k^* $ (exists as $w$ is infinite, $k^*$
will be shown to be as required in (d)).

Now $r^n_{i,k}\in N \cap \name \QQ_n[G_n]$ is well-defined, and any $q'$,
$p^i(n)\le q'\in \name \QQ'_n[G_n]$ which is $(N, \name
\QQ_n[G_n])$-generic is as required (note that $p^{k^*}(n) \le p^i(n)$).

We can assume without loss of generality that $\QQ_n$ is closed under
countable disjunction, so we can find $\name r_n$ compatible with
$p^n(m)$  for all $m$ such that:

\centerline{$(q_0,\ldots,q_{n-1},q'_n)\forces_{\PP_{n+1}}$ ``$\name
f_{n+1}\in\br\, (t_{n+1})$''.}

\noindent Now find $q_n\ge q_n'$ such that
$(q_0,\ldots,q_{n-1},q_n)$ is $(\Nnorm,\PP_{n+1})$-generic. This
completes the induction step.

[If this infinite disjunction bothers you, define by induction on $n$
sequences $\<q^n_\eta:\eta\in \up {n+1} \omega \>$ where $q^n_\eta\in
\name \QQ_n$ is such that for every $\eta \in \up m \omega $ the condition
$\<q^i_{\eta\restrictedto (i+1)}:i<n \>$ is generic for $N$ and
$q^n_\eta$ is above $p^{\eta(n)}(n)$.] \hfill\qed

\proclaim{A2.3\\The $\reals$-bounding property}
\endproclaim

We leave the successor case to the reader (see A2.6(2)).

A forcing notion $\PP$ is $\reals$-bounding if it forces every
function in $\reals$ in the generic extension to be bounded by one in
the ground model.  In \S1 we quoted the result that a countable
support iteration of proper $\reals$-bounding forcing notions is again
$\reals$-bounding, which is almost Theorem V.4.3 of [Sh b]. In Chapter VI,
\S2 of [Sh b] this result is shown to fit into the framework just
given.  Here $D$ is just a single collection $\TT$ of trees; to fit
$D$ into the general framework given previously, we would let $A$ be any
suitable partial order, $D=A\disjointunion \TT$, and $R=A\times \TT$.
The set $\TT$ will consist of all subtrees of $\up{\omega>}\omega$
with finite ramification (as we have no measure on how small $t\in
\TT$ is, so $<, R$ are degenerate). 

In a generic extension of the universe, the set $\TT$ (as defined in
the ground model) will cover $\reals$ if and only if every function in
$\reals$ is dominated by one in the ground model. In fact the only
relevant trees are those of the form $T_{f}=\{\eta\in
\up{\omega>}\omega  
  :\eta(i)\le f(i) \hbox{ for $i<\len \eta$}\}$ with $f$
in the ground model.
Thus the $\reals$-bounding property coincides
with the property of being $\DD$-preserving, where $\DD$ is
essentially $\TT$, more precisely $\DD=(A\times \TT;R,<)$ for a
suitable $R,<$ (which play no role in this degenerate case).  Thus to
see that the general iteration theorem applies, it suffices to check
that such a $\DD$ will be a fine* covering model. We  have to check
the final clause  
$(3)$ of the definition of fine* covering model.  In fact we will
prove a strong version of $(3)^+$. 

\item{} For any sequence of trees $T_n$ in $\TT$, there is a tree $T$
such that for all $\eta\in \reals$, if  $\eta\restrictedto i\in
\Union_{j\le i} T_j$  for all $i$, then $\eta\in \br(T)$.

We will verify that this property  holds in any generic extension $V^*$
of $V$ in which $\DD$ covers $\reals$.
Let $T^*=\{\eta\in\up{\omega>}\omega:$ for all
$i\le\len(\eta)$, $\eta\restrictedto i\in\Union_{j\le i} T_j\}$.  If
$T^*$ is in $V$ this will do, but since the sequence $(T_n)$ came from
a generic extension, this need not be the case.  On the other hand the
sequence $T^*\restrictedto n$ of finite trees is itself coded by a
real $f\in \reals$, and as $\DD$ covers $\reals$, there is
a tree $T^\bullet$ in $D$ which contains this code $f$; via a
decoding, $T^\bullet$ can be thought of as a
tree $T^o$ whose nodes $t$ are subtrees of $\up {n\ge}\omega$ with no
 maximal nodes below
level $n$, so that for any $s,t\in T^\bullet$ with $s\le t$, $s$ is
the restriction of $t$ to the level of $s$, and such that the sequence
$T^*\restrictedto n$ actually is a branch of $T^o$. Let $T$ be the subtree
of $\up{\omega>}\omega$ consisting of the union of all the nodes of
$T^o$. Then $T$ still has finite ramification, lies in the ground
model, and contains $T^*$.

\proclaim{A2.4\\Cosmetic Changes}
\endproclaim

(a) We may want to deal just with $\Br(T\up{*}\;)$, where 
$T\up{*}$ a subtree $^{\omega>} \omega$ (hence downward closed).
So $D$ is a set of subtrees of $T\up{*}$ , so 
we can replace $D$ by $\{\{\eta\in\up{\omega>}\omega:\eta\in T$ or
$(\exists\ell)[\eta\restrictedto\ell \in T
\and \eta\restrictedto (\ell+1)\notin T\up{*}\;\}: T\in D\}$.

(b) We may replace subtrees $T\up{*}$ of $\up{\omega>}\omega$
 by isomorphic trees.

(c) We may want to deal with some $(D_i;R_i,<_i)$ simultaneously;
 by renaming without loss of generality the $D_i$ are pairwise disjoint,
and even: $\bigwedge_{\ell=1,2} t_l\in D_{i_l}\&
i_1\not= i_2 \implies \Br T_1 \cap \Br t_2
= \emptyset$. Then we use  ($\Union D_i; \Union R_i,\Union{ <}_i)$
to get the result.

(d) We may want to have $(D;R)$ (i.e. no $<$); just use
 $(D\cup\Q\times D; R', <)$ where $R'(x,t)$ iff
 $x=(q,y)$, $q\in\Q,y\in D, R(y,t), (q_1,y_1)<(q_2,y_2)$ iff
 $q_1<q_2\and y_1=y_2$.

\proclaim{A2.5\\The $(f,g)$-bounding property}
\endproclaim

We leave the successor case to the reader (see A2.6(2)).

Let $\unF$ be a family of functions in $\reals$, and $g\in\reals$ with
$1<g(n)$ for all $n$.  We say that a forcing notion $\PP$ has the
$(\unF,g)$-{\sl bounding property} if:
$$\vcenter{\hsize=0.8\hsize \advance\rightskip0pt plus 1 cm minus 1 cm
\noindent
For any sequence $(A_k:k<\omega)$ in the ground model,
with $|A_k|\in\unF$ (as a function of $k$), and any
$\name\eta\in\prod_kA_k$ in the generic extension and $\epsilon>0$,
there is a ``cover'' $\BB=(B_k:k<\omega)$ in the ground model with
$B_k\includedin A_k$,
$[|B_k|>1\Rightarrow|B_k|<g(k)\up{\epsilon}\;$] and 
$\name\eta(k)\in B_k$ for each $k$.\par}\leqno(*)$$
This notion is only of interest if $g(n)\to\infty$ with $n$.

We will show that this notion is also covered by a case of
the general iteration theorem of \S A2.2.

Let $\TT_{f,g}$ $[\TT_{f,g}^{\epsilon}]$
be the set of those subtrees $T$ of $\Union_n\prod_{m<n}
f(n)$ of the form
$\Union_n\prod_{m<n} B_m$,
such that $|B_n| <\max \{g(k), 2\}$ [such that $|B_n| \le \max \{2,
g(k)^\varepsilon\}$],
where as usual
$f(n)$ is thought of as the set $\{0,\ldots,f(n)-1\}$.
Let $\TT_{\unF,g}$
be $\Union_{f\in\unF, \varepsilon\in \Q^+} \TT^\varepsilon_{f,g}$.
Our fine* covering model is essentially $\TT_{\unF,g}$, more accurately, 
it is the family of $\{(\TT_{f,g}\cup \Q^+; R, {<}): f\in \unF\}$,
 where $\Q^+$ is the set of positive rationals, $<$ is the order on $\Q^+$, and 
$R(\varepsilon,t) =: \varepsilon \in \Q^+ \and t\in
\TT^\varepsilon_{f,g}$.  See A2.4(c).

Call a family $\unF$ {\sl $g$-closed} if it satisfies the following
two closure conditions:

\item{1.} For $f\in \unF$, the function $F(n)=\prod_{m<n}(f(m)+1)$ lies in
$\unF$;
\item{2.} For $f\in \unF$, $f^g$ is in $\unF$.

If $\unF$ is $g$-closed, $f\in \unF$, and $(A_n)_{n<\infty}$ are sets
with $|A_n|=f(n)$, then the function $f'(n)=$ the number of trees
of the form $\prod_{m<n}B_m$ with $B_m\includedin A_m$ and
$|B_m|<g(m)$ is dominated by a function in $\unF$.

Using the formalism of \S A2.2, we wish to prove:

\proclaim{Theorem}
If $\unF$ is $g$-closed then
a countable support iteration of $(\unF,g)$-bounding proper forcing notions
is again an $(\unF,g)$-bounding proper forcing.
\endproclaim

Since the $\DD$-preserving forcing notions are the
same as the $(\unF,g)$-bounding ones, we need only check that $\DD$ is a
fine* covering model.  Again the nontrivial condition is $(3)^+$,
i.e.,

\item{} Let $f\in \unF$. For any sequence of trees $T_n$ in
$\TT_{f,g}$, $R(\varepsilon ' , T_n)$, $\varepsilon ' < \varepsilon$
(in $\Q^+$), there is a tree $T$ in $D$ satisfying $R(\varepsilon, T)$
and an increasing  sequence $n_i$ such that for all $\eta\in \reals$,
if  $\eta\restrictedto n_i\in
\Union_{j\le i} T_{n_j}$  for all $i$, then $\eta\in \br(T)$.

This must be verified in any generic extension $V^*$
of $V$ in which $\DD$ covers $\reals$.  Working in $V$, choose
$(n_i)_{i<\omega}$ increasing so that  $n_0=0$ and for $n_i\le n$ we have
$\min_{n\ge n_i} g(n)^{(\varepsilon  - \varepsilon ' ) / 2} > i+1$. 
For $n_i\le n<n_{i+1}$ set:
$$B_n=\{\eta(n):\eta\in \cup \{T_j:n_j \le n\}.$$
(For $n <n_0$ let $B_n = \{\eta(n): \eta \in T_0\}$.)
If the sequence $B_n$ was in
the ground model, we could take $T=\Union_n\prod_{m<n} B_m$. Instead
we have to think of the sequence $B_n$ as a possible branch through
the tree of finite sequences of subsets of $f(n)$ of size at most
(say) $\max\{1, g(n) - 1\}$. As $\unF$ is $g$-closed, $\TT_{\unF,g}$
contains a tree 
$T^\bullet$ which encodes a tree $T^o$ of such subsets, for which the
desired sequence $B_n$ is a branch in $V^*$, so that the number of members
of $T\up{0}$ of level $m$ is $\le g(m)\up{(\epsilon-\epsilon')/2}$
(or is $\le 1$).
  Let
$B^o_n=\Union_{b\in T^o} b(n)$.  Then $B_n\includedin B^o_n$,
$\lim_{n\to\infty} |B^o_n|/g^\epsilon(n)=0$
 and $\Union_n \prod_{m<n}(B^o_m)$ is in $V$.\qed

\proclaim{A2.6\\Discussion:}
\endproclaim

This was treated in [Sh-f,VI] [Sh-f, XVIII \S 3] too
(the presentation in [Sh-b, VI] was inaccurate).
The version chosen here goes for less generality
(gaining, hopefully, in simplicity and clarity) and is usually sufficient.
We consider below some of the differences.

\proclaim{A2.6(1)\\A technical difference}
\endproclaim

In the context as phrased here the preservation in  the successor case
of the iteration was trivial --- by definition essentially. 
We can make the fine* covering model (in A2.1) more similar to
[Sh-f, VI \S1] by changing (3$*$) to 
$$(*)'\qquad \vcenter{\hsize=0.8 \hsize\advance\rightskip0pt plus 1 cm
minus 1 cm 
\noindent For
$y_0<y_1<\ldots y<x$ in $\dom R$ and $t_n\in D$ such that
$R(y_n,t_n)$ for all $n$, there is $t\in D$ with $R(x,t)$ holding
and indices $n_0<n_1<\ldots$ such that
$[\eta\in\up{\omega>}\omega\and\bigwedge_i\eta\restrictedto\eta_i\in
\bigcup_{j\le i}t_j\Rightarrow\eta\in t]$.\par} $$ 
 We can use this version here.

\proclaim{A2.6(2) Two-stage iteration}
\endproclaim

We can make the fine* covering model  (in A2.1) more similar to [Sh-f,
VI \S1] by changing (3$*$).  In the context as presented here the
preservation by two step iteration is trivial --- by definition
essentially.  In [Sh-f VI, \S2] we phrase our  framework such that we
can have: 
if $Q_0\in V$ is $x$-preserving,
$\name Q_1$ is $X$-preserving (over $V\up{Q_0}$, $\name Q_1$ a $Q_0$-name)
then $Q_0*\name Q_1$ is $x$-preserving.
The point is that $X$-preserving means $(D,R,<)\up{V}$-preserving,
i.e. $(D,R,<)$ is a definition
(with a parameter in $V_0$).
The point is that if $V_1=V_0\up{Q_0}$, $V_2=V_1^{\supername Q_1}$
{\it then} for $\eta\in(\reals)\up{V_2}$ and $x\in \dom R$,
we choose $y<x$ and $t\in D\up{V_1}$ such that
$\eta\in \Br(t)$,
$R\up{V_1}(y, t)$, then we look in $V_0$ at the tree of possible
initial segments of $t$ getting $T\in D\up{V_0}$
such that $t\in \Br(T)$ , $R\up{V_0}(y, T)$.
If $y$ was chosen rightly, $\Union\Br(T)$ is as required.
Here it may be advantegous to use a preservation of several $(D,R,<)$'s
at once (see A2.4(c)).

\proclaim{A2.6(3) Several models --- the real case}
\endproclaim

We may consider a (weak) (fine*) covering family of models
$\<(D_\ell, R_\ell, <_\ell):\ell<\ell\up{*}\>$ (actually a sequence)
i.e. {\it not} that each one is a cover, but simultaneously.

\noindent (A) We say $(\overline{D}, \overline{R})=
\<(D_\ell, R_\ell):\ell<\ell\up{*}\>$
 is a weak c.f.m.\ if each $D_\ell$ is a set,
$R_\ell$ a binary relation, $\ell^*<\omega$ and

\item{1.} $R_\ell(x,t)$ implies that $t$ is a subtree of
$\up{\omega>}\omega$
(nonempty, no maximal models).

\item{2.} Every $\eta\in\reals$ is of kind $\ell$
for at least one $\ell<\ell\up{*}$ which means:
for every $x\in\dom R_\ell$ for some $t$,
 we have $R_\ell(x,t)\and \eta\in \Br(t)$.

\noindent (B) We say $(\overline D, \overline R, \overline <)$
is a fine*
c.f.m.\ if:

\item{0.} $(\overline D,\overline R)$ is a weak family.

\item{1.} If $x\in \dom R_{\ell}\Rightarrow (\exists z) z<_{\ell} x$
and $\forall y<_{\ell} x\;\exists z (y<_{\ell}z<_{\ell} x)$
 (and $D\neq \emptyset$).

\item{2.} $x<_\ell y\and R_\ell(x,t)\Rightarrow R_\ell(y,t)$.

\item {3.} For any generic extensions $V^*$ in which
$(\overline D, \overline R)$ is a weak c.f.m.

\itemitem{$(*)$} for every $\ell<\ell\up{*}$ and $y<_\ell x$ (from $\dom
R_\ell$)
and $t_n\in D_\ell$ with $R_\ell(y, t_n)$ for all $n$ there is $t\in D$ with
$R_\ell(x, t)$ and there are indices $n_0<n_1<\ldots$ such that for every
$\eta\in \reals$:
if $\eta\restrictedto n_i\in \bigcup_{j\le i} t_j$ for all
$i$ then $\eta\in\Br(t)$.

\itemitem{$\otimes$} if $\ell<\ell^*$, $\eta\in {}^\omega\omega$,
 $\eta_n\in {}^\omega\omega$,
$\eta_n\restrictedto n=\eta\restrictedto n$, $x\in \dom R_\ell$
 and $\eta, \eta_n$ are of kind $\ell$,
{\it then} for some $t\up{*}$, $R_\ell(x, t\up{*})$,
$\eta\in \Br(t\up{*})$ and for infinitely many $n<\omega$,
$\eta_n\in\Br(t\up{*})$.

\proclaim {Theorem} If $(\overline D; \overline R, \overline <)$
 is a fine$^*$ c.f.m., 
$\lk \PP_\alpha, \name\QQ_\beta:\alpha\le \delta, \beta<\alpha\rk$ 
is a countable support iteration of proper forcing notions with each factor 
$(\overline D; \overline R, \overline <)$-preserving. {\it Then} $\PP_\delta$
 is $(\overline D; \overline R, \overline<)$-preserving.  
\endproclaim
\proof:
Similar to the previous one, with the following change.
After saying that without loss of generality  $\delta=\omega$
and, above $p$, for every $n$,
$\name f(n)$ as a $P_n$-name, and choosing $x_n, x'$, we do the following.
For clarity think that our universe $V$ is countable in the true
universe or at least $\beth_3(|P_\omega|)\up{V}$ is.
We let $K=\{(n,p,G):n<\omega, p\in P_\omega$,
 $G\subseteq P_n$ is generic over $V$ and $p\restrictedto n \in G_n\}$.
 On $K$ there is a natural order $(n,p,G)
 \le(n', p', G')$ if $n\le n'$,
 $P_\omega\models p\le p'$ and $G\subseteq G'$.
 Also for $(n,p,G)\in G$ and $n'\in (n,\omega)$
 there are $p', G'$ such that
 $(n,p,G)\le(n',p',G')$.
 For $(n,p,G)\in K$ let
 $L_{(n,p,G)}=\{g:g\in (\reals)\up{V[G]}$
 and there is an increasing sequence
 $\<p_\ell:\ell<\omega\>$ of conditions in
 $P_\omega/ G$, $p\le p_0$, such that $p_\ell \forces
 \name f \restrictedto \ell= g\restrictedto \ell\}$.
 So:
 \nl
 $g\in L_{(n,p,G)}\Rightarrow
 \name f\restrictedto n=g\restrictedto n$
 \nl
 $(n,p,G)\le (n', p', G')$
 $\Rightarrow L_{(n', p', G')}\subseteq
 L_{(n,p,G)}$.

 \claim

  There are $\ell_*$ and $(n,p,G)\in K$ such that
 if $(n,p,G)\le(n',p',G')\in K$ then there is
 $g\in L_{(n', p', G')}$
 which is of the $\ell_*$'th kind.
\endproclaim

 \proof:
 Otherwise choose by induction $(n^{\ell},p^{\ell}, G^{\ell})$
 for $\ell\le \ell\up{*}$, in $K$, increasing such that:
  $L_{(n\up{\ell+1},p\up{\ell+1}, G\up{\ell+1})}$ has no member of the
  $\ell$'th kind.
  So $L_{(n\up{\ell\up{*}}, p\up{\ell\up{*}}, G\up{\ell\up {*}})}=\emptyset$
  contradiction.

  So without loss of generality for every $(n,p,G)\in K$,
  $L_{(n,p,G)}$ has a member of the $\ell_*$'th kind.
  Now we choose by induction on $n$,
  $A_n,\, \lk p_\eta:\eta\in {}^{n+1}\omega,
  \eta\uhr n\in A_n \rk$,
  $\lk \name f_\eta:\eta\in A_n \rk$,
  $\lk q_\eta:\eta\in A_n\rk$, and $t_n$ such that

  \item{(A)$'$} $A_n \subseteq {}^n\omega, A_0=\{\<\>\},
   \eta\in A_n\Rightarrow
    (\exists^{\aleph_0}\ell)(\eta\conc\lk\ell\rk\in A_{n+1})$
    $p_\eta\in \PP_\omega\cap N, p_{<>}=p, p_\eta\le p_{\eta\conc\lk \ell\rk},
    p_\eta\restrictedto n\le q_n$.

\item{(B)$'$} $q_\eta$ is $(N, \PP_{\lg \eta}$)-generic,
$q_\eta\in \PP_{\lg \eta}$ and 
[$ \ell<\lg \eta \Rightarrow q_\eta\uhr \ell=q_{\eta\uhr \ell}$].

\item{(C)$'$} $q_\eta\forces ``\name f_\eta\in \Br (t_n)$ is of the $\ell_*$'th
kind'' when $\eta\in A_n$ and $\name f_\eta$ is a $\PP_n$-name.

\item{(D)$'$} $R_{3n}(x_{3n}, t_n), t_n\subseteq t_{n+1}$.

\item {(E)$'$} $p_{\eta\conc\lk \ell\rk}\forces_{\PP_w}$``$\name f_\eta\restrictedto \ell
=\name f_{\eta\conc\lk \ell\rk}\restrictedto \ell=\name f \uhr \ell$''.

This suffices, as $x_n<x'$ so
$\bigwedge_n R(x', t_n)$
 hence for some $\langle n_i:i<\omega\rangle$ strictly increasing
and $t$ as guaranteed by $(*)$ of (3) we find
$\nu\in {}^\omega\omega$ increasing fast enough and let
 $q=\bigcup_{n<\omega}q_{\nu\restrictedto n}$.
 In the induction there is no problem for $n=0,1$.
 For $n+1$; first for each $\eta\in {}^{n+1}\omega$ we choose $\ell$,
 work in $V^{\PP_{n+1}}$ and find
 $\lk p_{\eta\conc\lk\ell\rk}:\ell<\omega\rk, \name f_\eta$,
 and without loss of generality they are in $N$.
 For $\eta\in {}^n\omega$ there is a $\PP_{n+1}$-name
 $\name t_\eta\in N$ of a member of $D$,
 $R_\ell(x_{3n}, \name t_\eta)$,
 $\name f_\eta\in \Br(\name t_\eta)$,
 $(\exists^\infty\ell)\name f_{\eta\conc \lk\ell\rk}\in \Br (\name t_\eta)$.
  Now we can replace $p_{\eta\conc\lk 1\rk}$ by $p_{\eta\conc\lk \ell'\rk}$,
 $\ell'=\Min\{m:m\ge \ell, \name f_{\eta\conc\lk \ell\rk}\in \Br (\name t_\eta)\}$.
 We continue as in A2.2.
 Note: it is natural to use this framework e.g.
  for preservation of $P$-points. 
\qed

 \proclaim{A3. Omitting types}
\endproclaim

\proclaim{A3.1\\Uniform partial orders.}
\endproclaim

In the proof of Theorem 3.1 given in \S3 we used the combinatorial
principle developed in [ShLH162].  (Cf.  [Sh107] for applications
published earlier.) This is a combinatorial refinement of forcing with
$\AP$ to get a $\Pp_3$-name $\name \FF$ with the required properties
in a generic extension.  We now review this material.

With the cardinal $\lambda$ fixed, a partially ordered set $(\PP,<)$ is
said to be {\sl standard $\lambda^+$-uniform} if $\PP\includedin
\lambda^+\times \PP_{\lambda}(\lambda^+)$ (we refer here to subsets of
$\lambda^+$ of size strictly less than $\lambda$), satisfying the
following properties (where we take e.g.  $p=(\alpha,u)$ and write
$\dom(p)$ for $u$):

\item {1.} If $p\le q$ then $\dom p\includedin \dom q$.

\item {2.}
For all $p,q,r \in \PP$ with $p,q \le r$ there is $r' \in \PP$ so that
$p,q \le r' \le r$ and $\dom r'=\dom p \union
\dom q$.

\item {3.}
If $(p_i)_{i < \delta} $ is an increasing sequence of length less than
$\lambda$, then it has a least upper bound $q$, with domain \
$\Union_{i<\delta} \dom p_i$; we will write $q=\Union_{i<\delta}
p_i$, or more succinctly: $q=p_{<\delta}$.

\item {4.} For all $p \in \PP$ and $\alpha < \lambda^+$ there exists a
$q \in \PP$ with $q \le p$ and $\dom q = \dom p \intersect \alpha$;
furthermore, there is a unique maximal such $q$, for which we write
$q=p \restrictedto\alpha$.

\item {5.}
For limit ordinals $\delta$, $p \restrictedto \delta =\Union_{\alpha <
\delta} p \restrictedto\alpha$.

\item {6.}
If $(p_i)_{i < \delta}$ is an increasing sequence of length less than
$\lambda$, then  $(\Union_{i < \delta} p_i)\restrictedto\alpha\ =
\
\Union_{i < \delta}
(p_i \restrictedto \alpha).$

\item {7.}
{\sl (Indiscernibility)} If $p=(\alpha, v) \in \PP$ and $h:v
\rightarrow v'\includedin \lambda^+$ is an order-isomorphism onto $V'$ then
$(\alpha, v') \in
\PP$.  We write $h[p]=(\alpha, h[v])$.
Moreover, if $q\le p$ then $h[q]\le
h[p]$.

\item {8.}
{\sl (Amalgamation)} For every $p,q \in \PP$ and $\alpha<\lambda^+$,
if $p\restrictedto\alpha \le q$ and $\dom p \intersect \dom q= \dom p
\intersect \alpha$, then there exists $r \in \PP$ so that $p,q \le r$.

It is shown in [ShHL162] that under a diamond-like
hypothesis, such partial orders admit reasonably generic objects.  The
precise formulation is given in A3.3 below.

\proclaim{A3.2 \\Density systems.}
\endproclaim

Let $\PP$ be a standard $\lambda^+$-uniform partial order.  For
$\alpha<\lambda^+$, $\PP_\alpha$ denotes the restriction of $\PP$ to
$p\in \PP$ with domain contained in $\alpha$.  A subset $ G$ of
$\PP_\alpha$ is an {\sl admissible ideal} (of $\unP_\alpha$)
 if it is closed downward, is
$\lambda$-directed (i.e.  has upper bounds for all small subsets), and
has no proper directed extension within $\PP_\alpha$.  For $G$ an
admissible ideal in $\PP_\alpha$, $\PP/G$ denotes the restriction of
$\PP$ to $\{p\in \PP: p\restrictedto \alpha\in G\}$.

If $G$ is an admissible ideal in $\PP_\alpha$ and
$\alpha<\beta<\lambda^+$, then an $(\alpha,\beta)$-{\sl density system}
for $G$ is a function $D$ from pairs $(u,v)$ in $P_\lambda(\lambda^+)$
with $u\includedin v$ into subsets of $\PP$ with the following
properties: \item{(i)} $D(u,v)$ is an upward-closed dense subset of
$\{p\in \PP/G:\dom(p)\includedin v\union \beta\}$; 
\item{(ii)} For pairs
$(u_1,v_1),(u_2,v_2)$ in the domain of $D$, if
$u_1\intersect\beta=u_2\intersect\beta$ and
$v_1\intersect\beta=v_2\intersect\beta$, and there is an order
isomorphism from $v_1$ to $v_2$ carrying $u_1$ to $u_2$, then for any
$\gamma$ we have $(\gamma,v_1)\in D(u_1,v_1)$ iff $(\gamma,v_2)\in
D(u_2,v_2)$.
\medskip

An admissible ideal $G'$ (of $\unP_\gamma$) is said to
{\sl meet} the $(\alpha,{\beta})$-density system $D$ for $G$ if
$\gamma\ge\alpha$,
$G'\ge G$
and for each $u\in P_\lambda(\gamma)$ there is $v\in P_\lambda(\gamma)$
containing $u$ such that $G'$ meets $D(u,v)$.

\proclaim{A3.3\\The genericity game.}
\endproclaim

Given a standard $\lambda^+$-uniform partial order $\PP$, the {\sl
genericity} game for $\PP$ is a game of length $\lambda^+$ played by
Guelfs and Ghibellines, with Guelfs moving first.  The Ghibellines build
an increasing sequence of admissible ideals meeting density systems set
by the Guelfs.  Consider stage $\alpha$.  If $\alpha$ is a successor,
we write $\alpha^-$ for the predecessor of $\alpha$;  if $\alpha$ is a
limit, we let $\alpha^-=\alpha$. Now at stage $\alpha$ 
 for every $\beta<\alpha$ an admissible ideal $G_\beta$
 in some $\PP_{\beta'}$ is given, and one can
check that there is a unique admissible ideal $G_{\alpha^-}$ in
$\PP_{\alpha^-}$ containing
$\Union_{\beta<\alpha}G_{\beta'}$ (remember A 3.1(5)) 
[Lemma 1.3, ShHL 162].  The Guelfs now
supply at most $\lambda$\
density systems $D_i$  over $G_{\alpha^-}$ for
$(\alpha,\beta_i)$ and also fix an element $g_\alpha$
in $\PP/G_\alpha^-$.
Let $\alpha'$ be minimal such that $g_\alpha\in \PP_{\alpha'}$ and
$\alpha'\ge\sup\,\beta_i$.  The Ghibellines then build an admissible ideal
$G_{\alpha'}$ for $\PP_{\alpha'}$ containing $G_\alpha^-$ as well as
$g_\alpha$, and meeting all specified density systems, or forfeit the
match; they let $G_{\alpha''}=G_{\alpha'}\cap \alpha''$ when
 $\alpha\le \alpha''<\alpha'$.
   The main result is that the Ghibellines can win with a little
combinatorial help in predicting their opponents' plans.

For notational simplicity, we assume that $G_\delta$ is an
$\aleph_2$-generic ideal on $\AP\restrictedto \delta$,
when $\cf \delta=\aleph_2 $ which is true
on a club in any case.
\proclaim{A3.4 \bf $\rm Dl_\lambda$.}
\endproclaim

The combinatorial principle $\rm Dl_\lambda$ states that there are subsets
$Q_\alpha$ of the power set of $\alpha$ for $\alpha<\lambda$ such that
$|Q_\alpha|<\lambda$, and for any $A\includedin \lambda$ the set
$\{\alpha:A\intersect\alpha\in Q_\alpha\}$ is stationary.  This follows
from $\diamond_\lambda$ or inaccessibility, obviously, and Kunen showed
that for successors, $\rm Dl$ and $\diamond$ are equivalent.  In addition
$Dl_\lambda$ implies $\lambda^{<\lambda}=\lambda$.  

\bigskip

\proclaim{A3.5\\A general principle}
\endproclaim

\proclaim{Theorem} 

Assuming $Dl_\lambda$, the Ghibellines can win any standard
$\lambda^+$-uniform $\PP$-game. 
\endproclaim

This is Theorem 1.9 of [ShHL 162].

In our application we identify $\AP$ with a standard
$\aleph_2^+$-uniform partial order via a certain coding. We first
indicate a natural coding which is not quite the right one, then
repair it. 
\medskip

\noindent {\bf First Try}

An approximation $q=(\AA,\name \FF,\name{\fakebf\epsilon},\name {\bf k},)$
will be identified with a pair $(\tau,u)$, where $u=\AA$, and $\tau$ is the
image of $q$ under the canonical order-preserving map $h:\AA\with
\otp(\AA)$.  
One important point is that the first
parameter $\tau$ comes from a fixed set $T$ of size $2^{\aleph_1}=\aleph_2$;
so if we enumerate $T$ as $(\tau_\alpha)_{\alpha<\aleph_2}$ then we
can code the pair $(\tau_\alpha,u)$ by the pair $(\alpha,u)$.  Under
these successive identifications, $\AP$ becomes a standard
$\aleph_2^+$-uniform partial order, as defined in \S A3.1.  Properties
1, 2, 4, 5, and 6 are clear, as is 7, in view of the uniformity in the
iterated forcing $\Pp$, and properties 3, 8 were, 
in essence but not formally,  stated in Claim 3.10.

The difficulty with this approach is that in this formalism, density
systems cannot express nontrivial information: any generic ideal meets
any density system, because for $q\le q'$  with $\dom q=\dom q'$, we
will have $q=q'$; thus $D(u,u)$ will consist of all $q$ with $\dom q =
u$, for any density system $D$.

So to recode $\AP$ in a way that allows nontrivial density systems to
be defined, we proceed as follows.
\medskip

\noindent{\bf Second Try}

Let $\iota:\aleph_2^+\with \aleph_2^+\times \aleph_2$ order preserving
where $\aleph_2^+ \times \aleph_2$ is ordered lexicographically.  Let
$\pi:\aleph_2^+\times \aleph_2\to \aleph_2^+$ be the projection on the
first coordinate. First encode $q$ by $\iota[q]=(\iota[\AA],\ldots)$,
then encode $\iota[q]$ by $(\tau,\pi[\AA])$, where $\tau$ is defined
much as in the first try -- a description of the result of collapsing
$q$ into $\otp\pi[\AA]\times\aleph_2$, after which $\tau$ is encoded
by an ordinal label below $\aleph_2$. The point of this is that now
the domain of $q$ is the set $\pi[\AA]$, and $q$ has many extensions
with the same domain.  After this recoding, $\AP$ again becomes a
$\aleph_2^+$-uniform partial ordering, as before.  We will need some
additional notation in connection with the indiscernibility condition.
It will be convenient to view $\AP$ simultaneously from
an encoded and a decoded point of view. 
One should now think of $q\in \AP$ as a
quintuple $(u,\AA,\name\FF, {\fakebf\epsilon},\name {\hbox{\bf k}})$
with $\AA\includedin u\times \aleph_2$.
If $h:u\with v$ is an order isomorphism, and $q$ is an approximation with
domain $u$, we extend $h$ to a function
$h_*$ defined on $\AA^q$ by letting it act as the
identity on the second coordinate.  
Then $h[q]$ is the transform of $q$ using $h_*$, and has domain $v$.

In order to obtain {\sl least} upper bounds for increasing sequences,
it is also necessary to allow some extra elements into $\AP$, 
 by adding formal least upper bounds to increasing sequences of length
 $<\aleph_2$. 

This provides the formal background for the discussion in \S3.  The
actual construction should be thought of as a match in the genericity
game for $\AP$, with the various assertions as to what may be
accomplished corresponding to proposals by the Guelfs to meet certain
density systems.  To complete the argument it remains to specify these
systems and to check that they are in fact density systems.

\proclaim{A3.6\\The Major Density Systems}
\endproclaim

The main density systems under consideration were introduced
implicitly in 3.11.  Suppose that $\delta <\aleph_2$,  $q\in \AP$ with 
$\delta\in\dom q\includedin
\aleph_2$, $q^*_\delta\le q$, and  $\name z$  is a $(\Pp\restrictedto \dom
q)$-name.  Define a density system \def\Ddqz{D^\delta_{q,\name z}}
$\Ddqz(u,v)$
for $u\includedin v\includedin \aleph_3$ with $|v|\le
\aleph_1$ as follows.  First, if $\otp u\le  \otp \dom q$ then let
$\Ddqz(u,v)$ degenerate to $\AP\restrictedto v$.   Now suppose that
$\otp u > \otp \dom q$ and that $h:\dom q\to u$ is an order
isormorphism from $\dom q$ to an initial segment of $u$.
Let $q^*=h[q]$.
Call an element $r$ of $\AP$ a $(u,v)$-{\sl witness} if:

\item{1.} $u\includedin \dom r \includedin  v$;
\item{2.} $r\ge q^*$;
\item{3.} for some $\pp\in \PP\restrictedto \AA^r$ with
$\pp\ge\pp^\delta$, and some
$(\Pp\restrictedto[\AA^r\intersect \delta])$-name $\name x$, $\name
F_\delta(\name x)$ is a $(\Pp\restrictedto[\AA^r\intersect \delta])$-name;
and:
\item{4.} $\pp'\forces_{\Pp\restrictedto \AA^r}$
``$\{n:[\Gamma^1_{k^1_\delta(n)}\models R(\name x(n),\name
 x_\delta(n))\iff \Gamma^2_{k^2_\delta(n)}
 \models\neg R(\name F_\delta(\name x)(n),\name z(n))]\}\in$
 \nl
 \phantom{$\pp'\forces_{\Pp\restrictedto \AA^r}P_{\restrictedto \AA^r}$}
  $\FF^r$.''

Let $\Ddqz(u,v)$ be the set of $r\in \AP$ with $\dom r
=v$ such that either $r$ is a $(u,v)$-witness, or else there is no
$(u,v)$-witness $r'\ge r$.

This definition has been arranged so that $\Ddqz(u,v)$  is trivially
dense.  In \S3 we wrote the argument as if no default condition had
been used to guarantee density, so that the nonexistence of
$(u,v)$-witnesses is called a ``failure of density''.  Here we adjust
the terminology to fit the style of [ShHL 162].

Now we return to the situation described in 3.12. We had
$\Pp$-names $\name F$, $\name k^1$, $\name k^2$, and a condition
$\pp\in\Pp$, satisfying conditions (3,4) as stated there, and we
considered the set $C=\{\delta<\aleph_3:\cof (\delta) = \aleph_2,\, \delta $ is
$(\name F,\name k^1,\name k^2,\pp)$-closed$\}$, and a stationary set
$S_C$ on which $\name F\restrictedto \delta$, $\pp$,
$\epsilon_\delta$, $\name k^1_\delta$ were guessed by $\diamond$. Then
$\name z=: \name F(\name x_\delta)$ is a $(\Pp\restrictedto
\AA^q)$-name for some $q\in G$. Let $u=\dom q$, $q_0=q\restrictedto
\delta$.  Now we consider the following condition used in 3.12:

\item{(iv)} For all $r\ge q$ in $\AP$ such that
$r\restrictedto \delta\in G_\delta$, and 
$\name x$ a $(\Pp\restrictedto \AA^{r\restrictedto \delta})$-name,
with $\name y=: \name F(\name x)$  a 
$(\Pp\restrictedto \AA^{r\restrictedto \delta})$-name, we have:
\itemitem{$(*)_{\name x,\name y}$\ } 
$\pp\forces$
	``The set $\{n:\Gamma^1_{\name k^1(n)}\satisfies
	R(\name x(n),\name x_\delta(n))\hbox{ iff }\Gamma^2_{\name
	k^2(n)}\satisfies 
	R(\name y(n),\name z(n))\}$
\itemitem{}\\\\is in $\name \FF^r$''. 

We argued in 3.12 that we could confine ourselves to the case in which
(iv) holds.  We now go through this more carefully.   Suppose on the
contrary that we have $r\ge q$ in $\AP$ with $r\restrictedto \delta\in
G_\delta$, and a $(\Pp\restrictedto \AA^{r\restrictedto\delta})$-name
$\name x$, so that $\name y=:\name F(\name x)$ is a $(\Pp\restrictedto
\AA^{r\restrictedto \delta})$-name, and a condition $\pp'\ge \pp$, so that
$$	\hbox{ $\pp'\forces$ ``The set $\{n:\Gamma^1_{\name k^1(n)}\satisfies
	R(\name x(n),\name x_\delta(n))$ iff $\Gamma^2_{\name
	k^2(n)}\satisfies 
	R(\name y(n),\name z(n))\}$ is not in $\name \FF^r$''.}$$ 

Let $\alpha>\sup (\dom r)$, $u=\{\delta\}\union \dom r\union\{\sup\, \dom
r\}$.  Let $q^*\in G$, $q^*\ge r\restrictedto \delta,q$, and let $\pi$
collapse $u$ to $\otp u$. Set $D=D^{\pi(\delta)}_{\pi[q^*],\pi[\name z]}$.
Fix
$v\includedin \alpha$, and $r'\in G_\alpha\intersect {\rm D}(u,v)$.
We can copy $r$ via an order-isomorphism 
inside $\alpha\times \aleph_2$, fixing $r\restrictedto \delta$, so
that the
result can be amalgamated with $r'$, to yield $r''$, which is then a
$(u,v)$-witness above $r'$.  Since $r'\in {\rm D}(u,v)$, this means
that $r'$ is itself a $(u,v)$-witness in $G_\alpha$. As this is all
that the construction in 3.12 was supposed to achieve, this case is
covered by the discussion there.
\medskip

\proclaim{3.7\\Minor Density Systems}
\endproclaim

In the course of the argument in 3.12, we require two further density
systems.  In the course of that argument we introduced the set 
$$S =\{\gamma\in S_C:\hbox{$\name F(\name x_\gamma)$ is a
$[\Pp\restrictedto (\gamma+1)]$-name}\},$$ 
and argued that $S$
is stationary.  
This led us to consider certain ordinals $\gamma<\delta$, with
$\delta$ of cofinality $\aleph_2$, and an
element $r_1\in G_\delta$, at
which point we claimed that we could produce 
a 1-1 order preserving function $h$ with domain $\AA^{r_1}$,
equal to the identity on $\AA^{r_1}\intersect (\gamma+1)$, with
$h(\min(\AA^{r_1}\backslash (\gamma+1)))>\sup \AA^{r_1}$, and
$h[r_1]\in G_\delta$.  More precisely, our claim was that this could be
ensured by meeting suitable density systems.

For $\alpha <\aleph_2$, $q\in \AP\restrictedto \aleph_2$, define
\def\Daq{{\rm D}^\alpha_q} $\Daq(u,v)$ as follows. If
$(\{\alpha\}\union\otp \dom q) \ge 
\otp u$ then let $\Daq(u,v)$ degenerate. Otherwise, fix
$k:(\{\alpha\}\union\dom q) 
\to u$ an order isomorphism onto an initial segment of
$u$, and let $\beta=\inf (u\setminus \rg k)$.  Let $\Daq(u,v)$ be the set of
$r\in \AP$ with domain $v$ such that $r\restrictedto v\setminus u$ contains the
image of $q$ under an order-preserving map $h_0$ which agrees with $k$ below
$\alpha$ and which carries $\inf(\AA^q\setminus (\alpha\times\aleph_2))$ above
$\beta$ (i.e., above $(\beta,0)$).   The density condition corresponds
to our ability to copy over part of $q$ onto any set of unused
ordinals in $(v\setminus \beta)\times \aleph_2$, recalling that $|\dom
r|<\aleph_2$ for 
any $r\in \AP$, and then to perform an amalgamation.

For our intended application, suppose that $\gamma, \delta, r_1$ are given
as above, and let $u=\{\gamma\}\union \dom r_1\union\{\sup\,\dom r_1)$.
Let $\pi$ be the canonical isomorphism of $u$
with $\otp u$, and $\alpha=\pi(\gamma)$, $q=\pi[r_1]$.  As $G_\delta$
meets $\Daq$, we 
have $v\includedin \delta$, and $r\in G_\delta\intersect \Daq(u,v)$. 
Then with $h=h_0\circ\pi$, we have $h[r_1]\le r$, and our claim is verified.

Finally, a few lines later in the course of the same argument we
mentioned that the claim proved in 3.14 can be construed as the
verification that certain additional density systems are in fact
dense, and that accordingly we may suppose that the condition $r$
described there lies in $G$.

\vfill\eject
{\bf References}
\bigskip
\frenchspacing
\newbox\papernumber
\def\ref#1.#2\par{\vskip .1 in\leavevmode\noindent
\setbox\papernumber=\hbox{[#1]\ }
\ifdim\wd\papernumber>2.8em\hangindent \wd\papernumber
\else\hangindent 2.8em\fi
\hangafter =1
\llap{\hbox to\hangindent{[#1]\hss}}#2\par}

\ref AxKo. J. Ax and S. Kochen, ``Diophantine problems over local
rings I,'' {\sl Amer. J. Math.} {\bf  87} (1965), 605-630.

\ref Bollobas. Bollobas, ``Random Graphs'', Academic Press,
London-New York 1985.

\ref JBSh 368. H. Judah, T.Bartoszynski, and S.Shelah, 
The Chico\'n diagram, {\sl Israel J. of Math.}

\ref Gregory.  J. Gregory, ``Higher Souslin trees and the generalized
continuum hypothesis,'' {\sl J. Symb. Logic} {\bf 41} (1976), 663-671.

\ref Jech. T. Jech, {\bf Multiple Forcing},  Cambridge Univ. Press,
Cambridge, 1986.

\ref Keisler. H. J. Keisler, ``Ultraproducts and saturated
models,'' {\sl Indag. Math.} {\bf 26} (1964), 178-186.

\ref Sh a. `` Classification theory and the number of non-isomorphic
models'', {\sl North Holland Publ. Co.},
Studies in Logic and the foundation of Math., vol. 92, 1978.

\ref Sh b. S. Shelah, {\bf Proper Forcing}, Lect. Notes Math. 940,
Springer-Verlag, Berlin Heidelberg NY Tokyo, 1982, 496 pp.

\ref Sh c.  Classification theory and the number of non-isomorphic models, 
revised, {\sl North Holland Publ. Co.} Studies in Logic 
and the foundation of Math., Vol.92, 1990, 705+xxxiv.

\ref Sh f. \same, {\bf Proper and Improper Forcing}, revised edition of [Sh
b], in preparation.

\ref Sh 13. \same, ``Every two elementarily equivalent models
have isomorphic ultrapowers,'' {\sl Israel J. Math} {\bf 10} (1971), 224-233.

\ref Sh 65. K. Devlin and S. Shelah, ``A weak form of the diamond
which follows from $2^{\aleph_0}<2^{\aleph_1}$,'' {\sl Israel J.
Math.} {\bf 29} (1978), 239-247.

\ref Sh 82. S. Shelah, ``Models with second order properties III.
Omitting types for $L(Q)$,'' {\sl in} {\bf Proceedings of a workshop
in Berlin, July 1977}, {\sl Archiv f. Math. Logik} {\bf 21} (1981), 1-11.

\ref Sh 107. \same, ``Models with second order properties IV. A
general method and eliminating diamonds,'' {\sl Annals Math. Logic}
{\bf 25} (1983), 183-212.

\ref ShHL 162. S. Shelah, B.  Hart, C.  Laflamme,
 ``Models with second order
properties V. A general principle,'' {\sl Annals Pure Appl. Logic}
{\bf}, to appear.

\ref Sh 177. S. Shelah, ``More on proper forcing,'' {\sl J. Symb. Logic} {\bf 49}
(1984), 1035-1038.

\ref Sh 345. \same, ``Products of regular cardinals and cardinal invariant
of Boolean Algebras'', {\sl Israel J. of Math} {\bf 70} (1990), 129-187. 

\ref Sh 405. \same, ``Vive la differance II-refuting  Kim's conjecture'', 
  Oct 89 (preprint).

\ref ShFr 406. S.Shelah and D.Fremlin,
``Pointwise compact and stable sets of measurable functions'',
Journal of Symbolic Logic, submitted.
\bye